\documentclass[11pt]{article}

\usepackage{amsmath}
\usepackage{amssymb} 
\usepackage{caption}
\usepackage{graphicx}
\usepackage{authblk}
\usepackage[hypertexnames=false,colorlinks=true,linkcolor=blue,citecolor=blue]{hyperref}
\usepackage[numbers,comma,square,sort&compress]{natbib}
\usepackage[a4paper,text={6.5in,10in},centering]{geometry}

\graphicspath{{eps/}{pdf/}{images/}}

\makeatletter\@addtoreset{equation}{section}\makeatother

 {\begin{trivlist} \item[]{\bf Proof. }}%
 {\hspace*{\fill}$\rule{.4\baselineskip}{.4\baselineskip}$\end{trivlist}}

 {\begin{trivlist}\item[]\textbf{Acknowledgments.}}{\end{trivlist}}
{\begin{trivlist}\item[]\textbf{Proof#1 }}%
 {\hspace*{\fill}$\rule{0.3\baselineskip}{0.35\baselineskip}$\end{trivlist}}

\usepackage{amsmath}
\usepackage{caption}
\usepackage{subcaption}
\usepackage{wrapfig}

\def\xb{\mathbf x}

\def\vb{\mathbf v}

\begin{document}

\title{Tracking and forecasting oscillatory data streams using Koopman autoencoders and Kalman filtering}
\author[1]{Stephen A Falconer}
\author[1]{David J.B. Lloyd}
\author[1]{Naratip Santitissadeekorn}
\affil[1]{\small Department of Mathematics, University of Surrey, Guildford, GU2 7XH, UK}
\date{\today}
\maketitle


\begin{abstract}
Data-driven modelling techniques provide a method for deriving models of dynamical systems directly from complicated data streams. However, tracking and forecasting such data streams poses a significant challenge to most methods, as they assume the underlying process and model does not change over time. In this paper, we apply one such data-driven method, the Koopman autoencoder (KAE), to high-dimensional oscillatory data
to generate a low-dimensional latent space and model, where the system's dynamics appear linear. This allows one to accurately track and forecast systems where the underlying model may change over time. States and the model in the reduced order latent space can then be efficiently updated as new data becomes available, using data assimilation techniques such as the ensemble Kalman filter (EnKF), in a technique we call the KAE EnKF. We demonstrate that this approach is able to effectively track and forecast time-varying, nonlinear dynamical systems in synthetic examples. We then apply the KAE EnKF to a video of a physical pendulum, and achieve a significant improvement over current state-of-the-art methods. By generating effective latent space reconstructions, we find that we are able to construct accurate short-term forecasts and efficient adaptations to externally forced changes to the pendulum's frequency.
\end{abstract}

\section{Introduction}

In the modern era, most systems of interest are too complex to be amenable to traditional derivations from first principles \cite{data-driven_science_and_engineering_preface}. To deal with this, data-driven modelling methods have been devised that allow for the derivation of dynamical systems models from complex data \cite{data-driven_modelling_science_and_engineering}, having been applied for the purpose of modelling, forecasting, and controlling dynamical systems \cite{data-driven_science_and_engineering_preface}, and with data-driven discovery being referred to as the fourth paradigm of science \cite{fourth_paradigm_data_science}. In fields such as epidemiology \cite{data-driven_epidemiology}, neuroscience \cite{data-driven_neuroscience} and finance \cite{data-driven_finance} where governing equations are partially or completely unknown, data-driven discovery provides a much-needed window of insight into how these systems operate. Typically, large-amounts of data can arrive sequentially (known as {\it data streams}) that need to be handled from a data-driven discovery perspective, as well as wanting to track and forecast in real-time.

Dynamic mode decomposition (DMD) is one of the most currently well-known data-driven modelling techniques, identifying a low dimensional, linear model of a dynamical system directly from high-dimensional data \cite{official_intro}. DMD has been successfully applied in a multitude of fields \cite{epidemiology_use, finance_use, neuroscience_use}, and due to its popularity, has been extended for use cases with multi-resolution \cite{mrdmd}, noisy \cite{ncdmd, tdmd, kfdmd} or streaming \cite{streaming_dmd, incremental_tdmd, online_dmd, dmdenkf} data.

A major limitation of these DMD based frameworks is that they produce linear models, whereas most complex systems that exhibit varied and interesting behaviour are in fact nonlinear \cite{nonlinear_systems_book}. In this paper, we are concerned with learning, tracking, and forecasting discrete, nonlinear dynamical systems given by
\begin{equation}
    \mathbf{x}_{k+1} = \mathbf{f}(\mathbf{x}_k),
    \label{eqn:xk+1=fxk}
\end{equation}
where $\mathbf{x}_k \in \mathbb{X}$ represents the system's state in the state space $\mathbb{X} \subseteq \mathbb{R}^n$ at time $k$. 
Koopman theory states that the nonlinear system in equation \eqref{eqn:xk+1=fxk}, can instead be represented by a linear operator $\mathbf{K}:\mathbb L^{2}(\mathbb{X}) \rightarrow \mathbb L^{2}(\mathbb{X})$ generated by $\mathbf{f}$ known as the Koopman operator. While the Koopman operator is linear, it comes with the caveat that it is infinite dimensional, propagating the set of measurements of a state $\mathbf{g}(\mathbf{x}_k)\in\mathbb{L}^2(\mathbb{X})$ forward in time \cite{koopman_og} via
\begin{equation}
    \mathbf{g}(\mathbf{x}_{k+1}) = \mathbf{K}\mathbf{g}(\mathbf{x}_k).
    \label{eqn:gxk+1=Kgxk}
\end{equation}
A nonlinear system's behaviour is completely defined by the spectral decomposition of its Koopman operator, hence to tractably model a system from a Koopman theoretic perspective, important eigenfunctions of the Koopman operator must be identified \cite{modern_koopman_theory}.

Methods for finding an approximation of the Koopman operator are often achieved by augmenting the state measurements in some way, for example DMD can be applied to a nonlinear system after altering the data via a nonlinear transformation in extended DMD \cite{extended_dmd}, including time-delay embeddings in Hankel-DMD \cite{hankel-dmd}, or using kernel functions in kernel DMD \cite{kernel_dmd}. Other methods of augmenting the state measurements include the use of Koopman autoencoders (KAE) \cite{dmd_autoencoder, autoencoder_dmd, consistent_koopman_autoencoder}, grading representations using the VAMP-score \cite{VAMP-score, VAMPnet}, adding an intermittent nonlinear forcing term \cite{havok}, or some combination of the above techniques \cite{deep_delay_autoencoders}.

When modelling high-dimensional data, the computational efficiency of extended DMD scales poorly \cite{dmd_book_nonlinear}, and the results produced by kernel DMD are highly sensitive to the kernel chosen by the user \cite{kernel_dmd}. As such, the preferred methods for identifying Koopman eigenfunctions of nonlinear systems often employ neural networks in some part of their architecture, embracing their ease of use, effectiveness \cite{deep_learning_book}, and expressivity as universal function approximators \cite{neural_net_universal_approximator_theory}. These models however can be difficult to train, as when using stochastic gradient descent to optimize Koopman eigenvalues across multiple forecast horizons, the optimization landscape becomes pitted with local minima \cite{henning_lange_fourier_to_koopman}, which can lead to the production of defective models \cite{local_minima_neural_net_cases}. Neural network based models, particularly the larger, more powerful models, also take significant amounts of time to train \cite{deep_learning_book}. Hence, there is a growing demand for models that can process data in real-time \cite{machine_learning_real-time}, and to maintain the most up-to-date models, they must incorporate relevant information from the new data as it is generated.

Data assimilation refers to the collection of methods that integrate new data with sophisticated mathematical models, to track and forecast systems that may evolve or change \cite{data_assimilation_book}. The majority of its applications lie in the earth sciences \cite{da_use_earth_sciences}, however due to the generality of its techniques they have also been successfully applied in a wide range of areas from medicine \cite{da_use_medicine} to ecology \cite{da_use_ecology}. The Kalman filter \cite{kf_og} is one such data assimilation technique widely used throughout industry \cite{kf_applications}, that optimally combines predictions from a linear model with Gaussian data. Whilst traditionally applied to a model's state, the parameters of the model can simultaneously be filtered, leading to what is known as the joint state-parameter estimation problem \cite{kf_state_param}. When the system being filtered is nonlinear, as in our case, alternative versions of the Kalman filter can be utilized such as the extended Kalman filter \cite{ekf}, unscented Kalman filter \cite{ukf} or ensemble Kalman filter (EnKF) \cite{enkf}. The EnKF represents the distribution of a system's state with an ensemble of random samples, that can then be used to estimate useful statistics like the state's covariance via the sample covariance or a point estimate of the state via the sample mean \cite{enkf}, and is well-suited for high-dimensional problems. By creating a model of the system using data-driven modelling techniques such as Koopman autoencoders, for which we filter the states and parameters as new data becomes available using data assimilation techniques like the EnKF, we are able to develop a generalizable framework for building and maintaining up-to-date models of dynamical systems, directly from data.

In this paper, we first propose a series of alterations to the Koopman autoencoder developed in \cite{dmd_autoencoder} for high-dimensional oscillatory data streams, with the aim of encouraging efficient training using the layer initialization technique from \cite{deep_delay_autoencoders}, identifying globally optimal eigenvalues via the optimization procedure in \cite{henning_lange_fourier_to_koopman}, and producing stable models through the inclusion of an additional loss term. We then detail how the Koopman autoencoder and EnKF can be combined, taking previous work that introduced the DMDEnKF \cite{dmdenkf}, and expanding it to nonlinear systems, in what we refer to as the KAE EnKF.

The paper is outlined as follows. In Section~\ref{s:Koopman_auto_filter}, we review Koopman autoencoders and the ensemble Kalman filter, followed by Section~\ref{s:KAEEnKF} detailing the KAE EnKF algorithm. In Section~\ref{s:synthetic}, we then demonstrate the KAE EnKF on synthetically generated datasets, analysing its effectiveness relative to its linear counterpart, the DMDEnKF, at tracking a system's eigenvalues for varying degrees of nonlinearity in the system. We then compare the KAE EnKF's efficacy when tracking and forecasting a time-varying, nonlinear system, with that of other pre-existing iterative, nonlinear, data-driven modelling techniques. The decision in the KAE EnKF framework to filter the Koopman autoencoder's latent states as opposed to the full state measurements is then thoroughly investigated, before the KAE EnKF is applied to a nonlinear system with multiple frequencies, and again its performance is evaluated against that of other contemporary iterative, nonlinear techniques. In Section~\ref{s:pendulum}, we then demonstrate the KAE EnKF applied to raw video footage of a pendulum in motion. We analyse the Koopman autoencoder's ability to generate informative latent state representations, and subsequent full state reconstructions of the data, compared with those produced by extended DMD. We then assess the KAE EnKF's forecasting skill on unseen data from the pendulum video, afterwards doubling the pendulum's frequency, and investigating the KAE EnKF's ability to efficiently adapt to the change in the system. Finally, in Section~\ref{s:conclusion}, we draw conclusions and discuss future directions.

\section{Koopman autoencoders and ensemble Kalman filters}\label{s:Koopman_auto_filter}

In this section, we briefly review Koopman autoencoders and ensemble Kalman filters that will be used for the KAE EnKF method. 

\subsection{Koopman Autoencoders}
\begin{figure}[htbp]
\begin{subfigure}[b]{\textwidth}
         \centering
         \includegraphics[width=\textwidth]{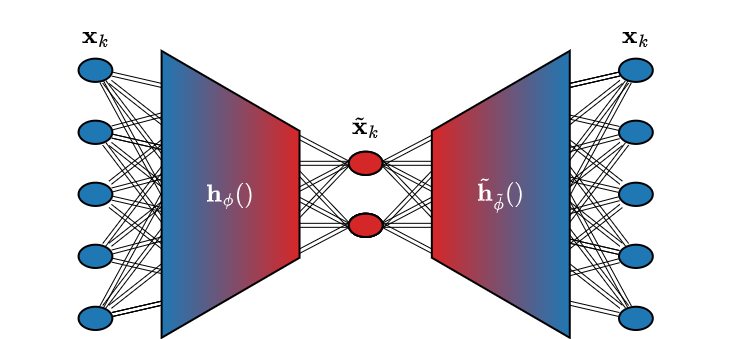}
     \end{subfigure}
     \caption{\centering\label{fig:litautoencoder}States $\mathbf{x}_k$ are encoded via neural network $\mathbf{h}_\phi()$ into latent states $\Tilde{\mathbf{x}}_k$, before being decoded back into the original state by the network $\mathbf{\Tilde{h}}_{\Tilde{\phi}}()$. By constraining the latent space in some way, then optimizing the trainable parameters $\phi$ and $\Tilde{\phi}$ of the encoder/decoder over the network's reconstruction loss, a useful latent representation of the data can be learnt by the autoencoder.}
\end{figure}

Autoencoders are a flexible architecture for encoding/decoding data into a latent space \cite{autoencoder_og}. Standard autoencoders operate on states $\mathbf{x}_k \in{\mathbb R^{n}}$ with $k=1,\hdots ,m$, and encode/decode them into the latent space $\mathbf{\Tilde{x}}_k \in{\mathbb R^{r}}$, where traditionally $r \ll n$, via the transformations
\begin{equation}
    \mathbf{h}_\phi(\mathbf{x}_k) = \mathbf{\Tilde{x}}_k, \quad \mathbf{\Tilde{h}}_{\Tilde{\phi}}(\mathbf{\Tilde{x}}_k) = \mathbf{x}_k.
\label{eqn:gxk=tildexk}
\end{equation}
Functions $\mathbf{h}_\phi$ and $\mathbf{\Tilde{h}}_{\Tilde{\phi}}$, the autoencoder's encoder and decoder functions, are usually neural networks \cite{neural_net_autoencoder}, with trainable parameters $\phi$ and $\Tilde{\phi}$ respectively, representing the networks weights and biases. A graphical representation of an autoencoder can be seen in Figure \ref{fig:litautoencoder}

Parameters $\phi$ and $\Tilde{\phi}$ are optimized using stochastic gradient descent based algorithms, and the overall loss function $\mathcal{L}$ to be minimized is
\begin{equation}
    \mathcal{L}(\mathbf{x}_k) = \mathbf{D}(\mathbf{\Tilde{h}}_{\Tilde{\phi}}(\mathbf{h}_\phi(\mathbf{x}_k)), \mathbf{x}_k) + \mathbf{R}(\phi, \Tilde{\phi}),
\label{eqn:autoencoder_loss}
\end{equation}
where $\mathbf{D}$ is some measure of difference between its two inputs $\mathbf{\Tilde{h}}_{\Tilde{\phi}}(\mathbf{h}_\phi(\mathbf{x}_k))$ and $\mathbf{x}_k$, and $\mathbf{R}(\phi, \Tilde{\phi})$ is a regularising term over the encoder/decoder network's parameters. The first term in equation \eqref{eqn:autoencoder_loss} is known as the reconstruction error, and encapsulates the loss of information caused by encoding and subsequently decoding the input data.

\begin{figure}[htbp]
\begin{subfigure}[b]{\textwidth}
         \centering
         \includegraphics[width=\textwidth]{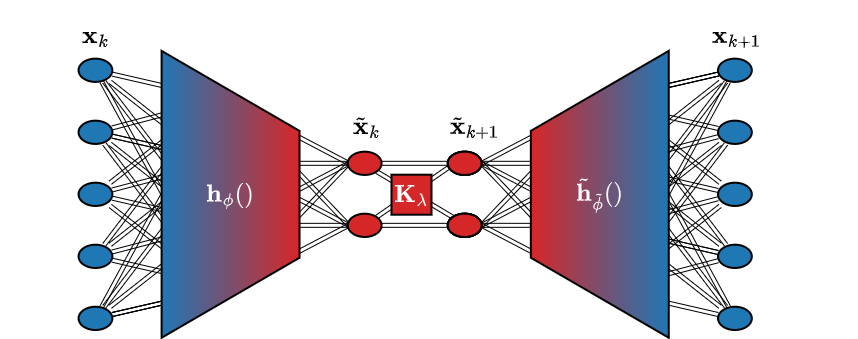}
     \end{subfigure}
     \caption{\centering\label{fig:kae}As in the autoencoder shown in Figure \ref{fig:litautoencoder}, functions $\mathbf{\Tilde{h}}_{\Tilde{\phi}}$ and $\mathbf{h}_\phi$ encode/decode states $\mathbf{x}_k$ to/from latent states $\Tilde{\mathbf{x}}_k$. Linear operator $\mathbf{K}_\lambda$ is applied in the latent space, to propagate the latent state $\Tilde{\mathbf{x}}_k$ one time step forward into $\Tilde{\mathbf{x}}_{k+1}$. By optimizing the encoder/Koopman approximator/decoder's trainable parameters $\phi$, $\lambda$, and $\Tilde{\phi}$, over the standard autoencoder's loss function with additional linearity/forecasting terms, the Koopman autoencoder is able to learn a latent representation of the data in which the system's dynamics act linearly.}
\end{figure}

Koopman autoencoders \cite{dmd_autoencoder, autoencoder_dmd} differ from the basic autoencoder setup in a number of ways. Firstly, the states $\mathbf{x}_k$ with $k=1,\hdots ,m$ are assumed to be measurements of a discrete, possibly nonlinear, dynamical system as generated by equation \eqref{eqn:xk+1=fxk}.
The aim is to find an encoding $\mathbf{h}_\phi(\mathbf{x}_k) = \mathbf{\Tilde{x}}_k$, that represents the important Koopman eigenfunctions of the system, usually for the purpose of system identification and forecasting future states. Hence, when encoding/decoding the data as in equation \eqref{eqn:gxk=tildexk}, we operate under the additional constraint that in the latent space $\mathbf{\Tilde{x}}_k \in{\mathbb R^{r}}$, the dynamics act linearly with
\begin{equation}
    \mathbf{\Tilde{x}}_{k+1} = \mathbf{K}_\lambda\mathbf{\Tilde{x}}_k.
    \label{eqn:xtk+1=klxtk}
\end{equation}
The matrix $\mathbf{K}_\lambda \in{\mathbb R^{r\times r}}$ is a finite dimensional approximation of the Koopman operator in the space of its most important eigenfunctions and corresponding eigenvalues, $\lambda$. 
The resultant structure of a Koopman autoencoder is shown in Figure \ref{fig:kae}.

As $\mathbf{K}_\lambda$ is represented in its eigenbasis, it could be denoted as a complex, diagonal matrix with eigenvalues $\lambda$ along its lead diagonal. For convenience, eigenvalues $\lambda$ are instead expressed in polar coordinates $\lambda_i = \tau_i e^{\theta_i}$, resulting in an effective parameterization of $\mathbf{K}_\lambda$ by
\begin{equation}
    \lambda=
    \left[\begin{array}{c c c c c c}
            \tau_1 & \hdots & \tau_r & \theta_1 & \hdots & \theta_r
             \end{array}\right]^T
             \label{eqn:lambda=tauthetas}.
\end{equation}
The corresponding eigenvalue modulus $\tau_i$ and arguments $\theta_i$ from $\lambda \in{\mathbb R^{2r}}$ are then stacked into their natural rotation matrices $\mathbf{B}_i \in{\mathbb R^{2\times 2}}$. Presented in this way, the matrix $\mathbf{K}_\lambda$ takes the real, block diagonal form
\begin{equation}
    \mathbf{B}_i = \left[\begin{array}{cc}
         \tau_i\cos(\theta_i)&-\tau_i\sin(\theta_i)  \\
         \tau_i\sin(\theta_i)&\tau_i\cos(\theta_i) 
    \end{array}\right], \quad
    \mathbf{K}_\lambda = \left[\begin{array}{cccc}
         \mathbf{B}_1& &\mathbf{0} \\
         &\ddots& \\
         \mathbf{0}& &\mathbf{B}_r 
    \end{array}\right].
    \label{eqn:kl=polar_eig_blocks}
\end{equation}

The additional requirements of the Koopman autoencoder to produce a latent space with linear dynamics, and be useful for future state forecasting, each create an extra term in the standard autoencoders loss function
\begin{equation}
\begin{split}
    \mathcal{L}(\mathbf{x}_k,\mathbf{x}_{k+1}) = &\|\mathbf{x}_k-\mathbf{\Tilde{h}}_{\Tilde{\phi}}(\mathbf{h}_\phi(\mathbf{x}_k)) \|_2^2 + a_1\|\mathbf{h}_\phi(\mathbf{x}_{k+1})- \mathbf{K}_\lambda\mathbf{h}_\phi(\mathbf{x}_k) \|_2^2 \\
    &+ a_2\|\mathbf{x}_{k+1}-\mathbf{\Tilde{h}}_{\Tilde{\phi}}(\mathbf{K}_\lambda\mathbf{h}_\phi(\mathbf{x}_k)) \|_2^2 + a_3\mathbf{R}(\phi, \Tilde{\phi}).
\end{split}
\label{eqn:kae_loss}
\end{equation}
The first and last loss terms remain relatively unchanged from equation \eqref{eqn:autoencoder_loss}, with the reconstruction and regularization errors fulfilling the same roles as before. The only caveat is that as dynamical systems data is unlikely to lie exclusively in the range of 0 and 1, we replace the general difference measure $\mathbf{D}$ with the standard mean squared error. The second term in equation \eqref{eqn:kae_loss} represents the Koopman autoencoders latent linear dynamics condition, penalizing how far the latent space and operator $\mathbf{K}_\lambda$ are from fulfilling equation \eqref{eqn:xtk+1=klxtk}. The third loss term is the forecasting error, which quantifies the model's ability to produce forecasts in the full state space from its linear $\mathbf{K}_\lambda$ latent space forecasts. Computational parameters $a_i$ are used to adjust the influence of each error term on the overall loss function, and are calibrated using hyperparameter optimization techniques, then kept constant as before.

To encourage more stable forecasts as in \cite{dmd_autoencoder}, the two additional loss terms in equation \eqref{eqn:kae_loss} can be modified to penalize inaccuracies when generating forecasts of up to $p$ time steps ahead, as opposed to only one step ahead. This is done by simply applying the approximation of the Koopman operator $\mathbf{K}_\lambda$ to the latent state $\mathbf{\Tilde{x}}_{k}$ $p$ times, as from equation \eqref{eqn:xtk+1=klxtk} we have
\begin{equation}
    \mathbf{K}_\lambda^p\mathbf{\Tilde{x}}_{k} = \mathbf{K}_\lambda^{p-1}\mathbf{\Tilde{x}}_{k+1} = \mathbf{\Tilde{x}}_{k+p}.
\end{equation}
The latent state $\mathbf{\Tilde{x}}_{k+p}$ can be used directly to evaluate the $p$ step ahead linearity error of the dynamics via 
\begin{equation}
\|\mathbf{h}_\phi(\mathbf{x}_{k+p})- \mathbf{K}_\lambda^p\mathbf{h}_\phi(\mathbf{x}_k) \|_2^2.
\end{equation}
It can then be decoded back into the original state space by applying $\mathbf{\Tilde{h}}_{\Tilde{\phi}}(\mathbf{\Tilde{x}}_{k+p}) = \mathbf{x}_{k+p}$ from equation \eqref{eqn:gxk=tildexk}, to assess the full state forecasting error given by
\begin{equation}
\|\mathbf{x}_{k+p}- \mathbf{\Tilde{h}}_{\Tilde{\phi}}(\mathbf{K}_\lambda^p\mathbf{h}_\phi(\mathbf{x}_k)) \|_2^2.
\end{equation}
The error terms are formed to balance the losses over a range of forecast horizons by calculating errors for all forecasts up to a horizon $p$, then using their average in the loss calculation. Applying this method for $dt \in{(1, \hdots, p )}$ produces the loss function for the Koopman autoencoder that we use given by
\begin{equation}
\begin{split}
    \mathcal{L}(\mathbf{x}_k,\hdots,\mathbf{x}_{k+p}) = &\|\mathbf{x}_k-\mathbf{\Tilde{h}}_{\Tilde{\phi}}(\mathbf{h}_\phi(\mathbf{x}_k)) \|_2^2 + 1/p\sum_{dt=1}^{p}(a_1\|\mathbf{h}_\phi(\mathbf{x}_{k+dt})- \mathbf{K}_\lambda^{dt}\mathbf{h}_\phi(\mathbf{x}_k) \|_2^2 \\
    &+ a_2\|\mathbf{x}_{k+dt}-\mathbf{\Tilde{h}}_{\Tilde{\phi}}(\mathbf{K}_\lambda^{dt}\mathbf{h}_\phi(\mathbf{x}_k)) \|_2^2) + a_3\mathbf{R}(\phi, \Tilde{\phi}).
    \end{split}
\label{eqn:kae_loss2}
\end{equation}
Calculating powers of $\mathbf{K}_\lambda^{dt}$ is computationally cheap, due to the block diagonal structure of the matrix given by
\begin{equation}
        \mathbf{K}_\lambda^{dt} = \left[\begin{array}{cccc}
         \mathbf{B}_1^{dt}& &\mathbf{0} \\
         &\ddots& \\
         \mathbf{0}& &\mathbf{B}_r^{dt}
    \end{array}\right],\qquad\mbox{where}\qquad \mathbf{B}_i^{dt} = \left[\begin{array}{cc}
         \tau_i^{dt}\cos(\theta_i dt)&-\tau_i^{dt}\sin(\theta_i dt)  \\
         \tau_i^{dt}\sin(\theta_i dt)&\tau_i^{dt}\cos(\theta_i dt) 
    \end{array}\right].
\end{equation}
Once training is complete, the Koopman autoencoder can produce forecasts $dt$ steps ahead from state $\mathbf{x}_k$ using
\begin{equation}
\mathbf{x}_{k+dt} = \mathbf{\Tilde{h}}_{\Tilde{\phi}}(\mathbf{K}_\lambda^{dt}\mathbf{h}_\phi(\mathbf{x}_k)).
\end{equation}

Full implementation details for how the Koopman autoencoder network is built and trained on the example applications are provided in the appendix.

\subsection{Ensemble Kalman Filters (EnKF)}

We consider the discrete-time, nonlinear dynamical system with a stochastic perturbation given by
\begin{equation}
    \mathbf{x}_{k} = \mathbf{F}(\mathbf{x}_{k-1}) + \mathbf{w}_k, \quad
    \mathbf{w}_k \sim \mathcal{N}(\mathbf{0},\mathbf{Q}_k),
    \label{eqn:kf_prop_eq}
\end{equation}
where $\mathbf{F}$ is a nonlinear function $\mathbf{F}:\mathbb{R}^{n} \to \mathbb{R}^{n}$, $\xb_k \in \mathbb{R}^{n}$ is the (observable) state variable, and the random variable $\mathbf{w}_k$ is normally distributed with mean $\mathbf{0}$ and covariance matrix $\mathbf{Q}_k$. 
An observation equation that relates the observation $\mathbf{y_k} \in \mathbb{R}^{l}$ to the unobserved state of the system is given by
\begin{equation}
    \mathbf{y}_k = \mathbf{H}(\mathbf{x}_k) + \mathbf{v}_k, \quad
    \mathbf{v}_k \sim \mathcal{N}(\mathbf{0},\mathbf{R}_k),
    \label{eqn:kf_update_eq}
\end{equation}
where $\mathbf{H}:\mathbb{R}^{n} \to \mathbb{R}^{l}$ is the system's observation operator, where the observational noise $\mathbf{v_k} \in \mathbb{R}^{l}$ has a normal distribution with mean $\mathbf{0}$ and covariance matrix $\mathbf{R}_k$. For our current work, $\mathbf{H}$ is linear, so can be represented by a matrix $\mathbf{H} \in \mathbb{R}^{l\times n}$.

The main objective of filtering is to use information of the observation up to the time step $k$ to make a statistical inference about the underlying state at time $k$. In particular, it aims to obtain the marginal posterior density $p(\mathbf{x_k} | \mathbf{Y_k})$ of the current state $\mathbf{x_k}$ give the observations up to time step $k$, where $\mathbf{Y_k}=(\mathbf{y_1}, ..., \mathbf{y_k}$). When $\mathbf{F}$ and $\mathbf{H}$ are both linear and both $\mathbf{w}_k$ and $\mathbf{v}_k$ are normal, the Kalman filter provides the exact Bayesian formulation for  $p(\mathbf{x_k} | \mathbf{Y_k})$ \cite{kf_og}. It can also be used to find the best linear unbiased estimator (BLUE) for linear non-normal systems. 

In many applications, such as those in geophysics and fluid dynamics, either the linear or normal assumption (or both) is violated and the dimension of the underlying state variable can be very large. The Ensemble Kalman filter (EnKF) is an approximation method to address these issues \cite{enkf} by using a set of samples (called ensembles) to represent the filtering distribution $p(\mathbf{x_k} | \mathbf{Y_k})$. In EnKF, the ensemble of $p(\mathbf{x_k} | \mathbf{Y_k})$ is updated by mimicking the Kalman filter formulation in such a way that the ensemble approximation would be consistent with the linear and normal system as the number of sample goes to infinity.
Any statistical inference such as the mean of $\mathbf{x_k}$ can then be achieved using the ensemble. For a high-dimensional problem, the number of samples is typically much smaller than the dimension of the state variable. However, the focus of our application here is not the high-dimensional issue but the non-linearity of $\mathbf{F}$. 

To begin the EnKF algorithm, an initial ensemble of $N$ state estimates $\hat{\xb}^{(1)}_0\text{,..., } \hat{\xb}^{(N)}_0$ is required. If an ensemble is not available, one can be generated from initial state estimates  $\hat{\xb}_0$ and covariance matrix $\mathbf{P}_0$ by taking $N$ independent draws from $\mathcal{N}(\hat{\xb}_0,\mathbf{P}_0)$. 

\paragraph{Algorithm}

The so-called stochastic EnKF algorithm is recursively carried out as below \cite{enkf}:

    \textbf{Step 1 (Forecast):} Propagate forward in time each ensemble member using the equation \eqref{eqn:kf_prop_eq} for $i = 1,...,N$ via
    \begin{equation}
         \hat{\xb}^{(i)}_{k|k-1} = \mathbf{F}(\hat{\xb}^{(i)}_{k-1|k-1}) + \mathbf{w}^{(i)}_k.
            \label{eqn:kf_k|k-1_eqs}
    \end{equation}
    The notation $\hat{\xb}^{(i)}_{k|k-1}$ denotes the state estimate at time $k$ of the $i$th ensemble member $\hat{\xb}^{(i)}_k$ using only information up to time $k-1$, and $\hat{\xb}^{(i)}_{k-1|k-1}$ represents the same ensemble member at time $k-1$ using information up to time $k-1$. Each $\mathbf{w}^{(i)}_k$ is independently drawn from $\mathcal{N}(\mathbf{0},\mathbf{Q}_k)$. The covariance matrix of $\hat{\xb}_{k|k-1} $ is estimated by the sample covariance, denoted by $\mathbf{\hat{P}}_{k|k-1}$. Note that $\mathbf{\hat{P}}_{k|k-1}$ will not be full-rank if $N<n$, which will create some bias. Some techniques such as covariance inflation can be used to ameliorate this issue. However, this is not a serious concern in our current application, since we will choose $n>N$. 
    An approximation of the Kalman Gain matrix $\mathbf{\hat{K}}_k$ can then be obtained by
    \begin{equation}
        \mathbf{\hat{K}}_k = \mathbf{\hat{P}}_{k|k-1}\mathbf{H}^T(\mathbf{H}\mathbf{\hat{P}}_{k|k-1}\mathbf{H}^T + \mathbf{R}_k)^{-1}.
        \label{eqn:kf_K}
    \end{equation}

    \textbf{Step 2 (Update):} The ensemble is updated when the observation $\mathbf{y}_k$ becomes available using the Kalman update formulation
        \begin{equation}
        \hat{\xb}^{(i)}_{k|k} = \hat{\xb}^{(i)}_{k|k-1} + \mathbf{\hat{K}}_k\mathbf{e}^{(i)}_k,
        \end{equation}
    where $\mathbf{e}^{(i)}_k  = \mathbf{y}_k + \mathbf{v}^{(i)}_k - \mathbf{H} \hat{\xb}^{(i)}_{k|k-1}$. Here each $\mathbf{v}^{(i)}_k$ is independently drawn from $\mathcal{N}(\mathbf{0},\mathbf{R}_k)$. It is this addition of $\mathbf{v}^{(i)}_k$ to the actual observation $\mathbf{y}_k$ that vitally distinguishes the stochastic EnKF from other deterministic EnKF algorithms, see \cite{Tippett2003,Bishop2001}.

\section{KAE EnKF}\label{s:KAEEnKF}

The Koopman autoencoder ensemble Kalman filter (KAE EnKF) builds on the Koopman autoencoder framework, by utilizing an ensemble Kalman filter on its latent states $\Tilde{\mathbf{x}}_k$ and parameters $\lambda$. This provides the Koopman autoencoder with the capability to efficiently assimilate new data, allowing for further refinement of the model or improved tracking of a system's states and non-stationary parameters as fresh data becomes available. Uncertainty estimates for the system's states and parameters are also produced, which are critical for making informed decisions based on a model's output \cite{uncertainty_quantification}.

The problem setup begins with a data stream consisting of states $\mathbf{x}_k \in{\mathbb R^{n}}$, measured at regular time intervals $k = 1, \hdots, m$, and then measured iteratively at $k= m+1, \hdots$. These states are all assumed to be generated by the system
\begin{equation}
    \mathbf{x}_{k+1} = \mathbf{f}_k(\mathbf{x}_k).
\end{equation}
This setup is similar to that of the Koopman autoencoder, with the only differences being the addition of a possible time dependence on evolution rule $\mathbf{f}_k$, unlike in equation \eqref{eqn:xk+1=fxk}, and the data's subsequent iterative generation. The KAE EnKF framework makes use of a Koopman autoencoder, which is trained during what is referred to as the algorithm's spin-up phase,  however the training procedure deviates in some ways from that of the original Koopman autoencoder described above. Hence, we first state these changes in what is referred to as the altered Koopman autoencoder, before detailing the KAE EnKF algorithm in full.

\subsection{Adapted Koopman Autoencoder}
The changes made to the original Koopman autoencoder's training procedure, that we then use in the KAE EnKF framework, are as follows.

\textbf{Layer Initialization: }
We design the autoencoder's neural network architecture such that the first layer of the encoder and last layer of the decoder are both linear layers. The weights of these layers are initialized using the truncated SVD and its transpose, calculated over matrix $\mathbf{X} \in{\mathbb R^{n\times m}}$, which is formed by stacking the spin-up data columnwise. This method of initializing the first and last layers of the network with the best linear approximations using the truncated SVD is taken directly from \cite{deep_delay_autoencoders}. All other weights and biases within the network are then initialized using the default Xavier initialization scheme \cite{xavier_init}.

\textbf{Sampling Scheme: }
The dimension $n$ of measured states $\mathbf{x}_k$ is assumed to be much larger than that of the latent space $r$, which we are able to effectively represent the system in $n \gg r$. For this reason, often only a small subset of the available measurements $\mathbf{x}_k$ are able to be loaded into memory at any one time. By averaging the loss for each state up to a forecast horizon $p$ as in equation \eqref{eqn:kae_loss2}, the loss calculation over one data point requires loading $p + 1$ states of dimension $n$ into memory. To mitigate this memory bottleneck, we instead calculate losses on a pairwise basis only, drawing an input $\mathbf{x}_k$ at random from the training set, and then selecting a forecast horizon from $dt \in{(1, \hdots, p)}$ also at random to determine the target output $\mathbf{x}_{k + dt}$. By repeating this procedure to generate a batch of appropriate size for the current machine's memory, we are able to simultaneously balance the optimization of Koopman approximator $\mathbf{K}_\lambda$ over multiple forecast horizons up to $p$, without relying too heavily on any one input $\mathbf{x}_k$ at any optimization step.

\textbf{Loss Function: }
The loss function from equation \eqref{eqn:kae_loss2} to be minimized over each input/target state pair in the altered Koopman autoencoder's batch thus becomes
\begin{equation}
\begin{split}
    \mathcal{L}(\mathbf{x}_k,\mathbf{x}_{k+dt},dt) =&
    \|\mathbf{x}_k-\mathbf{\Tilde{h}}_{\Tilde{\phi}}(\mathbf{h}_\phi(\mathbf{x}_k)) \|_2^2 + 
    a_1\|\mathbf{h}_\phi(\mathbf{x}_{k+dt})- \mathbf{K}_\lambda^{dt}\mathbf{h}_\phi(\mathbf{x}_k) \|_2^2 \\
    &+ a_2\|\mathbf{x}_{k+dt}-\mathbf{\Tilde{h}}_{\Tilde{\phi}}(\mathbf{K}_\lambda^{dt}\mathbf{h}_\phi(\mathbf{x}_k)) \|_2^2 
    + a_3\sum_{i=1}^r(\tau_i - 1) \\
    &+ a_4 (\| \phi \|_2^2 + \| \Tilde{\phi} \|_2^2 + \| \phi \|_1 + \| \Tilde{\phi} \|_1).
    \end{split}
\label{eqn:kae_enkf_loss2}
\end{equation}
We include an additional stability loss, that biases the Koopman approximator towards stable frequencies, by penalizing the $L^1$ norm of the deviation of $\mathbf{K}_\lambda$'s eigenvalues from the unit circle $\sum_{i=1}^r(\tau_i - 1)$. The regularization term we use for the neural network parameters $\phi$ and $\Tilde{\phi}$, is a combination of their $L^1$ and $L^2$ norms, often referred to in machine learning literature as elastic net. This form of regularization was chosen as the $L^1$ term encourages sparse solutions, while the $L^2$ term stops the regularization loss from saturating, in the event that the state size $n$ exceeds the number of available data points \cite{elastic_net}. 

\begin{figure}[htbp]
\begin{subfigure}[b]{\textwidth}
         \centering
         \includegraphics[width=\textwidth]{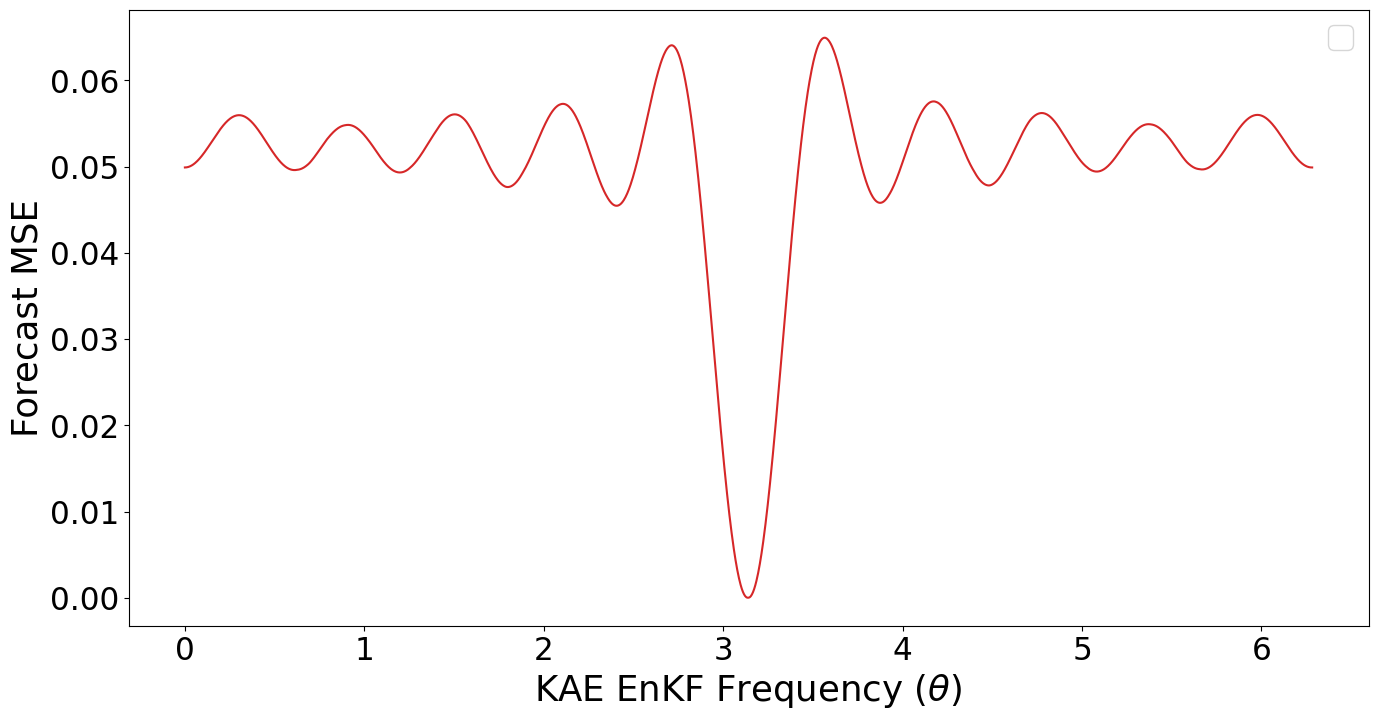}
     \end{subfigure}
     \caption{\centering\label{fig:local_minima} An example plot of the KAE forecast's mean square error (MSE) against  KAE EnKF’s frequency estimate. Data is generated for a discrete dynamical system, representing a simple sine wave with frequency $\theta = \pi$. Forecast for 10,000 states of this system are produced, with each forecast horizon being randomly selected from between 1 and 10 steps ahead. The resultant optimization landscape for the KAE EnKF's frequency estimate $\theta$ over the forecast's mean squared error, contains a strong global minimum at $\theta = \pi$ as expected, however also exhibits a further 9 local minima across the frequency's spectrum due to the use of 10 different forecast horizons in the optimization procedure.}
\end{figure}

\textbf{Optimization Procedure: }
The optimization procedure is still performed predominantly using stochastic gradient descent, however as noted in \cite{henning_lange_fourier_to_koopman}, when optimizing frequencies across multiple forecast horizons, the optimization landscape becomes pitted with local minima. An example of this can be seen in Figure \ref{fig:local_minima}, where the optimization landscape for forecasts errors produced by data from a simple sine wave is shown.

Using stochastic gradient descent over an optimization landscape with many local minima can produce defective models \cite{local_minima_neural_net_cases}, as the training procedure can get trapped in suboptimal local minima. The authors of \cite{henning_lange_fourier_to_koopman} developed a solution to this problem, by alternating training between iterations of stochastic gradient descent, followed by a global optimization of each of the Koopman approximators frequencies $\theta$. 

The global frequency optimization procedure in \cite{henning_lange_fourier_to_koopman} was applied over the entire dataset, however due to memory restrictions we instead apply it over a sample batch, generated in the same way as the batches used in the stochastic gradient descent steps. Frequencies $\theta_i$ are then selected one at a time, with all the model's other parameters kept fixed, so each frequency's global loss function can be calculated. For each triple $(\mathbf{x}_k, \mathbf{x}_{k+dt},dt)$ taken from the batch, the frequency's loss function will be periodic $2\pi/dt$, hence only values of $\theta_i$ in this domain need be calculated for the given triple. An appropriate number of intervals $s$ to represent the loss function over the domain $\lbrack 0, 2\pi/dt)$ is determined, and as a default we set $s=100$. We then calculate the mean squared error in the forecasts for the frequency value at each interval
\begin{equation}
    \mathcal{L}_{\theta_i}(2\pi j/dt s) = 
    \|\mathbf{x}_{k+dt}-\mathbf{\Tilde{h}}_{\Tilde{\phi}}(\mathbf{K}_{\lambda, \theta_i=  2\pi j/dt s}^{dt}\mathbf{h}_\phi(\mathbf{x}_k)) \|_2^2,
    \label{eqn:freq_loss}
\end{equation}
where $j\in{(0,\hdots,s-1)}$, and $\mathbf{K}_{\lambda, \theta_i= 2\pi j/dt s}$ is the current Koopman approximator with frequency $\theta_i=  2\pi j/dt s$.

This process is repeated for each triple in the batch, generating a frequency loss function for each, which then need to be summed into an overall loss function for the total batch. This amalgamation however is non-trivial, as each triple has a possibly different forecast horizon $dt \in{(1, \hdots, p)}$, meaning their frequency losses $\mathcal{L}_{\theta_i}$ from equation \eqref{eqn:freq_loss} will have been calculated at the different intervals $2\pi j/dt s$. By applying the Fourier transform \cite{fourier_transform} to each triple's frequency loss, we can efficiently resample to combine the losses, before using the inverse Fourier transform to return a total frequency loss function $\mathcal{L}_{\theta_i}^{total}$ over the domain $\lbrack 0, 2 \pi)$. The model's new estimate $\theta_i^{new}$ of frequency $\theta_i$ is then determined from the total frequency loss function by finding
\begin{equation}
    \operatorname{argmin}_{\theta_i^{new}}(-|\mathcal{L}_{\theta_i}^{total}(\theta_i^{new})| - \operatorname{median}(\mathcal{L}_{\theta_i}^{total}(\theta_i^{new})))).
    \label{eqn:freq_opt_target}
\end{equation}

Our aim is to find a new frequency estimate $\theta_i^{new}$ that is strongly present within the Koopman representation of our system. As this global frequency optimization occurs during the training process, it will likely be run before good estimates of the Koopman eigenfunctions $\mathbf{h}_\phi$ and modes $\Tilde{\mathbf{h}}_{\Tilde{\phi}}$ have been found. This results in what was described by \cite{henning_lange_fourier_to_koopman} as the unknown phase problem. Poor Koopman eigenfunction estimates mean the Koopman autoencoder would be unable to determine the correct phase for a given state $\mathbf{x}_k$ when represented in the frequency latent space $\mathbf{h}_\phi(\mathbf{x}_k)$. When fitting for the frequency of a wave, while an in-phase wave of the same frequency would minimize the loss $\mathcal{L}_{\theta_i}^{total}$, a wave of the same frequency but significantly out of phase would instead maximize the loss. By centring, then taking the negative absolute value of the loss function in equation \eqref{eqn:freq_opt_target}, the possibility of an incorrect latent space phase is mitigated, as the new frequency can be either the minimizer or maximizer of $\mathcal{L}_{\theta_i}^{total}$. Before the result of equation \eqref{eqn:freq_opt_target} is set as the new estimate for frequency $\theta_i$, it is checked that it is not equal to zero, as constant values should be handled via the bias terms in the model's network. Also, if $\theta_i^{new}$ is within a user-defined tolerance of the other current Koopman approximator frequency estimates $\theta_j$ for $j \neq i$, then it is not used to avoid all frequencies in $\mathbf{K}_\lambda$ converging to the same dominant frequency. The intervals in the domain $\lbrack 0, 2 \pi)$ are ordered by their result in equation \eqref{eqn:freq_opt_target}, and the first frequency that meets all the required conditions is set as the new frequency $\theta_i$ for eigenvalue $\lambda_i$ of Koopman approximator $\mathbf{K}_\lambda$. The global frequency optimization procedure is then repeated for the rest of the model's frequencies.

Training cycles between stochastic gradient descent for a number of iterations, before reapplying the global frequency optimization until both optimization schemes converge, or the loss on a separate validation dataset begins to consistently increase with each epoch, signalling overtraining.

\subsection{KAE EnKF Algorithm}
The algorithm, to train and forecast using the KAE EnKF, proceeds as follows.

\textbf{Step 1:}
First, the Koopman autoencoder model must be trained using the alterations detailed in the above subsection, optimizing parameters $\phi$, $\Tilde{\phi,}$ and $\lambda$, to construct encoder $\mathbf{h}_\phi$, decoder $\Tilde{\mathbf{h}}_{\Tilde{\phi}}$, and Koopman approximator $\mathbf{K}_\lambda$ respectively. Training is performed on states $\mathbf{x}_k$ for $k = 1, \hdots, m$, and this is referred to as the spin-up phase.

\textbf{Step 2: }
Before we can use the Koopman autoencoder model trained in step 1 to assimilate new data produced by the system, we must refactor the problem in terms of the EnKF equations \eqref{eqn:kf_prop_eq} and \eqref{eqn:kf_update_eq}.

\textbf{Ensemble Initialization: }
The KAE EnKF's main aim is to efficiently improve the current estimates of the system's model and state as new data becomes available. State estimates are assimilated efficiently by performing all calculations in the $r$ dimensional latent space $\Tilde{\mathbf{x}}_k$, as opposed to the much larger, $n$ dimensional, full state space $\mathbf{x}_k$. Our Koopman autoencoder model of the system is iteratively improved by fixing the encoder and decoder networks $\mathbf{h}_\phi$ and $\Tilde{\mathbf{h}}_{\Tilde{\phi}}$, but modelling the eigenvalues $\lambda$ of the Koopman approximator $\mathbf{K}_\lambda$ as time varying parameters. This leads us to define an initial state $\mathbf{z}_0\in{\mathbb R^{3r}}$ for the ensemble Kalman filter
\begin{equation}
            \mathbf{\mathbf{z}_0}
         = \left[\begin{array}{ccccccccc}
            - & \Tilde{\mathbf{x}}_m & - 
            &\tau_1 & \hdots & \tau_r
            &\theta_1 & \hdots & \theta_r
             \end{array}\right]^T.
    \label{eqn:z=xtmtautheta}
\end{equation}
By explicitly filtering the modulus and argument of the eigenvalues, we are able to work with entries of the filter's state in the Reals. To generate an initial ensemble from this state, we require an initial covariance matrix $\mathbf{P}_0\in{\mathbb R^{3r\times3r}}$. We define this as
\begin{equation}
    \mathbf{P}_{0} = \left[\begin{array}{c|c|c}
         \alpha_1 I_{r} &\mathbf{0} &\mathbf{0}  \\
         \hline
         \mathbf{0}& \alpha_2 I_{r} &\mathbf{0} \\
         \hline
         \mathbf{0} &\mathbf{0} & \alpha_3 I_{r}
    \end{array}\right],
    \label{eqn:kae_init_cov}
\end{equation}
with $\alpha_1 \gg \alpha_3 \gg \alpha_2$. Constant $\alpha_1$ governs the uncertainty over the latent space state $\Tilde{\mathbf{x}}_m$, while $\alpha_2$ and $\alpha_3$ define the uncertainty over the system's eigenvalue's modulus and arguments respectively. A system's state generally changes significantly more quickly than its parameters, leading to the conditions $\alpha_1 \gg \alpha_2$ and $\alpha_1 \gg \alpha_3$. We also know that a small change in a system's eigenvalue modulus can lead to vast changes in its future behaviour, hence for relatively stable systems we assume their eigenvalue's arguments are more prone to change than their modulus, hence $\alpha_3 \gg \alpha_2$. More advanced uncertainty initialization schemes for the system's eigenvalues are possible, for example by training multiple models via bagging in a similar manner to \cite{bagging_DMD}, however for our purposes this simple, three parameter initialization scheme is sufficient. Independent draws are then taken from $\mathcal{N}(\mathbf{z}_0,\mathbf{P}_0)$, until a sufficiently large ensemble has been generated. The requisite ensemble size $N$ will vary between specific applications, and is strongly but not exclusively influenced by the filter state's dimension $3r$. A default ensemble size $N = 100$ is used, which is increased/decreased to balance efficiency against the risk of filter degeneration as required.

\textbf{Propagation: }
For convenience, we introduce notation $\mathbf{z}^{i:j}_{k} \in \mathbb{R}^{j-i+1}$ to denote the $i$th through to the $j$th element of $\mathbf{z}_k$ where $i \leq j$
\begin{equation}
    \mathbf{z}^{i:j}_{k} =
    \left[\begin{array}{c c c}
            \mathbf{z}^{i}_{k} & \hdots & \mathbf{z}^{j}_{k}
             \end{array}\right]^T
             \label{eqn:dmdenkfzki:j=zki...zkj}.
\end{equation}
By rewriting the matrix $\mathbf{K}_\lambda$ in terms of the newly assimilated eigenvalue modulus and arguments $\mathbf{z}_k^{(r+1):3r}$, we have that at time $m+k$, $\tau_i = \mathbf{z}_k^{r+i}$ and $\theta_i = \mathbf{z}_k^{2r+i}$. The action of the Koopman approximator over the latent space then becomes
\begin{equation}
    \mathbf{B}_i = \left[\begin{array}{cc}
         \mathbf{z}_k^{(r+i)}\cos(\mathbf{z}_k^{(2r+i)})
         &-\mathbf{z}_k^{(r+i)}\sin(\mathbf{z}_k^{(2r+i)})  \\
         \mathbf{z}_k^{(r+i)}\sin(\mathbf{z}_k^{(2r+i)})
         &\mathbf{z}_k^{(r+i)}\cos(\mathbf{z}_k^{(2r+i)}) 
    \end{array}\right], \quad
    \mathbf{K}_{\mathbf{z}_k} = \left[\begin{array}{cccc}
         \mathbf{B}_1& &\mathbf{0} \\
         &\ddots& \\
         \mathbf{0}& &\mathbf{B}_r 
    \end{array}\right].
    \label{eqn:kl=enkf_polar_eig_blocks}
\end{equation}
The propagation equation \eqref{eqn:kf_prop_eq} in the EnKF framework can now be written as
\begin{equation}
    \mathbf{z}_{k+1} = \left[\begin{array}{c|c}
         \mathbf{K}_{\mathbf{z}_k}&\mathbf{0}  \\
         \hline
         \mathbf{0}& I_{2r}
    \end{array}\right]
    \mathbf{z}_k + \mathbf{w}_k.
    \label{eqn:kae_zk+1=z+w}
\end{equation}
Operator $\mathbf{K}_{\mathbf{z}_k}$ advances the assimilated latent space $\mathbf{z}_k^{1:r}$ forward by 1 time step in the system's dynamics, as it did in the Koopman autoencoder framework. The identity matrix $I_{2r}$ sends the Koopman eigenvalues to themselves, as although we anticipate them to vary over time, how they will change is assumed to be unknown. Vector $\mathbf{w}_k\in{\mathbb R^{3r}}$ is normally distributed with $\mathbf{w}_k \sim \mathcal{N}(\mathbf{0},\mathbf{Q}_k)$, where $\mathbf{Q}_k$ determines the uncertainty when propagating the system 1 time step forward. Covariance matrix $\mathbf{Q}_k$ is constructed similarly to the initial covariance matrix $\mathbf{P}_0$ from equation \eqref{eqn:kae_init_cov}, with
\begin{equation}
    \mathbf{Q}_{k} = \left[\begin{array}{c|c|c}
         \alpha_4 I_{r} &\mathbf{0} &\mathbf{0}  \\
         \hline
         \mathbf{0}& \alpha_2 I_{r} &\mathbf{0} \\
         \hline
         \mathbf{0} &\mathbf{0} & \alpha_3 I_{r}
    \end{array}\right].
    \label{eqn:kae_prop_noise}
\end{equation}
Noise when propagating the system is assumed to be uncorrelated, hence $\mathbf{Q}_{k}$ is constructed to be diagonal. Constant $\alpha_4$ shares the same properties as $\alpha_1$ from equation \eqref{eqn:kae_init_cov}, with $\alpha_4 \gg \alpha_2$ and $\alpha_4 \gg \alpha_3$, as again $\alpha_4$ governs variance over the system's state, which is assumed to deviate over a smaller timescale than the system's parameters.

\textbf{Measurement: }
The filter's measurement equation \eqref{eqn:kf_update_eq} is then
\begin{equation}
    \mathbf{y}_{k} = \left[\begin{array}{c|c}
         I_{r}&\mathbf{0}  \\
    \end{array}\right]
    \mathbf{z}_k
    +
    \mathbf{v}_k,
    \label{eqn:kae_yk=z+v}
\end{equation}
where measurements $\mathbf{y}_k\in{\mathbb R^{r}}$ are observations at time $m+k$ of the system's state encoded into the latent space $\mathbf{y}_k = \mathbf{h}_\phi(\mathbf{x}_{m+k}) = \Tilde{\mathbf{x}}_{m+k}$. The latent space representation of new state observations $\mathbf{x}_{m+k}$ is employed to improve efficiency, as the latent space dimension is assumed to be significantly smaller than the full state space dimension ($r \ll n$). Vector $\mathbf{v}_k \in {\mathbb R^{r}}$ governs the noise in these measurements, and we set it to be a normally distributed random variable $\mathbf{v}_k \sim \mathcal{N}(\mathbf{0},\mathbf{R}_k)$, with $\mathbf{R}_k = \alpha_5 I_{r}$. The matrix $\mathbf{R}_k$ is diagonal for the sake of simplicity, and the ratios between constants $\alpha_5$ and $\alpha_4$ from equation \eqref{eqn:kae_prop_noise} govern the relationship between the filter's confidence in new measurements of the state relative to the model's forecasts. As the filter observes new measurements over the latent space, it is possible that even if the measurement noise in the full state variables $\mathbf{x}_{m+k}$ is normally distributed, after applying the non-linear transform $\mathbf{h}_\phi$ to produce $\mathbf{y}_k$, the noise profile may have deviated significantly from the normal distribution. Methods to constrain the latent space in a way that guaranteed normal measurement noise in its variables, for example by augmenting the Koopman autoencoder's loss function, were unsuccessful. We have found empirically however, that approximating the latent space measurement noise $\mathbf{v}_k$ as an independent, normally distributed variable, is sufficient for our desired applications.

\textbf{Step 3: }
As new state measurements $\mathbf{x}_{m+k}$ become available, we encode them into the latent space to produce observations in the format the filter expects $\mathbf{y}_k = \mathbf{h}_\phi(\mathbf{x}_{m+k})$. The EnKF then assimilates $\mathbf{y}_k$, producing an updated ensemble of estimates for the system's latent state and parameters $\mathbf{z}_k$. This new ensemble estimate combines the model's forecast from the filter's prior state $\mathbf{z}_{k-1}$, with the information from the fresh noisy measurement $\mathbf{y}_k$, as detailed in the previous chapter's subsection on the inner workings of the EnKF.

\textbf{Step 4: }
The first $r$ elements of each ensemble member's state $\mathbf{z}_k^{1:r}$ then represents the latent state of the original system $\Tilde{\mathbf{x}}_{m+k}$. The second $r$ elements $\mathbf{z}_k^{r+1:2r}$ represent the current estimate of the system's eigenvalues modulus, and the third $\mathbf{z}_k^{2r+1:3r}$ their corresponding arguments. We can thus recover an ensemble of estimates for the full state of the system at time $m+k$, by calculating $\mathbf{x}_{m+k} = \Tilde{\mathbf{h}}_{\Tilde{\phi}}(\mathbf{z}_k^{1:r})$ for each ensemble member. Forecasting an ensemble of future states for the full system $dt$ steps ahead can then be generated as
\begin{equation}
    \mathbf{x}_{m+k+dt} = \Tilde{\mathbf{h}}_{\Tilde{\phi}}
    (\mathbf{K}_{\mathbf{z}_k}^{dt} \mathbf{z}_k^{1:r}).
\end{equation}
A point estimate of the future state $\mathbf{x}_{m+k+dt}$ can be calculated using the ensemble's mean, and the spread of the ensemble's distribution can be used to quantify uncertainty within the forecast.

A major difference between the KAE EnKF and a more typical fusion of a system model with a data assimilation filter, is that the KAE EnKF applies noise directly to the KAE's latent state $\Tilde{\mathbf{x}}_{m+k}$ (represented within the filter as $\mathbf{z}_k^{1:r}$), as opposed to the full system's state $\mathbf{x}_{m+k}$. This ensures the majority of the filter's resources are employed representing uncertainty in the latent states, as these are the states the model has determined to be the most dynamically relevant description of the system.

By modelling uncertainty as additive Gaussian noise applied directly to the latent states, the KAE EnKF ensemble's distribution in the latent space has a typically Gaussian profile. Hence, the system's true latent state value will likely be well captured within the ensemble's spread. When full state estimates are required, the approximately Gaussian latent state ensemble is transformed via the potentially highly nonlinear KAE decoder, back into the full state space as $\mathbf{x}_{m+k} = \Tilde{\mathbf{h}}_{\Tilde{\phi}}(\mathbf{z}_k^{1:r})$. The decoder $\Tilde{\mathbf{h}}_{\Tilde{\phi}}$ will deform the input Gaussian distribution to a new, arbitrarily complexly distribution of full states $\mathbf{x}_{m+k}$. This distribution can then be utilized for uncertainty quantification within the system's full state space. Ensemble methods are capable of representing arbitrarily complex distributions, provided sufficient ensemble members are included, and by ensuring uncertainty is effectively modelled in the system's most dynamically relevant directions, we attempt to reduce the number of ensemble members required to properly represent the full state output distribution $\Tilde{\mathbf{h}}_{\Tilde{\phi}}(\mathbf{z}_k^{1:r})$.

\section{Synthetic Applications}\label{s:synthetic}

\subsection{Synthetic data generation}

To test the performance of the KAE EnKF, we first generate synthetic data, to apply it and other comparable modelling techniques over. To make the synthetic data as specific as possible to the types of system the KAE EnKF is designed to model, we require high-dimensional, nonlinear measurements produced from an underlying, low-dimensional system with time varying parameters.

The latent dynamics of our system are governed by a 2-dimensional rotation matrix, with the angle of rotation $\theta_k$ increasing linearly from $\pi/128$ to $\pi/16$ over the course of the experiment's 1000 time steps. The underlying linear evolution of the system's latent state $\Tilde{\mathbf{x}}_k$ is thus
\begin{equation}
    \Tilde{\xb}_{k+1} = 
    \left[\begin{array}{cc}
         \cos({\theta_k})&-\sin({\theta_k})  \\
         \sin({\theta_k})&\cos({\theta_k})
    \end{array}\right]
    \Tilde{\xb}_k,
    \quad
    \Tilde{\xb}_1 =     
    \left[\begin{array}{c}
         1\\
         0
    \end{array}\right],
    \label{eqn:kaexk+1=rotx}
\end{equation}
with $\theta_k = \pi/128 + \frac{(k-1)(7\pi/128)}{999}$ and $k = (1,...,1000)$.

\begin{figure}[htbp]
\begin{subfigure}[b]{\textwidth}
         \includegraphics[width=\textwidth]{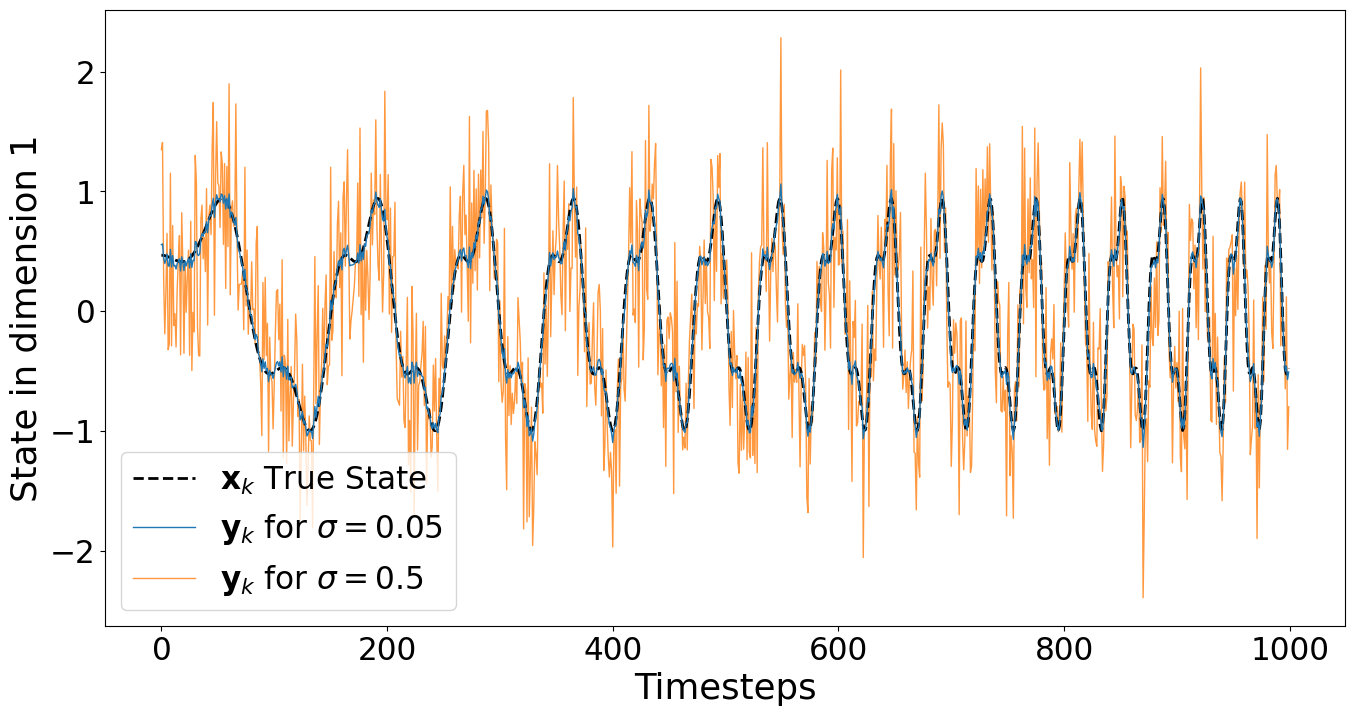}
     \end{subfigure}
     \caption{\centering\label{fig:kae_sin5_bothnoise}An example trajectory for the first dimension of the synthetic system's full state $\mathbf{x}_k$, with nonlinearity parameter $\nu = 3$, and both low $\sigma=0.05$ and high $\sigma=0.5$ levels of measurement noise, over the course of 1000 time steps.}
\end{figure}

A high-dimensional, nonlinear full state $\mathbf{x}_k$ is then generated from each latent state in equation \eqref{eqn:kaexk+1=rotx}, by raising $\Tilde{\mathbf{x}}_k$ to the power $\nu$ element-wise, before multiplying by matrix $\mathbf{A}\in{\mathbb R^{100\times2}}$, such that
\begin{equation}
    \mathbf{x}_k = \mathbf{A}{\Tilde{\mathbf{x}}_k}^\nu.
    \label{eqn:kaexk=Axk}
\end{equation}
The matrix $\mathbf{A}$ is randomly generated, with its entries being selected from a Uniform distribution $\mathcal{U}[0,1)$, and we select $\nu = 3$. Noisy measurements of the full system state at each time step $\mathbf{y}_k\in{\mathbb R^{100}}$ are assumed to be available, with
\begin{equation}
    \mathbf{y}_{k} = 
    \xb_k + \vb_k,
    \quad
    \vb_k \sim \mathcal{N}(\mathbf{0},\sigma^{2}I_{100}).
    \label{eqn:kaey=x+vvN0sI}
\end{equation}
Variable $\sigma$ is set to $0.05$ or $0.5$, to simulate low/high levels of measurement noise respectively. Example trajectories of $\mathbf{y}_k$ for each noise level over the course of all 1000 time steps can be seen in Figure \ref{fig:kae_sin5_bothnoise}.
The first 200 time steps are used for the spin-up phase of the KAE EnKF model, as described in step 1 of the KAE EnKF's algorithm. This synthetic data generation process is used in all synthetic experiments, with any modifications clearly stated at the start of the relevant experiment's subsection. The neural network architectures of the Koopman autoencoder's encoder and decoder used in each experiment are listed in in full, in the appendix.

\subsection{Comparison between full state and latent state Koopman autoencoder filtering}

Before comparing against other existing modelling techniques, we first investigate the decision in the KAE EnKF framework to filter the system's latent state, as opposed to its full measured state. The first benefit of latent state filtering is a vast improvement in the filter's efficiency. When assimilating a new data point with the EnKF, the algorithm's computational complexity if often dominated by the forming of intermediate covariance matrices \cite{enkf_comp_complexity}.This operation takes order $\mathcal{O}(n^2N)$, where $n$ is the filter's state size, $N$ is the number of ensemble members, and the filter's state size is traditionally much larger than its ensemble size $n \gg N$ \cite{understanding_the_enkf}. By applying the filter over the system's latent state of dimension $r \ll n$, this operation reduces to $\mathcal{O}(r^2N)$, precipitating a quadratic reduction in the algorithm's computational complexity.

To evaluate the differences in performance between the latent and full state filter KAE EnKFs, we utilize data generated from the synthetic, non-stationary sine wave as described in the previous subsection. The full state filter KAE EnKF is formed by first training the Koopman autoencoder on the spin-up data as normal. The filter's ensemble state estimate is then initialized using the full state estimate from the end of the spin-up phase as opposed to its encoded latent representation, this constitutes replacing $\Tilde{\mathbf{x}}_m$ with $\mathbf{x}_m$ in the KAE EnKF's equation \eqref{eqn:z=xtmtautheta}. To propagate the ensemble using the KAE EnKF, we apply the latent dynamics as specified in equation \eqref{eqn:kl=enkf_polar_eig_blocks}, however must first encode our full state to the latent space, and subsequently decode them to return the full state estimate forecast for the following time step. This transforms the propagation of the system's state section of the filter's state in equation \eqref{eqn:kae_zk+1=z+w} from $\mathbf{z}_{k+1} = \mathbf{K}_{\mathbf{z}_k} \mathbf{z}_k$ to 
\begin{equation}
    \mathbf{z}_{k+1} = \Tilde{\mathbf{h}}_{\Tilde{\phi}}(\mathbf{K}_{\mathbf{z}_k} \mathbf{h}_{\phi}(\mathbf{z}_k)).
    \label{eqn:fullstate_prop}
\end{equation}
Measurements provided to the filter for assimilation are simply the observations at time $m+k$ of the full system's state $\mathbf{x}_{m+k}$. The latent and full state filter KAE EnKF's were applied to data generated by the previously described synthetic system, with all computational parameters kept the same, and notably an ensemble size of $N = 100$.

\begin{figure}[!htbp]
\begin{subfigure}[b]{.49\textwidth}
         \includegraphics[width=\textwidth]{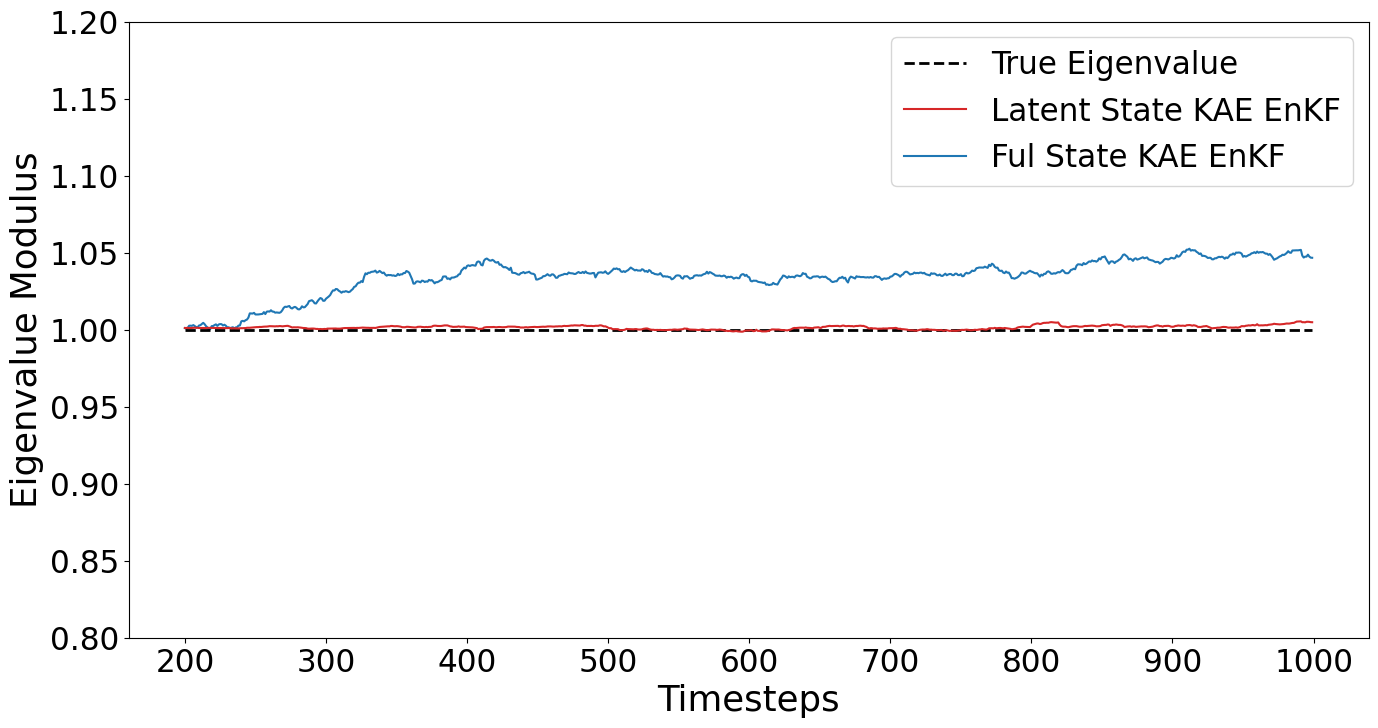}
                \caption{\centering\label{fig:low_latentfull_eigmod}The eigenvalue modulus estimates, for latent state and full state KAE EnKF's over an example data trajectory, with low $\sigma=0.05$ measurement noise.}
\end{subfigure}
\hfill     
\begin{subfigure}[b]{.49\textwidth}
         \includegraphics[width=\textwidth]{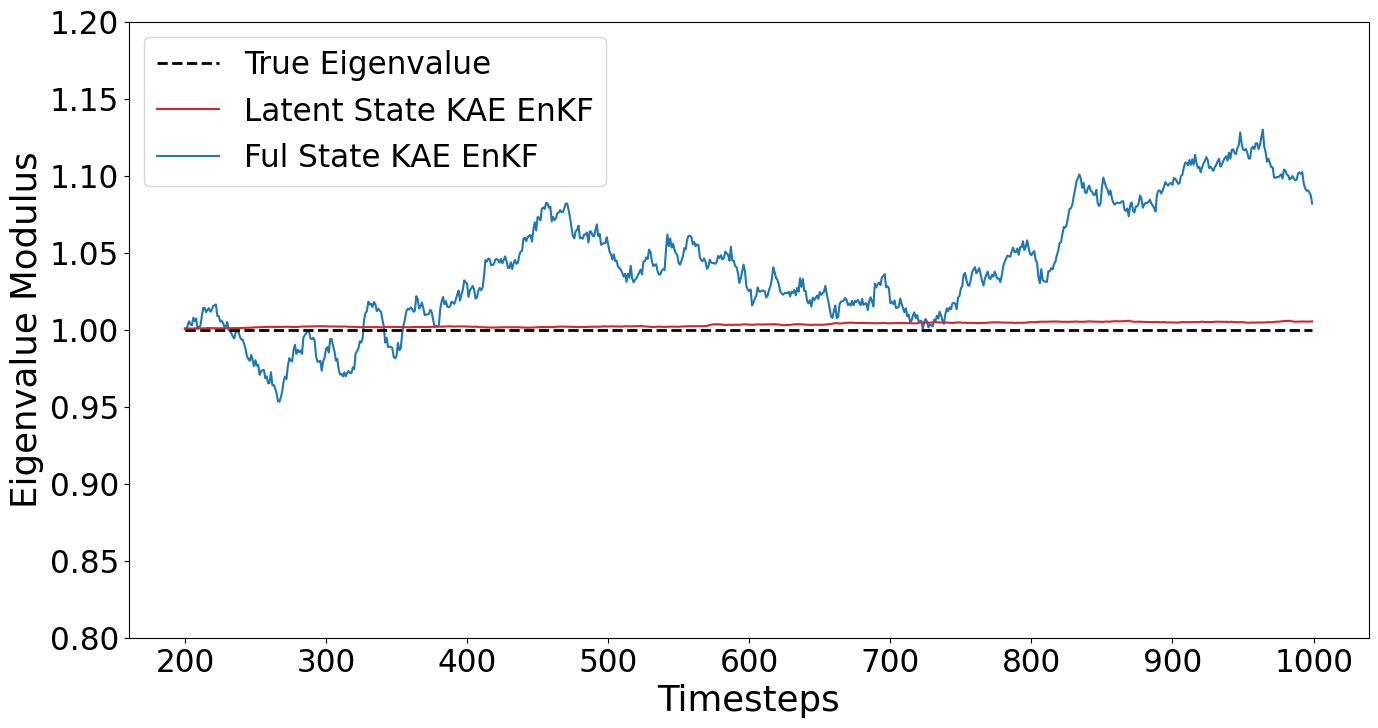}
                \caption{\centering\label{fig:high_latentfull_eigmod}The eigenvalue modulus estimates, for latent state and full state KAE EnKF's over an example data trajectory, with high $\sigma=0.5$ measurement noise.}
\end{subfigure}

\begin{subfigure}[b]{.49\textwidth}
         \includegraphics[width=\textwidth]{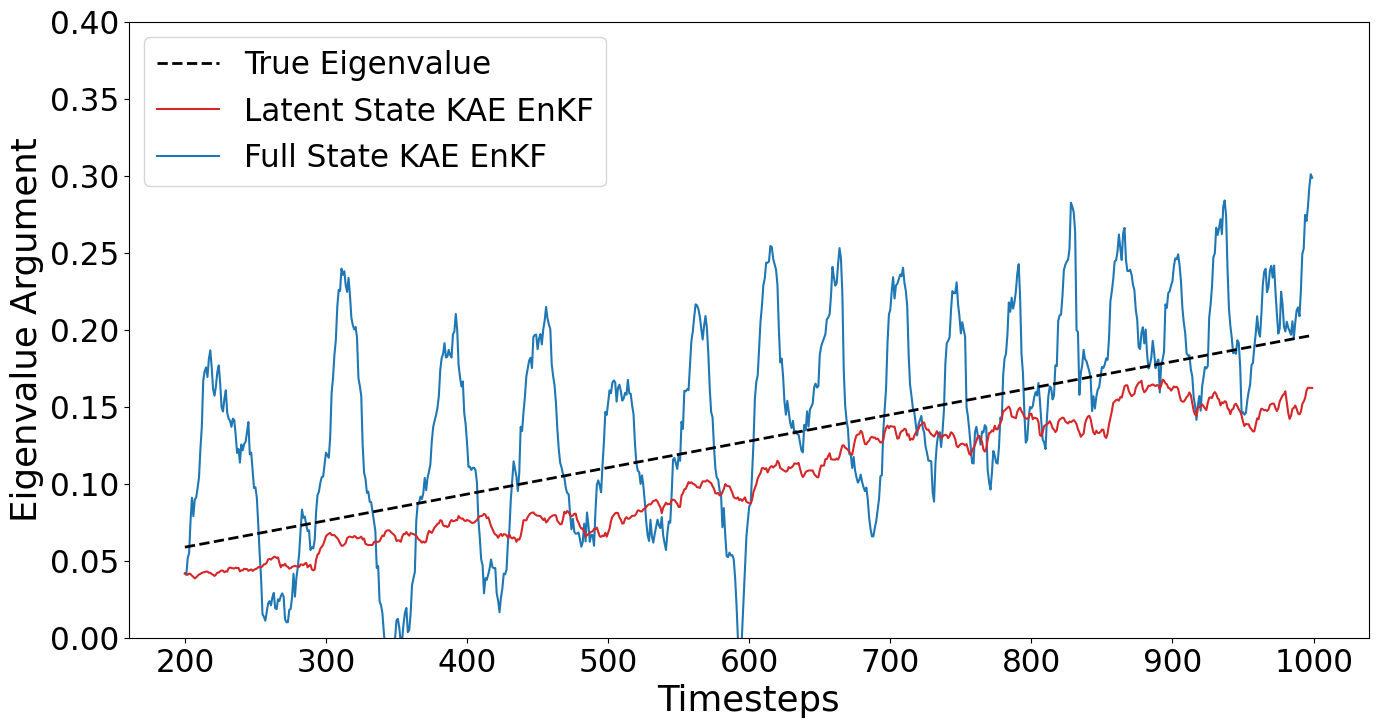}
                \caption{\centering\label{fig:low_latentfull_eigarg}The eigenvalue argument estimates, for latent state and full state KAE EnKF's over an example data trajectory, with low $\sigma=0.05$ measurement noise.}
\end{subfigure}
\hfill     
\begin{subfigure}[b]{.49\textwidth}
         \includegraphics[width=\textwidth]{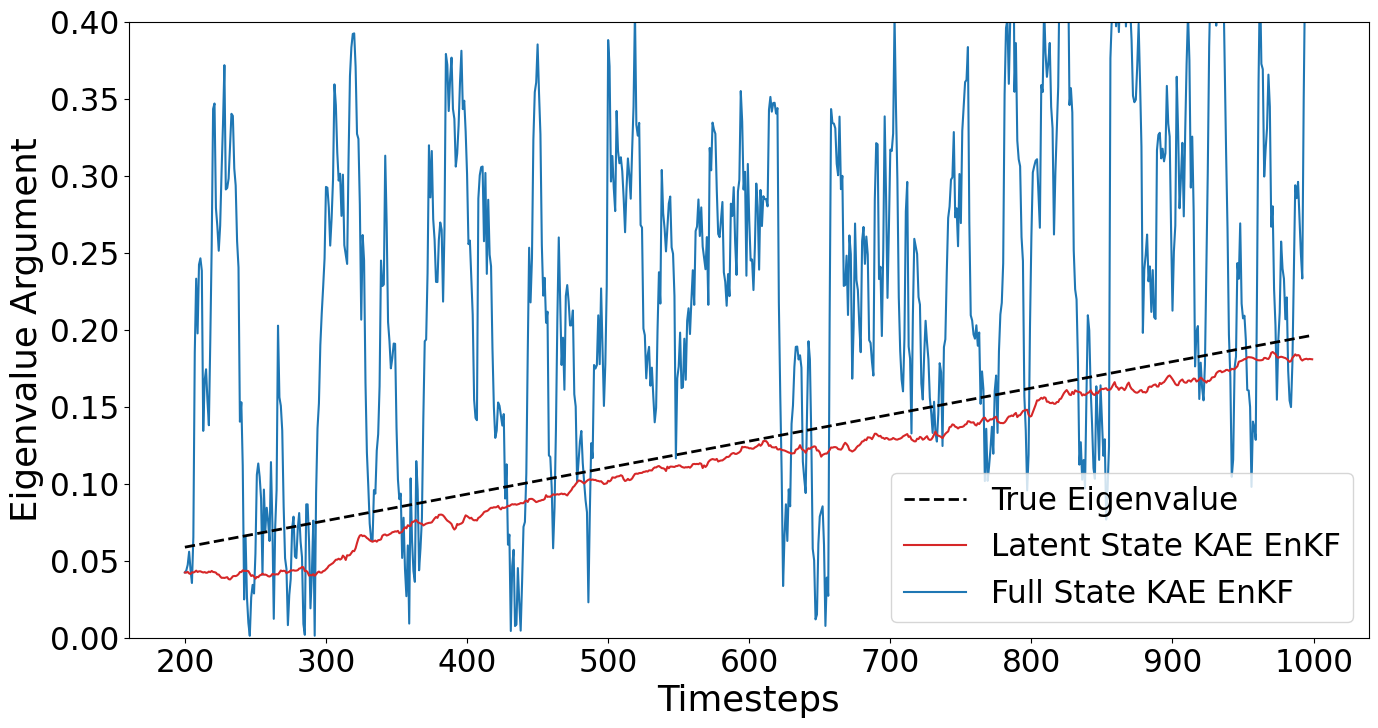}
                \caption{\centering\label{fig:high_latentfull_eigarg}The eigenvalue argument estimates, for latent state and full state KAE EnKF's over an example data trajectory, with high $\sigma=0.5$ measurement noise.}
\end{subfigure}

\begin{subfigure}[b]{.49\textwidth}
         \includegraphics[width=\textwidth]{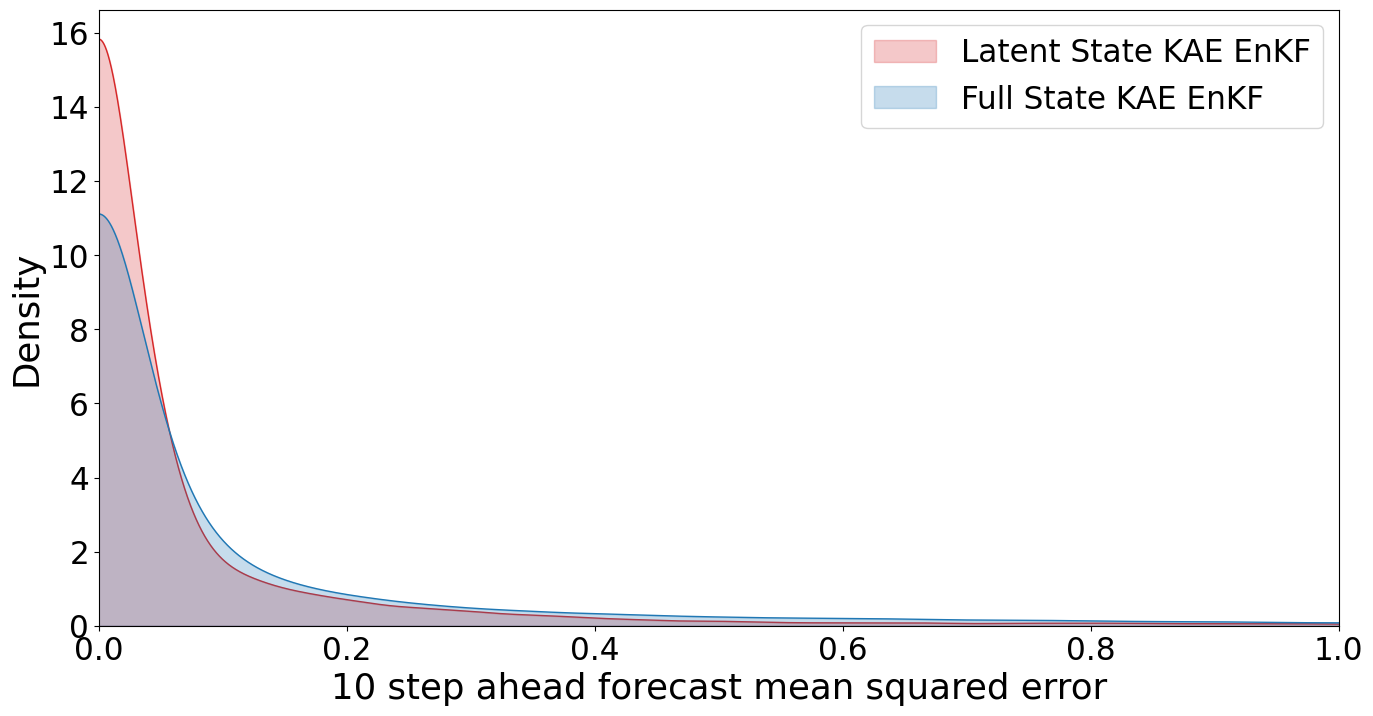}
                \caption{\centering\label{fig:low_latentfull_mse_dist}The mean squared error distributions for latent state and full state KAE EnKF's 10-step ahead forecasts of the system's full state, at low $\sigma=0.05$ measurement noise, calculated over the course of 100 runs.}
\end{subfigure}
\hfill     
\begin{subfigure}[b]{.49\textwidth}
         \includegraphics[width=\textwidth]{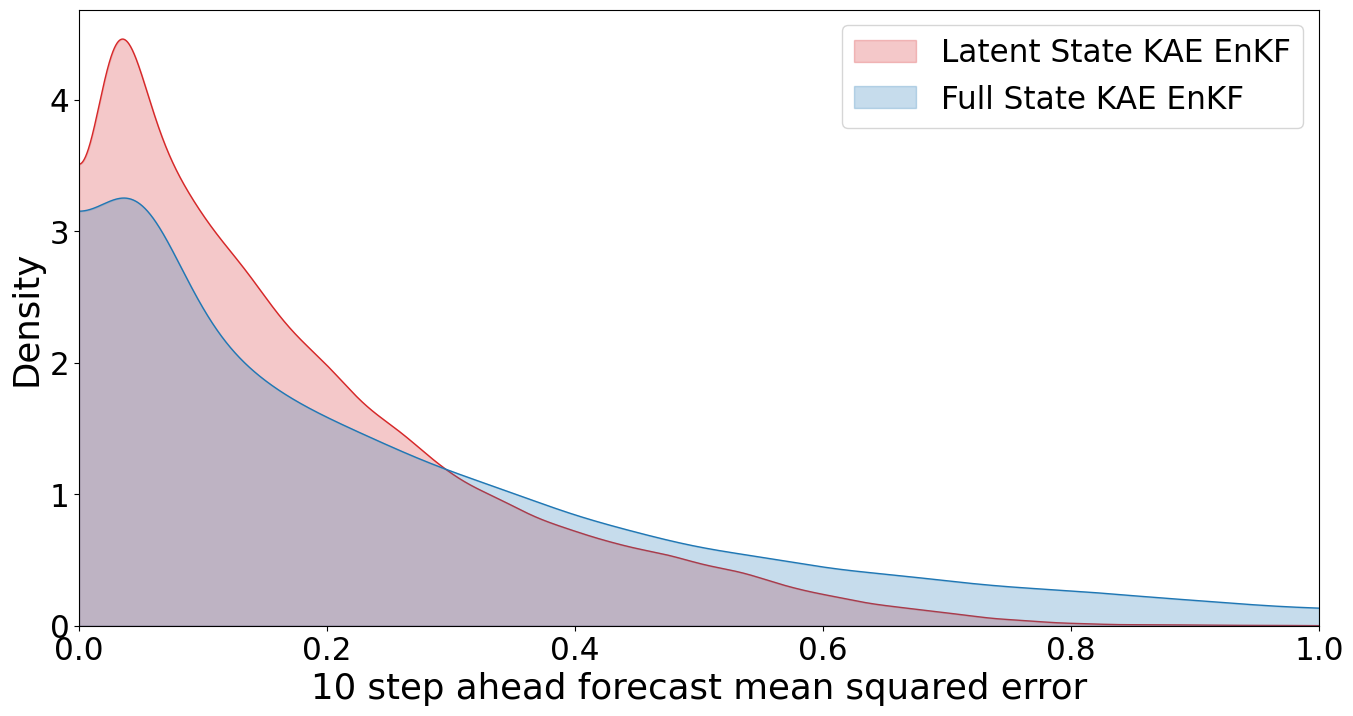}
                \caption{\centering\label{fig:high_latentfull_mse_dist}The mean squared error distributions for latent state and full state KAE EnKF's 10-step ahead forecasts of the system's full state, at high $\sigma=0.5$ measurement noise, calculated over the course of 100 runs.}
\end{subfigure}

     \caption{\centering\label{fig:latentfull_kaeenkf_performance}
     The eigenvalue modulus estimates,(top), argument estimates (middle), and mean squared 10-step ahead forecast error distributions (bottom), for latent state and full state KAE EnKF's, with low (left) and high (right) levels of measurement noise. The latent state KAE EnKF's parameter estimates and state forecasts are more stable and accurate than those of the full state KAE EnKF over all noise levels.}
\end{figure}

Figure \ref{fig:latentfull_kaeenkf_performance} shows the latent and full state filter KAE EnKF's performance when tracking the system's stationary eigenvalue modulus, linearly increasing eigenvalue argument, and forecasting 10 steps ahead its full state. At all noise levels and in all three metrics, the latent state KAE EnKF significantly outperforms its full state counterpart. The full state KAE EnKF is consistently unstable, producing erratic parameter estimates and wide error distributions, and this is strongly exacerbated as noise levels are increased. The latent state KAE EnKF experiences a relatively small degradation in performance as noise increases, with the main decline being in the 10-step ahead forecasts seen in Figures \ref{fig:low_latentfull_mse_dist} and \ref{fig:high_latentfull_mse_dist}. The high noise error distribution peaks at a greater value than that of the low noise, due to the increased irreducible error that further noise introduces. 

A major issue with the full state KAE EnKF, is that when tracking the system parameters shown in Figures \ref{fig:low_latentfull_eigmod}, \ref{fig:high_latentfull_eigmod}, \ref{fig:low_latentfull_eigarg} and \ref{fig:high_latentfull_eigarg}, the model's parameter estimates appear at times to drift away from their true values. A strongly motivating factor in applying data assimilation techniques to the Koopman autoencoder, is to develop a flexible data-driven framework, that is capable of adapting to changes within the system in real-time. Hence, we now investigate the cause behind this apparent lack of responsiveness within the full state KAE EnKF.

\begin{figure}[htbp]
\begin{subfigure}[b]{\textwidth}
         \includegraphics[width=\textwidth]{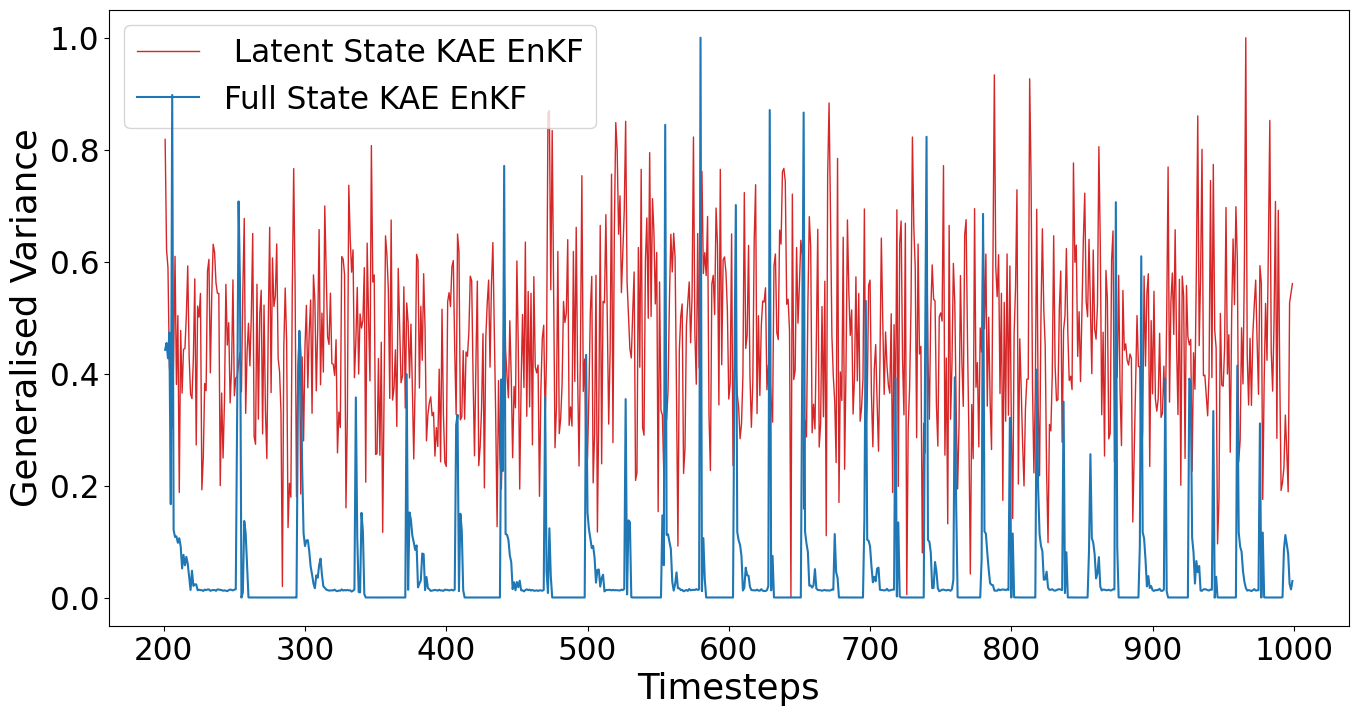}
     \end{subfigure}
     \caption{\centering\label{fig:fullstate_kaeenkf_genvar}The normalized, generalized sample variance of the latent state KAE EnKF's ensemble, and that of the full state KAE EnKF's ensemble when encoded into latent space. The latent state KAE EnKF's ensemble maintains a relatively stable spread, whereas the full state KAE EnKF's ensemble collapses for the majority of time steps, then develops sudden unstable spikes in its generalized variance at semiregular intervals.}
\end{figure}

The full state KAE EnKF's inability to effectively track the system, indicates that the filter's internal ensemble is in some way unable to properly represent the likely range of possible values the system's state and parameters could take. To investigate how the spread of the full state KAE EnKF's ensemble develops compared to that of the latent state KAE EnKF, we run both filters again under low measurement noise, this time with $N = 1000$ ensemble members, and plot the generalized covariance of their ensembles at each time step in latent space. For the full state KAE EnKF, this requires each ensemble member to be encoded into the latent space using $\mathbf{h}_\phi$, before its generalized covariance is then calculated. The graphs of these trajectories, with each set of generalized variances normalized to take values between 0 and 1, are displayed in Figure \ref{fig:fullstate_kaeenkf_genvar}.

\begin{figure}[htbp]
\begin{subfigure}[b]{\textwidth}
         \includegraphics[width=\textwidth]{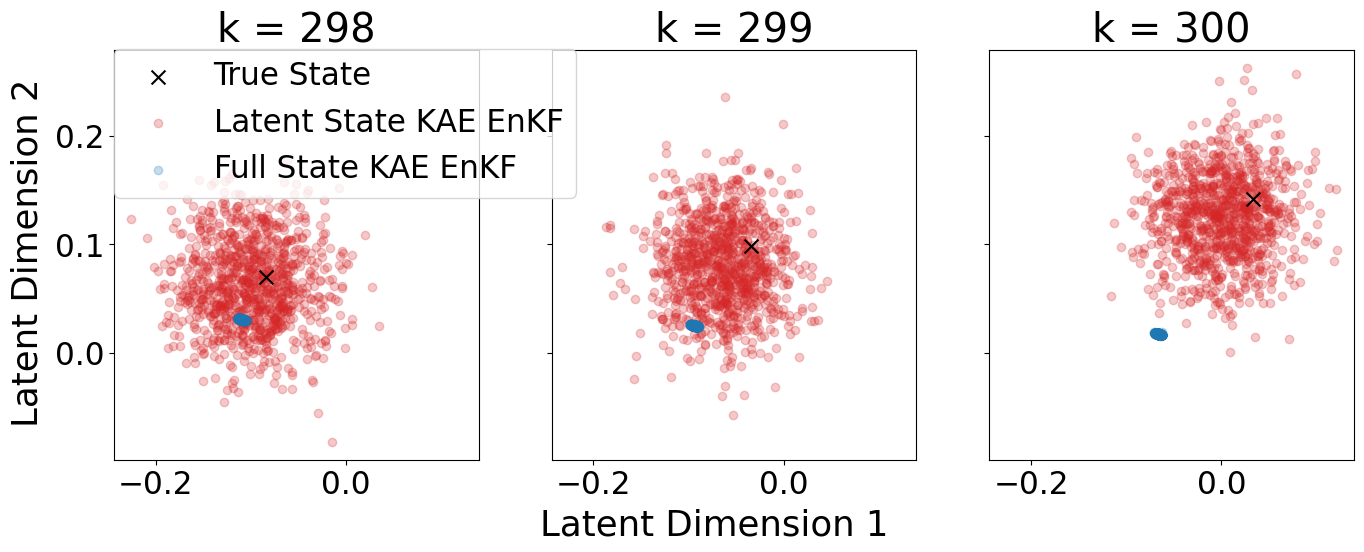}
     \end{subfigure}
     \caption{\centering\label{fig:fullstate_kaeenkf_ensemble1}The latent and full state KAE EnKF's 1000 ensemble members, displayed in latent space for the time steps just preceding the generalized covariance spike around $k = 300$ seen in Figure \ref{fig:fullstate_kaeenkf_genvar}. The latent state KAE EnKF's ensemble distribution is approximately centred on the true state, and is sufficiently wide in both latent directions to always include the true state when its centroid deviates from the true state’s value. The full state KAE EnKF’s ensemble is collapsed almost to a single point, and drifts further away from the true state as time progresses.}
\end{figure}

\begin{figure}[htbp]
\begin{subfigure}[b]{\textwidth}
         \includegraphics[width=\textwidth]{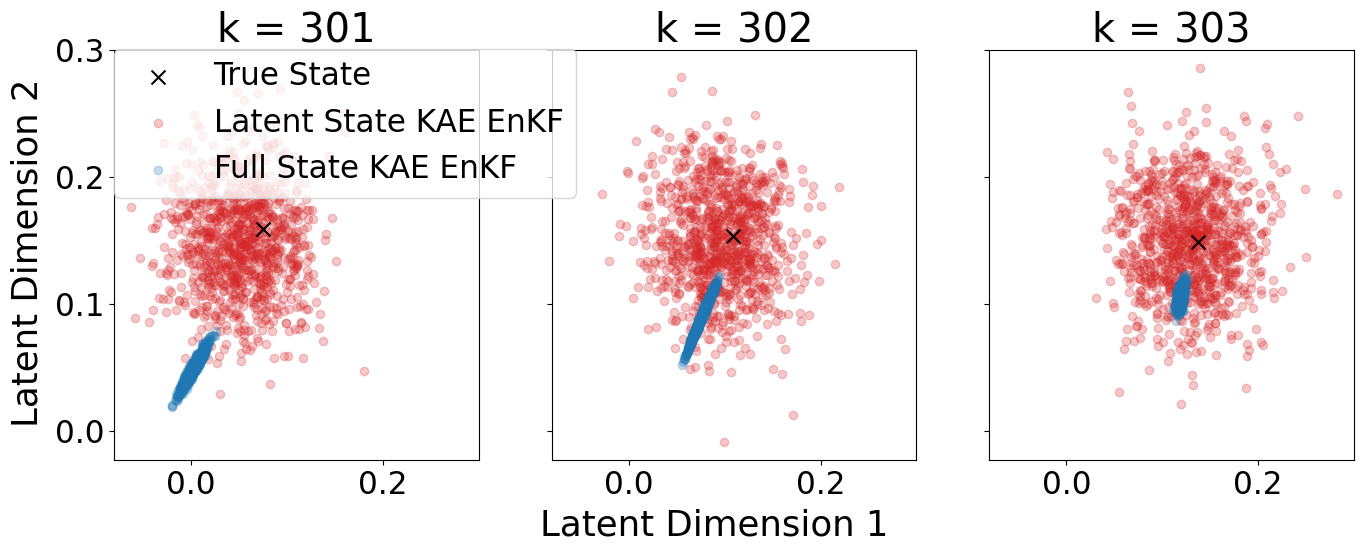}
     \end{subfigure}
     \caption{\centering\label{fig:fullstate_kaeenkf_ensemble2}The latent and full state KAE EnKF's 1000 ensemble members, displayed in latent space for the time steps during and just following the generalized covariance spike around $k = 300$ seen in Figure \ref{fig:fullstate_kaeenkf_genvar}. The latent state KAE EnKF's ensemble maintains similar properties as in Figure \ref{fig:fullstate_kaeenkf_ensemble1}. The full state KAE EnKF's ensemble drifts a significant distance from the true state's value, inducing its variance to rapidly grows in the direction between the ensemble and the true state at $k=301$. This allows for sufficient uncertainty in the model to adjust the ensemble's position towards the measured state, shrinking the disparity between the system propagation and measurement values, and hence tightening the ensemble's distribution again by $k=303$.}
\end{figure}

Here, we see the latent state KAE EnKF's ensemble maintains a relatively stable spread throughout the course of the experiment, as would be expected in the event of constant levels of system and measurement noise. The full state KAE EnKF's ensemble however, when encoded into latent space, completely collapses at the majority of time steps. It is interspersed with large, sudden spikes in the variance, that occur with a periodicity attuned to that of the synthetic system. These plateaus and sharp spikes in the generalized variance of the full state KAE EnKF's encoded ensemble are likely responsible for its poor performance, hence we view each ensemble's complete latent distributions around the spike at approximately $k = 300$, in Figures \ref{fig:fullstate_kaeenkf_ensemble1} and \ref{fig:fullstate_kaeenkf_ensemble2}.

The latent state KAE EnKF models uncertainty as additive Gaussian noise, applied directly to the latent space of the system. As such, the ensemble's point cloud in the latent space has a typically Gaussian profile in Figures \ref{fig:fullstate_kaeenkf_ensemble1} and \ref{fig:fullstate_kaeenkf_ensemble2}, with a variance heavily determined by the system and measurement covariance matrices input into the filter by the user. As such, it is relatively straightforward for the user to ensure uncertainty is effectively represented in the dynamically relevant latent states of the system, by appropriately calibrating the relevant covariance matrices provided to the filter. The full state KAE EnKF instead applies noise over the system's full state space, which is then encoded, to produce the ensemble's latent representation. This transition of the ensemble through the potentially highly nonlinear encoder $\mathbf{h}_\phi$, develops a more complex relationship between the noise profile supplied to the filter, and the subsequent ensemble's latent distribution.

Figure \ref{fig:fullstate_kaeenkf_ensemble1} shows the latent and full state KAE EnKF's ensembles in latent space, at the time steps preceding the sharp spike in the full state KAE EnKF's generalized covariance. The latent state KAE EnKF's ensemble has an adequate spread, however the full state KAE EnKF's ensemble has collapsed almost to a single point, resulting in the complete plateau of generalized variance visible before $k = 300$ in Figure \ref{fig:fullstate_kaeenkf_genvar}. The contractive effect of the encoder on the ensemble's mass is likely a by-product of the Koopman autoencoder's training procedure during the spin-up phase. By employing multiple regularization techniques within the model's optimization scheme, the Koopman autoencoder is conditioned to produce not only accurate, but stable forecasts of the system's state. To facilitate stable predictions, the Koopman autoencoder develops a strong denoising effect. Hence, when the intentionally whitened ensemble members are passed through the Koopman autoencoder, a large proportion of the added noise is removed by the encoder, and the system's state propagation equation \eqref{eqn:fullstate_prop} yields a densely packed ensemble. This causes the filter to disproportionately trust the model's forecast over the value of the newly assimilated measurement. Due to the imperfections in the model of the system, the overconfident full state KAE EnKF's latent state estimates drift away from the system's true value in Figure \ref{fig:fullstate_kaeenkf_ensemble1}, as the relatively insufficient trust in new measurements mean errors in the model's forecasts cannot be properly rectified.

Once the disparity between model forecasts and freshly assimilated measurements become significant, the full state KAE EnKF's ensemble uncertainty becomes sufficiently large to manifest in the ensemble's latent representation, as shown in Figure \ref{fig:fullstate_kaeenkf_ensemble2}. The ensemble's increased latent spread presents in Figure \ref{fig:fullstate_kaeenkf_genvar} as the sharp spikes in generalized variance. Increased uncertainty in the latent ensemble allows for the filter to more effectively adapt its state estimates to the new data, however induces instability within the model. Once the model's forecast and measurement values become more closely aligned, the ensemble's latent uncertainty reverts to its previous, insufficient levels, and remains here until imperfections in the model lead to significant forecast/measurement drift, and the cycle repeats.

\subsection{Comparison of the KAE EnKF against other nonlinear, iterative DMD variants}

\subsubsection{Single frequency system}
We now investigate how the KAE EnKF's performance compares against that of other nonlinear, iterative DMD variants. There are two main methods to alter DMD for systems that act nonlinearly in their measured states, extended \cite{extended_dmd} and kernel \cite{kernel_dmd} DMD. Extended DMD uses a collection of nonlinear basis functions $\mathbf{g}$, which act on the observed values that would traditionally represent the system's state $\mathbf{x}_k$ in the DMD framework, replacing them with what is referred to as a dictionary of observables $\mathbf{g}(\mathbf{x}_k)$.
The correct choice of basis functions $\mathbf{g}$ for a given application is an open problem, and this method's main drawback is an often exponential increase in the new DMD state's $\mathbf{g}(\mathbf{x}_k)$ dimension with respect to the number of basis functions. Kernel DMD addresses this issue with Extended DMD by employing the kernel trick \cite{kernel_trick}, representing the nonlinearity in the system by a choice of inner product over the dataset. This allows for data points to be represented in up to infinite dimensional spaces, with no additional computational complexity cost over that of the standard DMD algorithm. The optimal inner product space to use in a given application is also often not known.

We test the KAE EnKF against these methods using the same synthetic system as described above. Extended DMD is applied to the synthetic dataset, utilizing dictionaries of observables including the monomial basis, and Hermite polynomials as recommended by \cite{extended_dmd}, with $p-q$ quasi-norm order reduction as suggested in \cite{p-q_quasi_norm_order_reduction}. Kernel DMD is also implemented, with the standard $L^2$ inner product, as well as radial basis functions which are popular in other machine learning contexts \cite{kernel_machine_learning_methods}. The best performing basis functions however were time-delay embeddings, also known as Hankel-DMD \cite{hankel-dmd}. Augmenting a state with time-delay embeddings is a popular, uncomplicated technique for including additional information in a system's state \cite{dmd_time_delay_coords}, with a strong theoretical underpinning from Taken's theorem \cite{takens_embedding}. Time-delay embeddings have seen significant use across a variety of applications, for example in the eigensystem realization algorithm (ERA) \cite{era}, singular spectrum analysis (SSA) \cite{ssa}, and for chaotic systems the Hankel alternative view of Koopman (HAVOK) \cite{havok} and machine learning Hankel-DMD \cite{ml_hankel_dmd}.

\begin{figure}[!htb]
\begin{subfigure}[b]{.49\textwidth}
         \includegraphics[width=\textwidth]{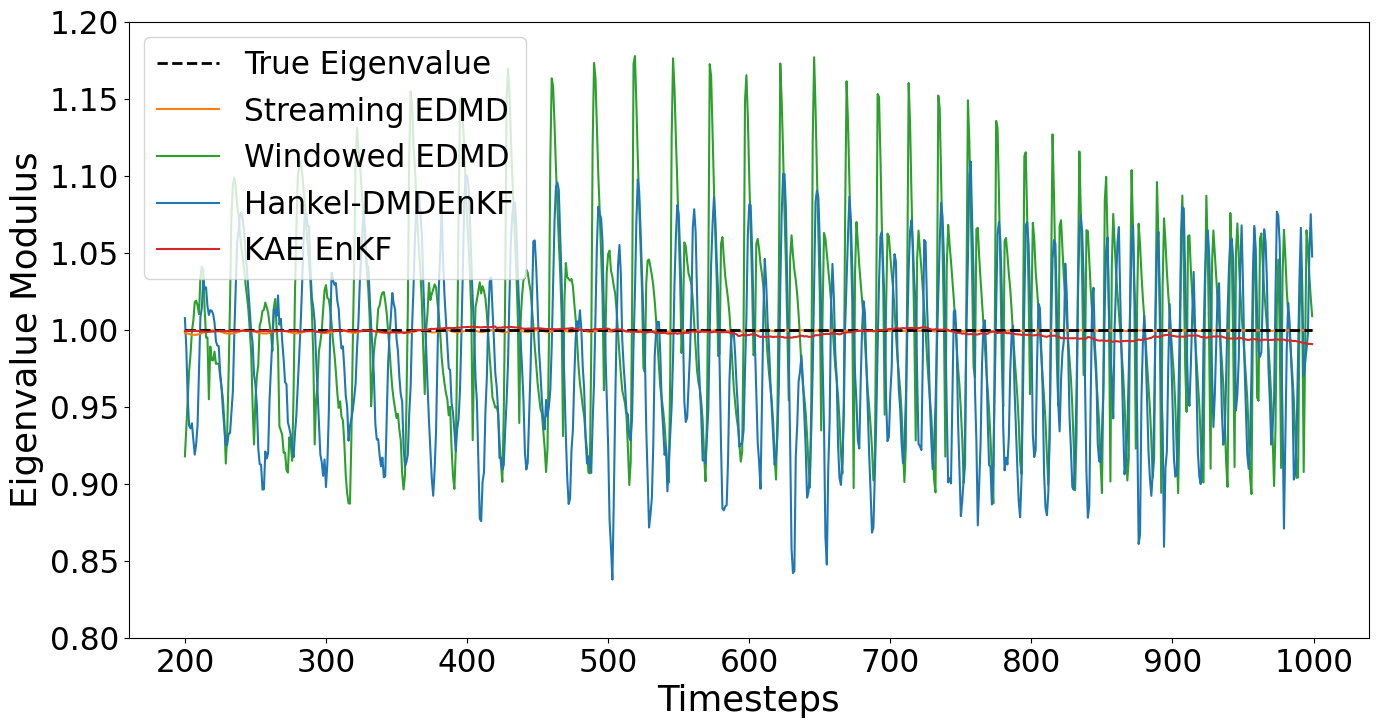}
                \caption{\centering\label{fig:low_edmd_mods}The eigenvalue modulus estimates, for each nonlinear iterative DMD variant over an example data trajectory, with low $\sigma=0.05$ measurement noise.}
\end{subfigure}
\hfill     
\begin{subfigure}[b]{.49\textwidth}
         \includegraphics[width=\textwidth]{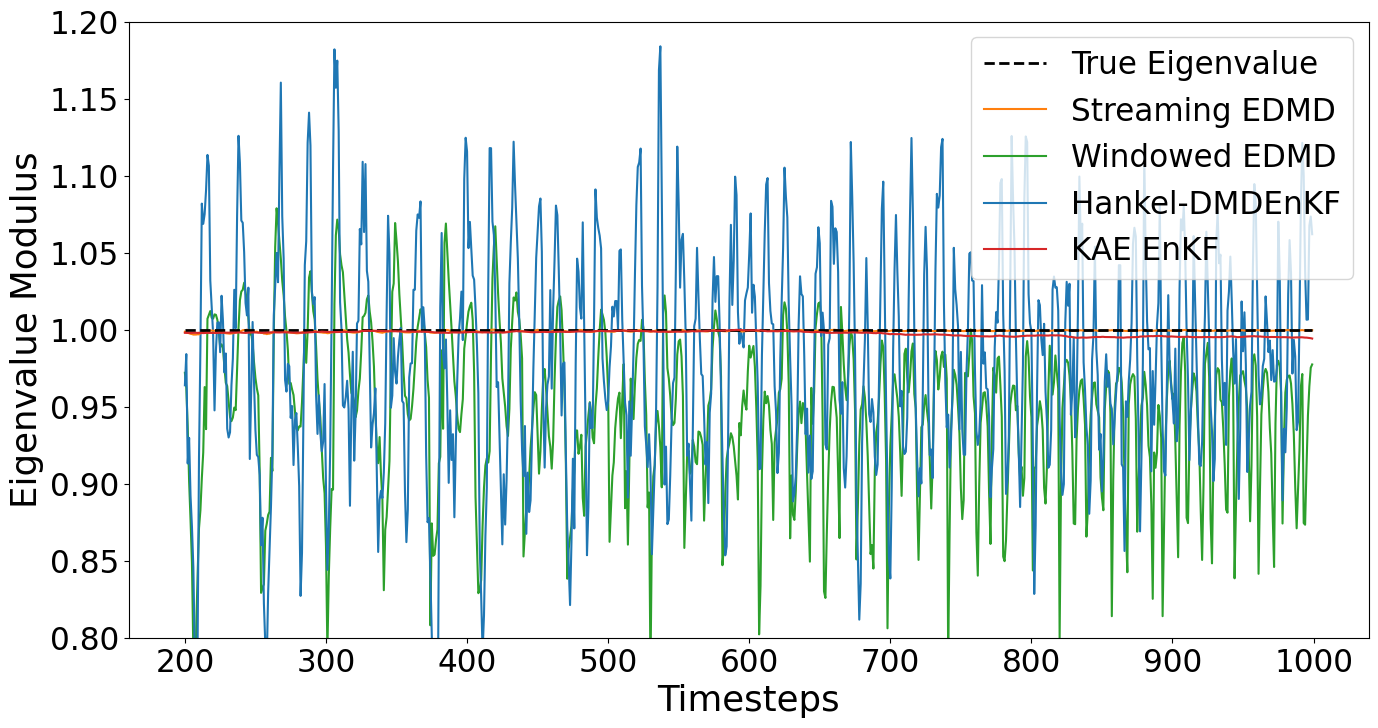}
        \caption{\centering\label{fig:high_edmd_mods}The eigenvalue modulus estimates, for each nonlinear iterative DMD variant over an example data trajectory, with high $\sigma=0.5$ measurement noise.}
     \end{subfigure}

\begin{subfigure}[b]{.49\textwidth}
         \includegraphics[width=\textwidth]{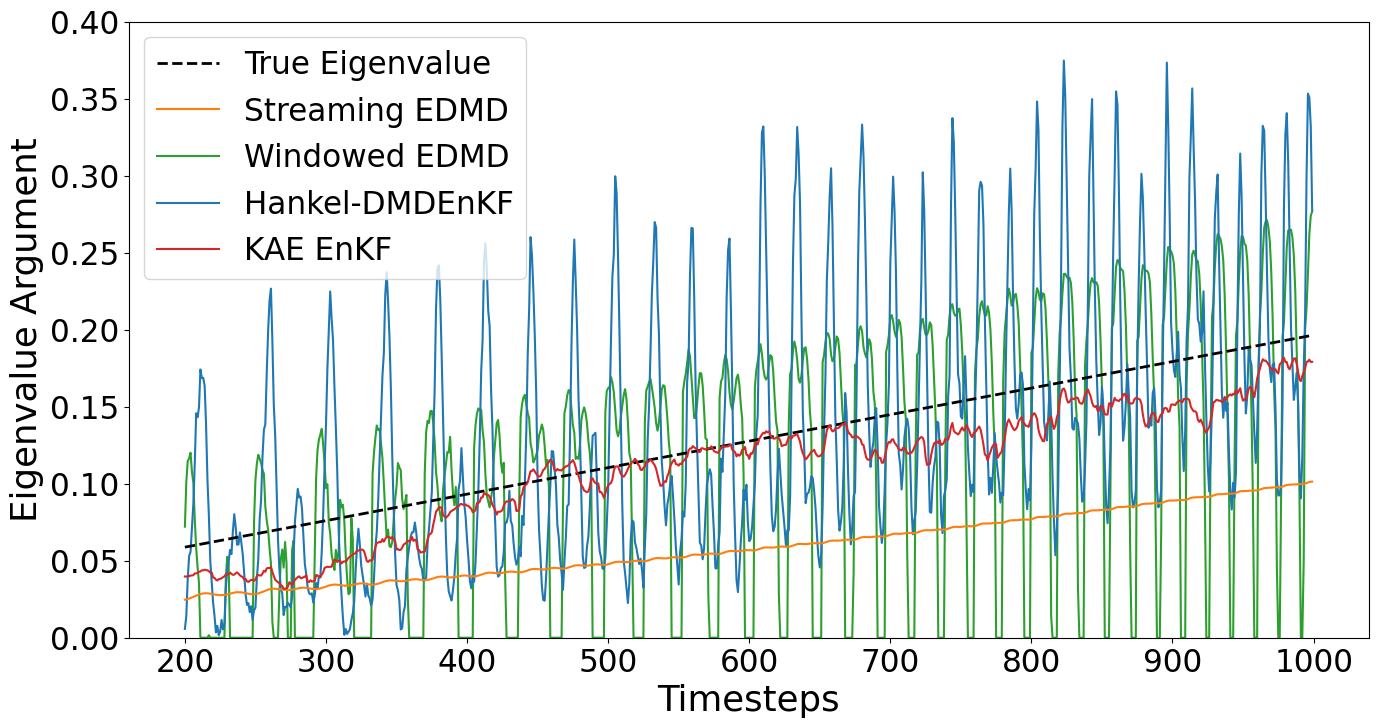}
                \caption{\centering\label{fig:low_edmd_args}The eigenvalue argument estimates, for each nonlinear iterative DMD variant over an example data trajectory, with low $\sigma=0.05$ measurement noise.}
\end{subfigure}
\hfill     
\begin{subfigure}[b]{.49\textwidth}
         \includegraphics[width=\textwidth]{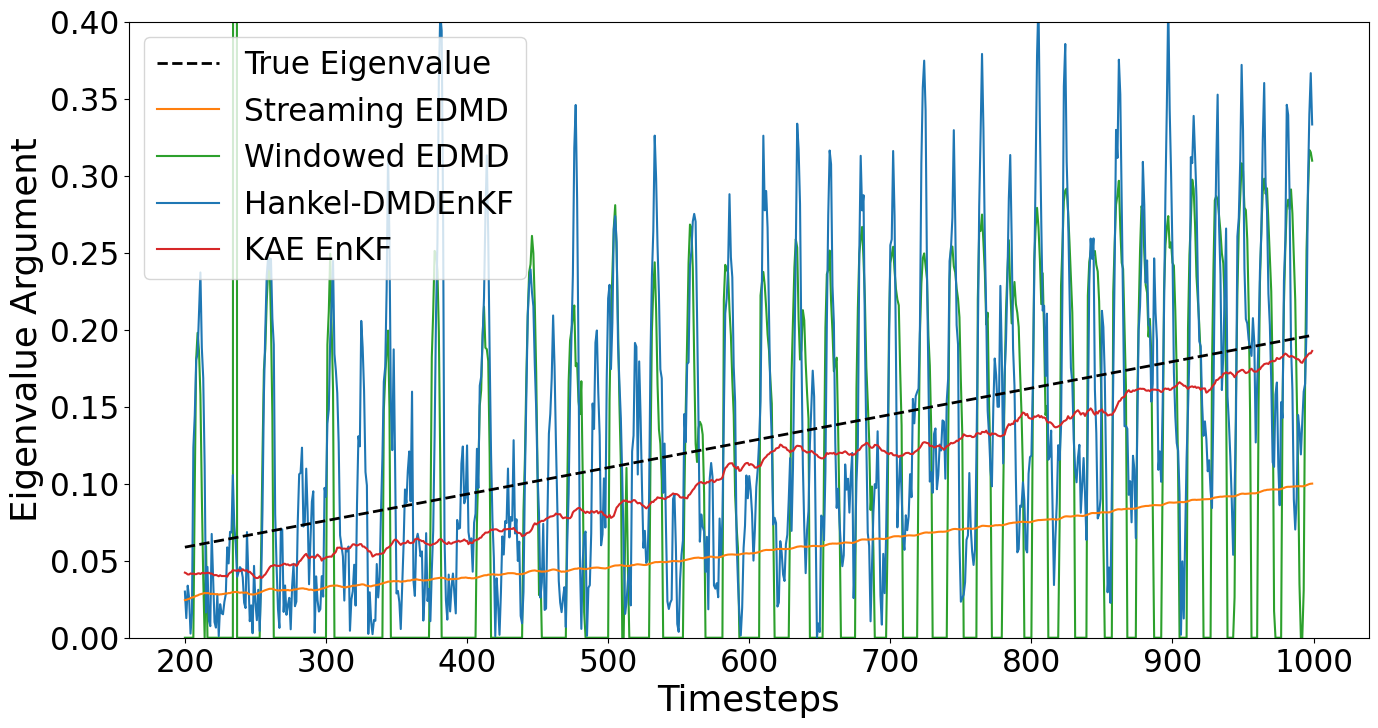}
        \caption{\centering\label{fig:high_edmd_args}The eigenvalue argument estimates, for each nonlinear iterative DMD variant over an example data trajectory, with high $\sigma=0.5$ measurement noise.}
     \end{subfigure}
     
     \caption{\centering\label{fig:edmd_arg_mod_trajs}The eigenvalue modulus (top) and argument (bottom) estimates, for each nonlinear iterative DMD variant over an example data trajectory, with low (left) and high (right) levels of measurement noise. Streaming EDMD achieves the most accurate eigenvalue modulus estimates, closely followed by the KAE EnKF. Streaming EDMD's argument estimates are stable but include a significant lag, whereas the KAE EnKF's estimates are slightly less stable but follow the changing parameter's value more effectively. Windowed EDMD and the Hankel-DMDEnKF are unstable over all parameters and noise levels, however produce less biased estimates of the non-stationary eigenvalue argument than streaming EDMD.}
\end{figure}

To form iterative variants of Hankel-DMD for comparison against the tracking and forecasting performance of the KAE EnKF, we replace the DMD models within the Streaming DMD and Windowed DMD algorithms with Hankel-DMD. These are henceforth referred to as Streaming EDMD and Windowed EDMD, to emphasize the time-delay embeddings purpose in these examples as the best performing set of nonlinear basis functions. We also utilize the Hankel-DMDEnKF, as the nonlinear version of the DMDEnKF described in \cite{dmdenkf}. Computational parameters for each experiment are set as follows; window size $w = 10$ for Windowed EDMD, delay-embedding dimension $d = 5$ for the Hankel component of Streaming EDMD, Windowed EDMD and the Hankel-DMDEnKF, and 200 time steps for the Hankel-DMDEnKF and KAE EnKF's spin-up phase as before. Figure \ref{fig:edmd_arg_mod_trajs} shows each method's tracking of the synthetic system's latent eigenvalue over the course of typical experiment trajectories.

In Figures \ref{fig:low_edmd_mods} and \ref{fig:high_edmd_mods}, we see the efficacy of each model, when tracking the eigenvalue's modulus, for low and high levels of measurement noise respectively. Streaming EDMD and the KAE EnKF both produce highly accurate modulus estimates over all noise levels, with Streaming EDMD slightly outperforming the KAE EnKF. Streaming EDMD's success here is due to its assumption that the underlying system producing the data is stationary, as this holds true for the eigenvalue's modulus, which is kept at a constant value of 1 throughout this experiment. As more data becomes available, Streaming EDMD's belief that the modulus is 1 intensifies, resulting in eigenvalue modulus estimates of almost exactly 1 by the time the final data point has been processed. The KAE EnKF also produces excellent modulus estimates, which we attribute to its flexible variance settings, allowing for different levels of uncertainty to be assigned for the model's eigenvalue modulus and arguments. By setting the modulus uncertainty parameter significantly below its argument counterpart ($\alpha_2 \ll \alpha_3$), we encourage the model to produce stable forecasts, while also permitting it to adjust its modulus estimates if the need for this change is evidenced in the newly assimilated data. Windowed EDMD and the Hanekl-DMDEnKF both produce significantly less accurate modulus estimates. As noise levels are increased, Windowed EDMD develops a bias to underestimate the modulus, however neither method's performance becomes significantly less stable. This is likely due to the known noise robustifying effect time-delay embeddings have, during the data rank truncation stage within each algorithm.

Figures \ref{fig:low_edmd_args} and \ref{fig:high_edmd_args} show each model's ability to track the latent system's eigenvalue argument $\theta_k$ as it linearly increases, over the course of a typical experiment's trajectory.
Here, Streaming EDMD's assumption of a stationary system works against it, producing stable eigenvalue argument estimates that lag significantly behind the true system's value. By weighting all available data points equally, the models estimates are heavily biased by older measurements, generated at a time when the system's non-stationary parameters were very different to their current values. The KAE EnKF performs excellently at both noise levels, with a small bias to underestimate the argument, as expected in the case of tracking an increasing parameter. As the filtering stage of the KAE EnKF incorporates the measurement noise explicitly into its estimates, the model places less trust on each individual new measurement in the high noise case. This further stabilizes the KAE EnKF's eigenvalue estimates, as the model's current system parameter estimates are held more strongly. This comes at the cost of a reduction in responsiveness to changes within the system, leading in this example, to a small increase in bias to its eigenvalue argument estimates. Windowed EDMD and the Hankel-DMDEnKF perform similarly to one another, producing estimates that are significantly less stable than for the other two methods. The stability of both methods only slightly degrades between the low and high measurement noise cases, which we again attribute to the denoising effect of employing a time-delay basis. The Hankel-DMDEnKF exhibits the least bias in its eigenvalue argument estimates of all the methods tested, over both noise levels. At low levels of noise, Windowed EDMD's estimates contain a small bias comparable to those of the KAE EnKF, however as noise is increased the Windowed EDMD estimates become more sporadic, and harder to draw useful insights from. At both noise levels, Windowed EDMD is regularly unable to identify that the system even contains a complex conjugate pair of eigenvalues, and this occurs more frequently at higher noise levels.

For all graphs in Figure \ref{fig:edmd_arg_mod_trajs}, the Hankel-DMDEnKF and Windowed EDMD exhibit highly periodic errors in their eigenvalue estimates. This is likely caused by both methods producing imperfect models that capture some, but not all, of the nonlinearity in the synthetic system. The system's latent state, raised to the power $\nu=3$ element-wise, experiences sharper peaks and flatter plateaus than those of an equivalent linear system. Thus, if the Hankel-DMDEnKF and Windowed EDMD generate models closer to the linear case, their eigenvalue modulus and argument estimates will increase at each peak, and decrease at each plateau. This causes a systematic pattern of the methods over then underestimating the system's eigenvalues, attuned to the increasing period of the true latent system's eigenvalues.

\begin{figure}[!htbp]
\begin{subfigure}[b]{.49\textwidth}
         \includegraphics[width=\textwidth]{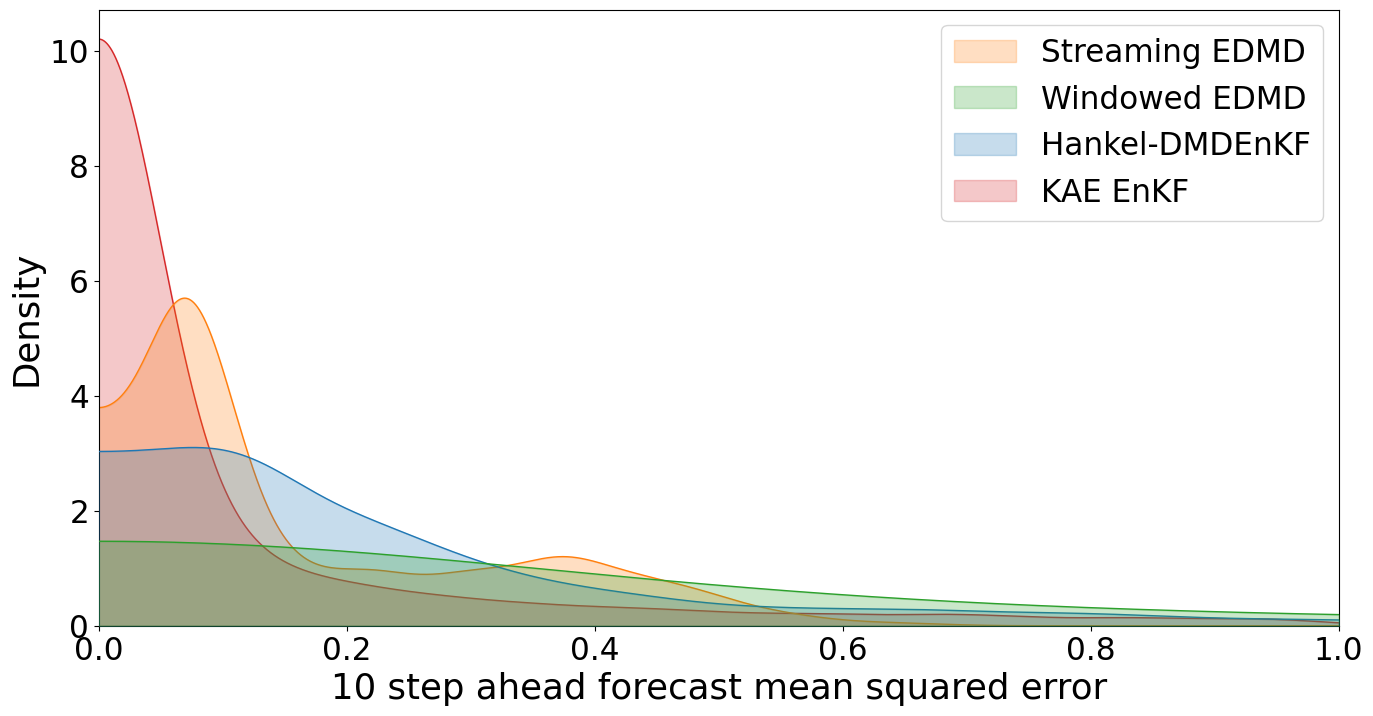}
                \caption{\centering\label{fig:low_edmd_mse_dist}The mean squared error distributions for each nonlinear iterative DMD variant's 10-step ahead forecasts of the full state, at low $\sigma=0.05$ levels of measurement noise, calculated over the course of 10 runs.}
\end{subfigure}
\hfill     
\begin{subfigure}[b]{.49\textwidth}
         \includegraphics[width=\textwidth]{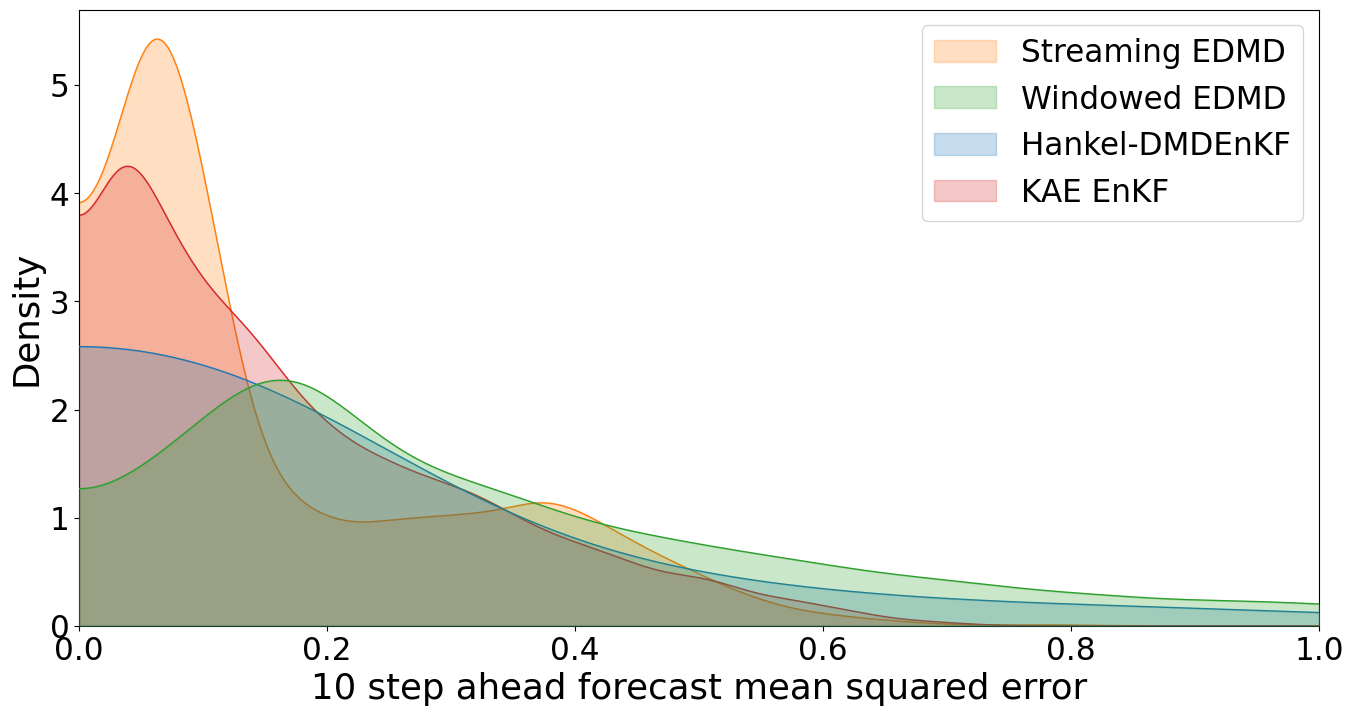}
        \caption{\centering\label{fig:high_edmd_mse_dist}The mean squared error distributions for each nonlinear iterative DMD variant's 10-step ahead forecasts of the full state, at high $\sigma=0.5$ levels of measurement noise, calculated over the course of 10 runs.}
     \end{subfigure}
     
     \caption{\centering\label{fig:edmd_forecast_mse_dist}The mean squared error distributions for each nonlinear iterative DMD variant's 10-step ahead forecasts of the full state, at low (left) and high (right) levels of measurement noise, calculated over the course of 10 runs. The KAE EnKF produces significantly more accurate forecasts than all other methods at low noise, however for high noise streaming EDMD's error distribution is tighter than the KAE EnKF's, but peaks at a greater error value. The Hankel-DMDEnKF and windowed EDMD are both unstable and produce wide error distributions, with the Hankel-DMDEnKF being the better performing of the two methods.}
\end{figure}

Having evaluated each method's ability to track the system's eigenvalues, we now investigate their capacity for forecasting the full state of the system. Each model generates a 10-step ahead forecast $\hat{\mathbf{x}}_k$ of the 100-dimensional system state $\mathbf{x}_k$, and the mean squared errors $\|\hat{\mathbf{x}}_k - \mathbf{x}_k\|_2^2$ are calculated for each time step after the spin up phase. The error distributions for the results over 10 runs of the experiment are displayed in Figure \ref{fig:edmd_forecast_mse_dist}, with the error's non-negative boundary condition enforced using a reflection of their KDE's density about 0 \cite{kde_boundary_reflection}.

At low levels of noise, the KAE EnKF's forecast are by far the most accurate, with the majority of errors, over all time steps, densely populated around 0. As noise increases, the errors still peak at a value close to 0, however are more spread out than before. The only difference in how the KAE EnKF is trained for each noise level is the filter's observation covariance is set to reflect the measurement noise in each dataset. Due to the framework's flexibility, there are many potential ways the KAE EnKF could be tweaked to mitigate the model's reduction in forecast accuracy as noise increases. It could be robustified by, using time-delay embeddings on the input data, training the spin-up Koopman autoencoder model over more epochs with a slower learning, adding artificial noise as in denoising autoencoders, or penalizing the network activation's Jacobian as done in contractive autoencoders. Streaming EDMD achieves consistent, relatively low forecast errors, that are somewhat independent of the level of noise in the data, due to the model's stability. This noise agnosticism is a boon in the high noise case, however the model's unavoidable delay bias when modelling any system with a non-stationary parameter imposes a hard upper limit on its forecast's performance. As such, even when forecasting a system from noise-free data, the model will produce similar forecast error distributions to those in Figure \ref{fig:edmd_forecast_mse_dist}. As noise-free forecasting is a significantly easier problem, this would result in a significant underperformance, compared to forecasts generated by the other methods. The Hankel-DMDEnKF produces unstable forecasts for both noise levels, although as noise increases the forecasts become slightly more unstable. Windowed EDMD performs the worst of all the methods, over all noise levels. Its forecasts are highly unstable, and as noise increases, we see a peak in its error distributions arise at $\sim 0.2$, which corresponds to the increased incidence of the model failing to identify a complex conjugate pair of eigenvalues in the system.

\subsubsection{Multiple frequency system}
\begin{figure}[htbp]
\begin{subfigure}[b]{\textwidth}
         \includegraphics[width=\textwidth]{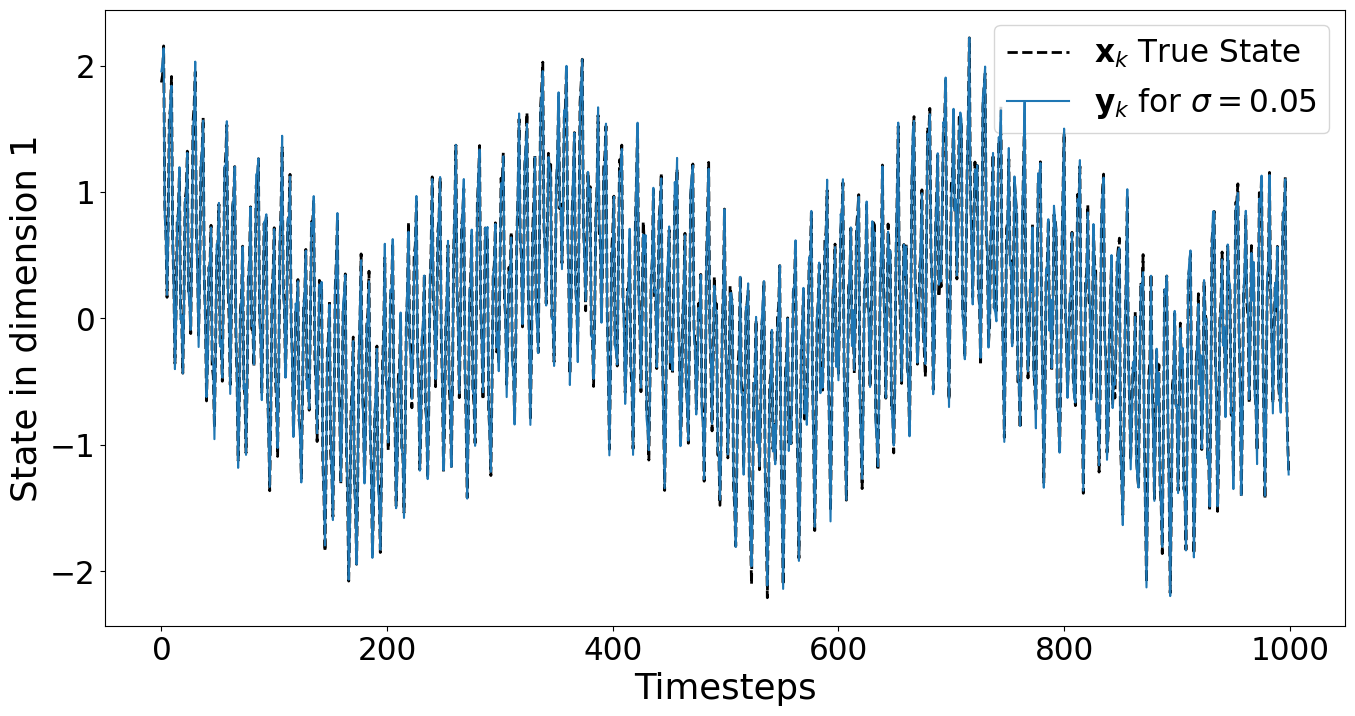}
     \end{subfigure}
     \caption{\centering\label{fig:multifreq_sin5_data_time_series}An example trajectory for the first dimension of the full state $\mathbf{x}_k$, generated by a synthetic system that contains multiple frequencies, one of which is non-stationary and linearly increasing over the course of the experiment. Measurement noise is kept to low levels at $\sigma = 0.05$, for all 1000 time steps.}
\end{figure}

For completeness, we now compare the KAE EnKF against the same nonlinear, iterative DMD variants, when tracking and forecasting a system that contains multiple underlying frequencies. To evaluate the performance of each method, we generate data for a synthetic system with the same latent rotation matrix structure defined by equation \eqref{eqn:kaexk+1=rotx}, used in the previous synthetic examples. In this case however, we run equation \eqref{eqn:kaexk+1=rotx} three times with three different sets of frequencies $\theta_k$, generating three sets of 1000 time step, 2-dimensional latent state data. Two of the frequencies $\theta_k = 2\pi/7$ and $\theta_k = 2\pi/365$, we keep constant for all time steps, and the third frequency linearly increases over the course of the experiment, with $\theta_k = 2\pi/30 + \frac{(k-1)(4\pi/30)}{999}$ for $k = (1,...,1000)$. These frequencies were chosen to mimic weekly, monthly, and yearly timescales often present in real time-series data, with the monthly frequency acting as a non-stationary system parameter to evaluate each method's system tracking capabilities. To form full states $\mathbf{x}_k$, the values of the latent states for each of the three frequency's rotational systems at time $k$ are stacked to form a latent state $\Tilde{\mathbf{x}}_k \in{\mathbb R^{6}}$. The latent state is then cubed to introduce nonlinearity, and multiplied by randomly generated matrix $\mathbf{A}\in{\mathbb R^{100\times6}}$ to form a full state $\mathbf{x}_k$ of much larger dimension than its latent state, as per equation \eqref{eqn:kaexk=Axk}. Low levels of measurement noise $\sigma =0.05$ only are added to the full states to produce measurements $\mathbf{y}_k$, and an example of the data produced by the synthetic system is demonstrated in Figure \ref{fig:multifreq_sin5_data_time_series}.

\begin{figure}[htb!]
\begin{subfigure}[b]{\textwidth}
         \includegraphics[width=\textwidth]{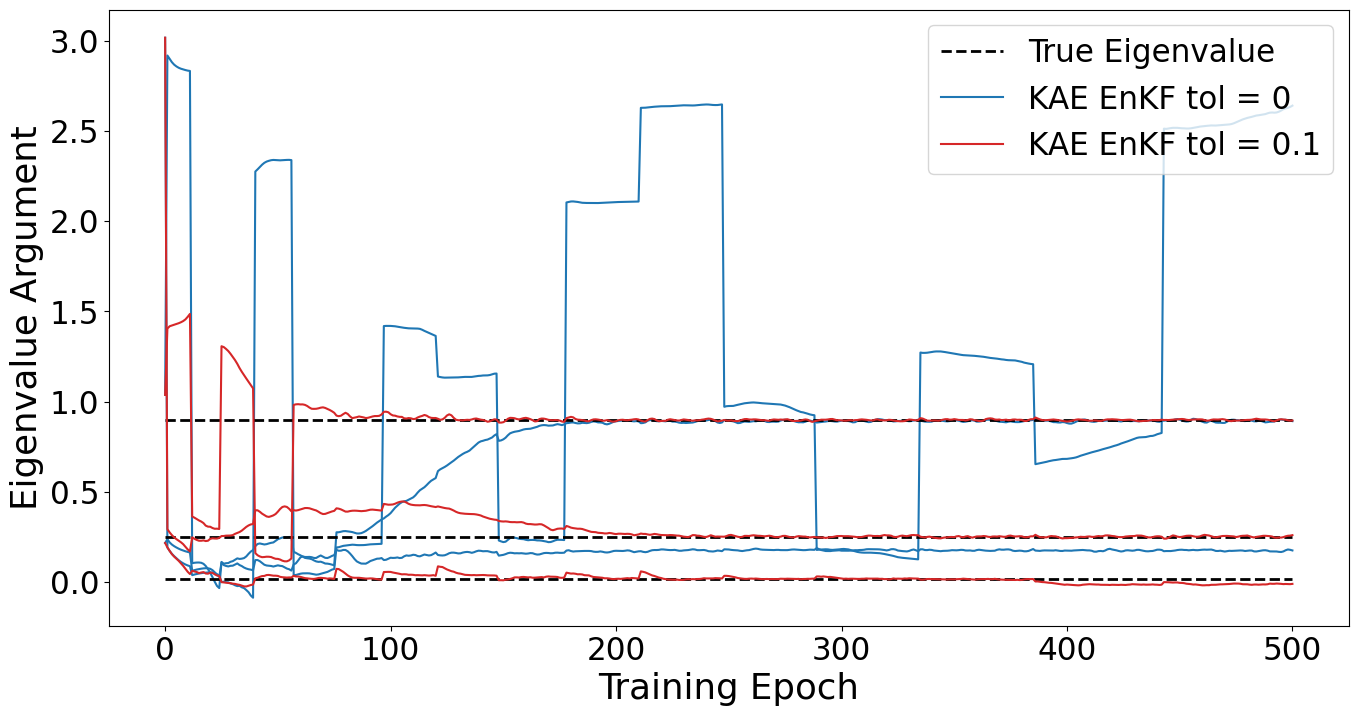}
     \end{subfigure}
     \caption{\centering\label{fig:multifreq_kae_training_freqs}The Koopman autoencoder model’s estimated eigenvalue arguments, over the course of training for 500 epochs on the spin-up data, with a multi-frequency tolerance enabled of $(\text{tol}=0.1)$, and without a multi-frequency tolerance $(\text{tol}=0)$. We mark the system’s true, non-stationary eigenvalue argument as its average value over the course of the training data. The Koopman autoencoder with a suitable multi-frequency tolerance is able to identify all the system’s eigenvalue arguments correctly within 200 training epochs, however the Koopman autoencoder with no multi-frequency tolerance produces eigenvalue argument estimates that cluster closely together during the early training epochs. This permanently biases the encoder/decoder networks which are being trained simultaneously with the model’s eigenvalue estimates, resulting in the Koopman autoencoder with $\text{tol}=0$ failing to identify all the system’s frequencies by the end of the training process.}
\end{figure}

The KAE EnKF was trained for more epochs than in the previous synthetic tests and all other parameters were kept the same, however as the data being fitted contains multiple frequencies, the Koopman autoencoder's multi-frequency tolerance parameter is now critically important. When training the Koopman autoencoder during the spin-up phase, it is possible for frequencies of the model's Koopman approximator $\mathbf{K}_\lambda$ to converge to the same frequency, and hence miss other less dominant frequencies also present within the data. To address this, the user sets a multi-frequency tolerance, which specifies a minimum distance from the model's pre-existing frequencies that newly identified frequencies during the global frequency optimization step must adhere to. Examples of how the Koopman autoencoder model's frequencies develop over the 500 training epochs of the spin-up phase, for multi-frequency tolerances of 0 and 0.1 respectively, are shown in Figure \ref{fig:multifreq_kae_training_freqs}.

\begin{figure}[htbp]
\begin{subfigure}[b]{\textwidth}
         \includegraphics[width=\textwidth]{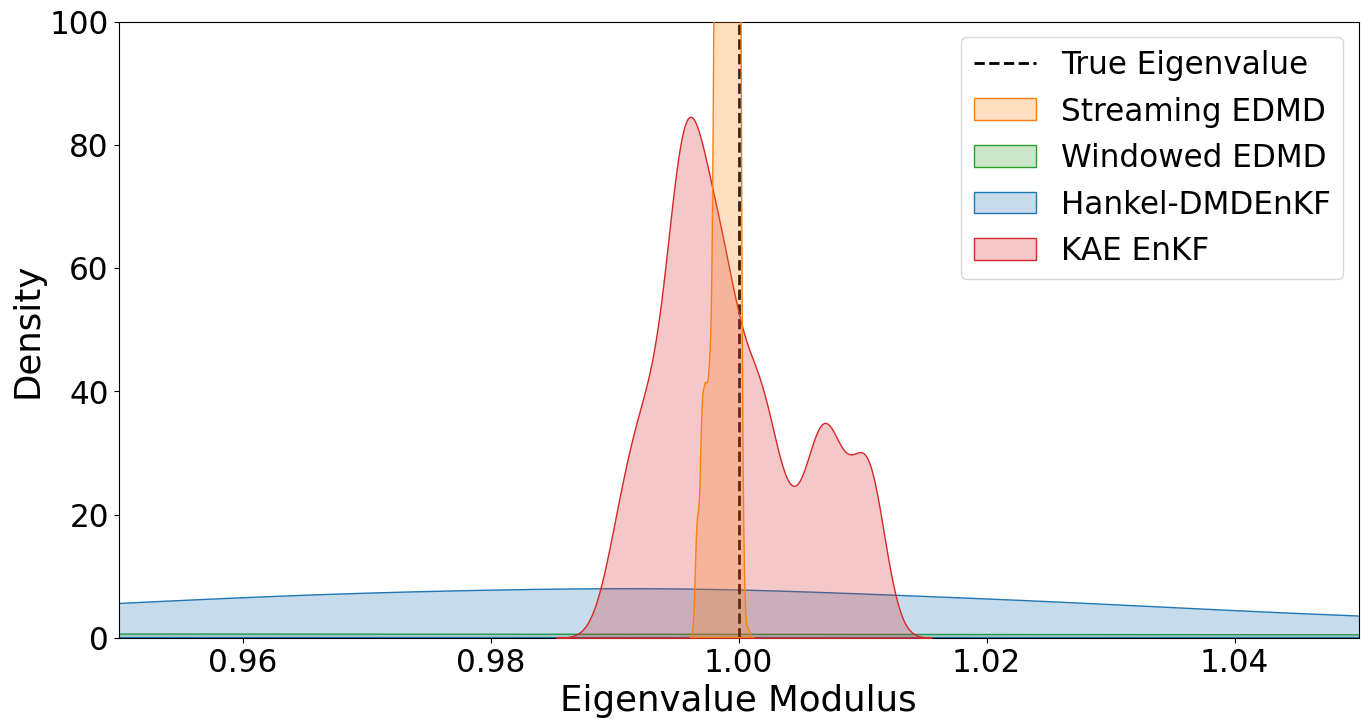}
     \end{subfigure}
     \caption{\centering\label{fig:multifreq_edmd_eigmod_dist}The combined distribution of all 3 of the  system's eigenvalue's moduli, estimated by each of the nonlinear iterative DMD variants over an example trajectory of the data. Combining the distributions is a  reasonable choice, as all eigenvalues have a constant modulus of 1 throughout the experiment. Streaming EDMD produces by far the most accurate modulus estimate due to its enhanced stability, with KAE EnKF also producing tight, relatively accurate estimates. The Hankel-DMDEnKF and windowed EDMD produce wider distributions, with windowed EDMD performing significantly worse than the Hankel-DMDEnKF, due to the highly unstable nature of the method.}
\end{figure}

Both of the KAE EnKF's experience occasional sharp changes in their estimated frequencies at epochs where the global frequency optimization algorithm is applied, which is to be expected as the stochastic gradient descent optimizer acts instead by taking small iterative steps towards the loss's minimizer. The KAE EnKF with a multi-frequency tolerance of 0.1 initially jumps around with its frequency estimates, due to factors like sample variance between batches, however by epoch 200 is able to correctly identify all three frequencies present within the data. Alternatively, the KAE EnKF with a multi-frequency tolerance of 0 is only able to successfully identify one of the system's frequencies by the end of the 500 epochs. Early on in the training process, all of this model's approximated frequencies group together, due to the lack of a sufficient multi-frequency tolerance forcing them apart. In other tests this significantly slowed training, however the correct frequencies were eventually able to be identified. In this case however, the encoder/decoder networks, which are simultaneously being fitted to the data, were too heavily biased by the initial frequency clustering to facilitate discovery of the system's true frequencies within the allocated training time. Given these results from Figure \ref{fig:multifreq_kae_training_freqs}, we set the KAE EnKF's multi-frequency tolerance to 0.1, and begin testing its ability to track the system's eigenvalue modulus and arguments against those of the previous nonlinear, iterative DMD variants. All computational parameters for streaming EDMD, windowed EDMD and the Hankel-DMDEnKF were kept the same as before, with the model's ranks $r$ being raised from 2 to 6, to accommodate the latent space's increased dimension. Each method's performance when tracking the system's 3 eigenvalue moduli can be viewed in Figure \ref{fig:multifreq_edmd_eigmod_dist}.

\begin{figure}[!htb]
\begin{subfigure}[b]{.49\textwidth}
         \includegraphics[width=\textwidth]{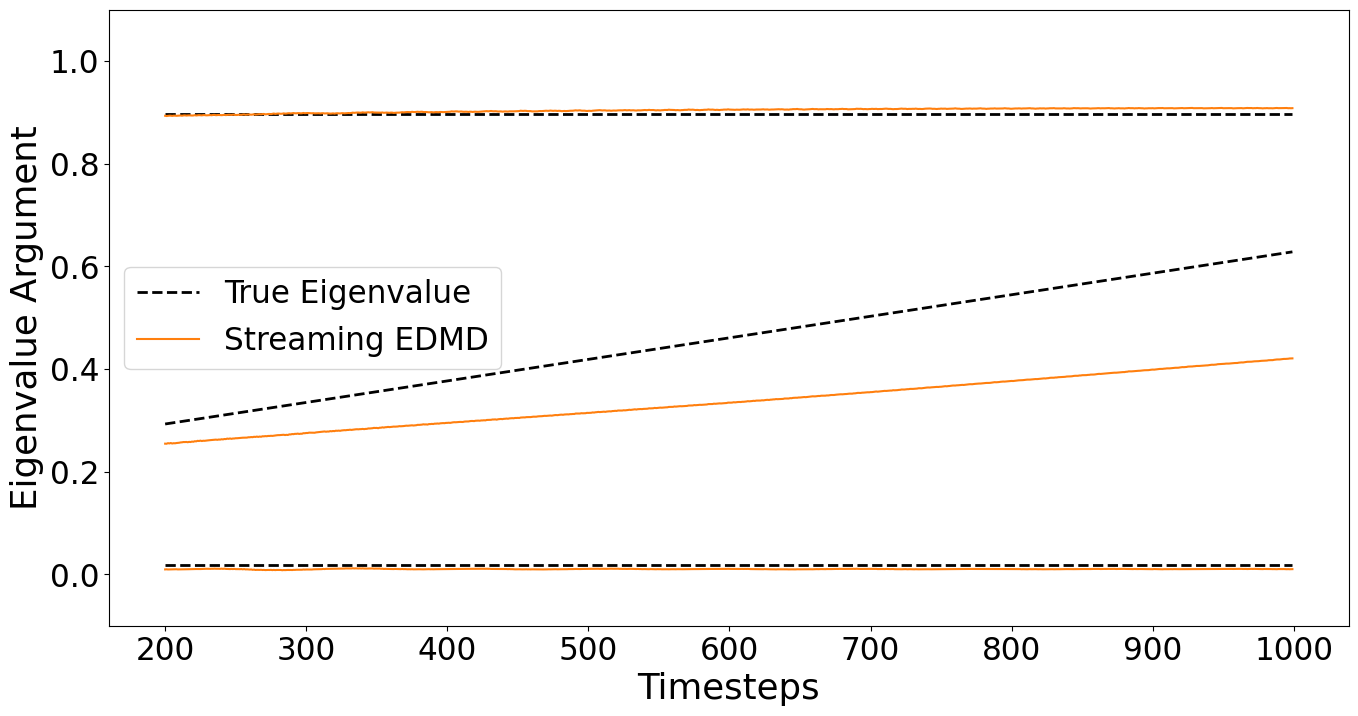}
                \caption{\centering\label{fig:multifreq_sedmd_eigarg_trajectory}Streaming EDMD's eigenvalue argument estimates.}
\end{subfigure}
\hfill     
\begin{subfigure}[b]{.49\textwidth}
         \includegraphics[width=\textwidth]{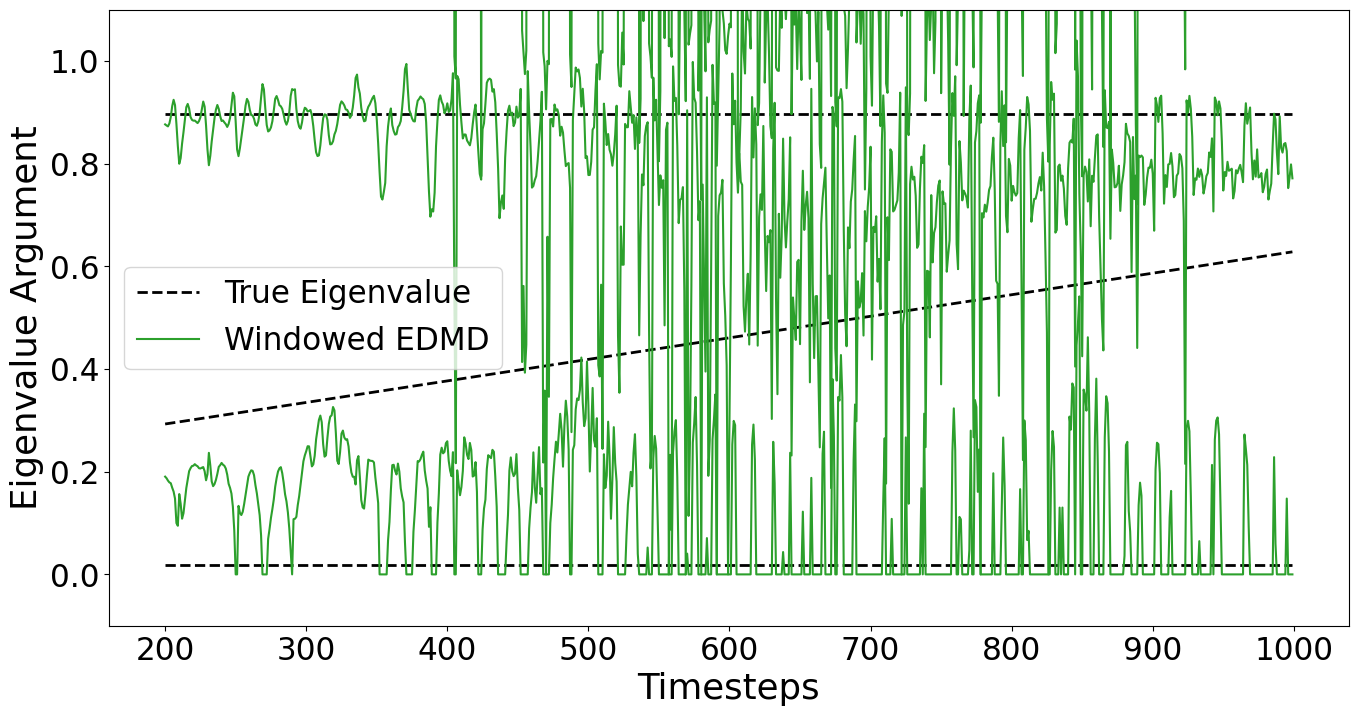}
                \caption{\centering\label{fig:multifreq_wedmd_eigarg_trajectory}Windowed EDMD's eigenvalue argument estimates.}
\end{subfigure}

\begin{subfigure}[b]{.49\textwidth}
         \includegraphics[width=\textwidth]{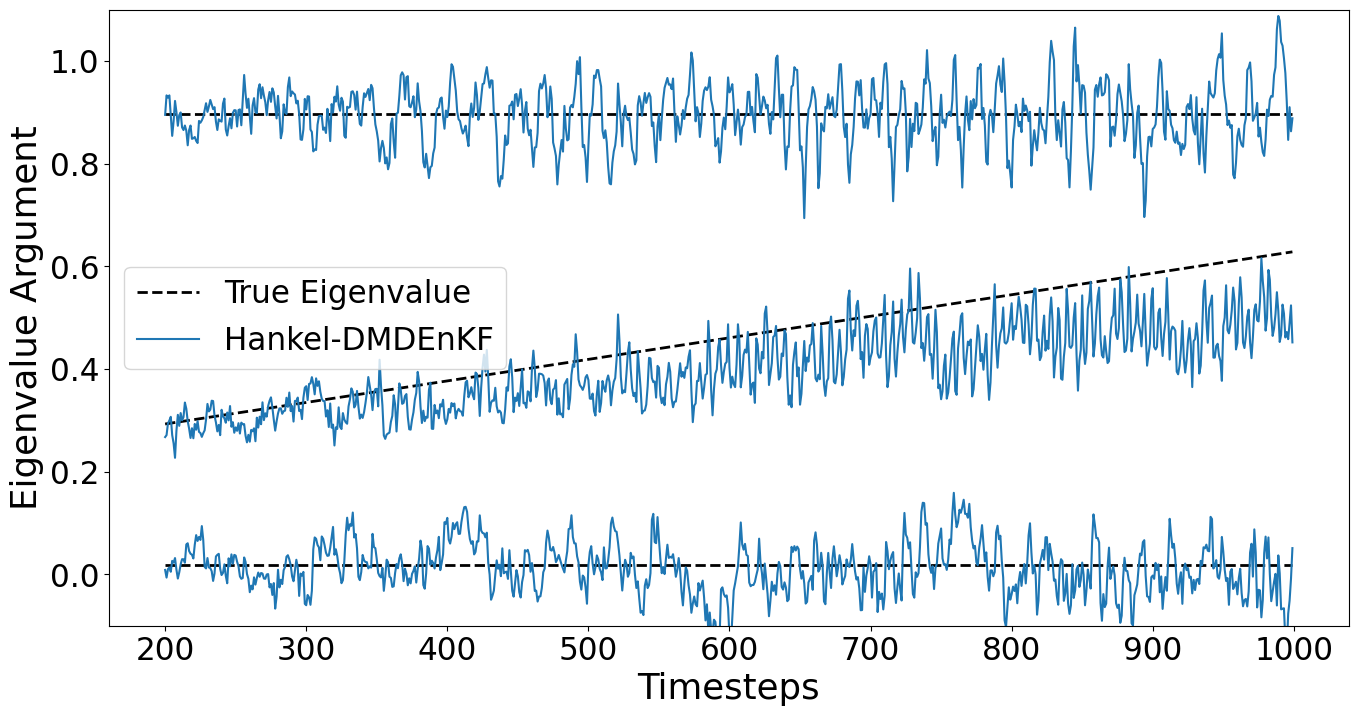}
                \caption{\centering\label{fig:multifreq_dmdenkf_eigarg_trajectory}Hankel-DMDEnKF's eigenvalue argument estimates.}
\end{subfigure}
\hfill     
\begin{subfigure}[b]{.49\textwidth}
         \includegraphics[width=\textwidth]{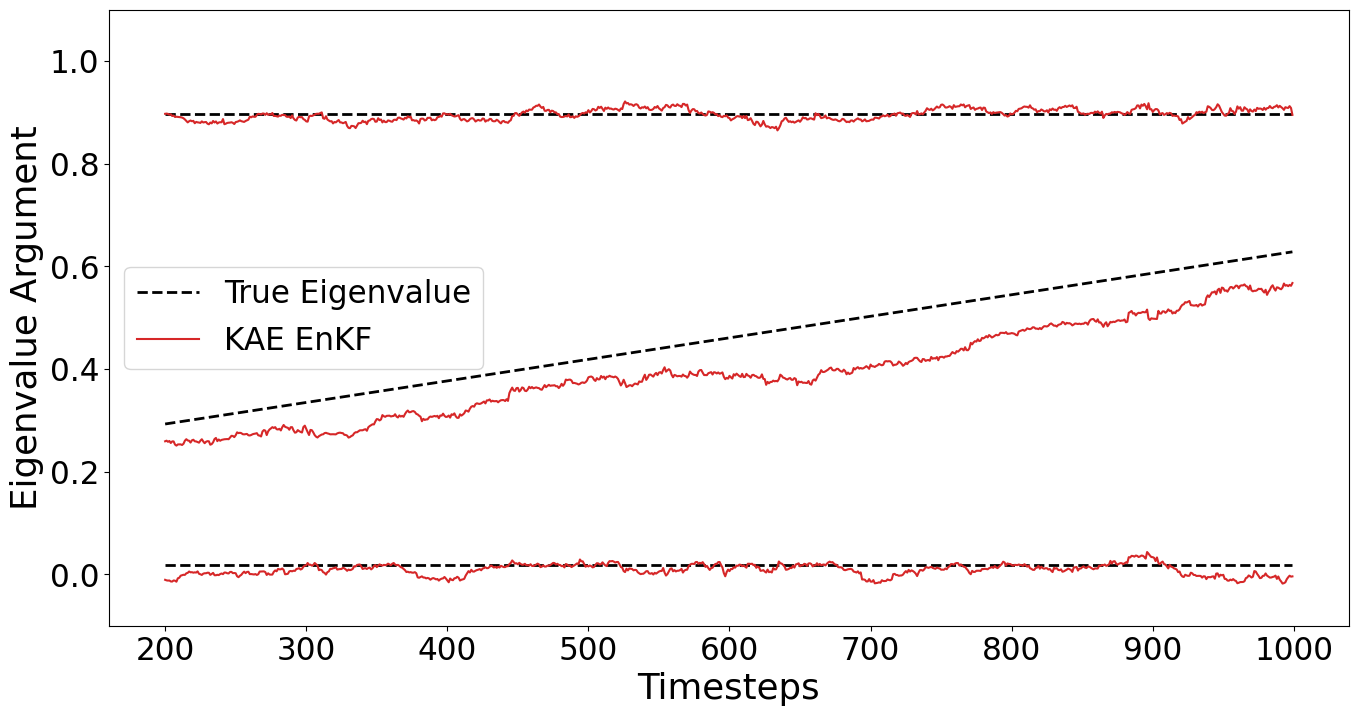}
                \caption{\centering\label{fig:multifreq_kaeenkf_eigarg_trajectory}The KAE EnKF's eigenvalue argument estimates.}
\end{subfigure}

     \caption{\centering\label{fig:multifreq_edmd_eigarg_trajectory}
     The eigenvalue argument estimates by each of the nonlinear iterative DMD variants over an example trajectory of the data. Streaming EDMD (top left) produces stable estimates, with a significant lag when tracking the linearly increasing eigenvalue argument. Windowed EDMD (top right) does not correctly capture all 3 eigenvalue’s arguments, and for those it does its estimates are highly unstable. The Hankel-DMDEnKF (bottom left) produces estimates with only a minor bias when tracking the system’s non-stationary eigenvalue, and experiences some small general estimates instability. The KAE EnKF (bottom right) tracks the system’s stationary parameters with a slightly larger level of instability than streaming EDMD, however produces significantly less lagged estimates than it for the increasing eigenvalue argument.}
\end{figure}

All three of the synthetic system's eigenvalues have modulus 1 throughout the experiment, hence we plot the moduli of all the system's eigenvalues estimated by each method together as a histogram in Figure \ref{fig:multifreq_edmd_eigmod_dist}. As in the previous experiments, streaming EDMD produces the tightest distribution around the true system's eigenvalue modulus, due to its superior stability over the other methods when estimating stationary parameters. The KAE EnKF is the next most accurate, as the eigenvalue stabilizing term in the Koopman autoencoder's loss function, alongside the distinct eigenvalue modulus and argument uncertainty parameters inside the KAE EnKF's filter, lead to enhanced stability over the Hankel-DMDEnKF. The worst performing method is windowed EDMD, with its recency prioritizing window leaving it highly susceptible to destabilization by small, random fluctuations in the data. We now analyse each method's ability to track the two stationary and 1 non-stationary eigenvalue arguments of the systems, shown graphically in Figure \ref{fig:multifreq_edmd_eigarg_trajectory}.

Streaming EDMD again produces very stable estimates of the system's eigenvalues, as seen in Figure \ref{fig:multifreq_sedmd_eigarg_trajectory}, however at the cost of a significant lag in the estimation of the non-stationary eigenvalue argument. Only two of the system's frequencies are ever successfully identified, and as the experiment progresses and the system's fastest two frequencies approach one another, the model's corresponding frequency estimate is forced into the interval between them. The Hankel-DMDEnKF is significantly more stable than windowed EDMD, however exhibits more instability than streaming EDMD or the KAE EnKF in its eigenvalue argument estimates displayed in Figure \ref{fig:multifreq_dmdenkf_eigarg_trajectory}. The Hankel-DMDEnKF's estimates of the two stationary eigenvalue arguments are unbiased, and for the non-stationary eigenvalue argument its estimates contain the least lag for any of the methods. Figure \ref{fig:multifreq_kaeenkf_eigarg_trajectory} shows the KAE EnKF's highly stable estimates for the system's eigenvalue arguments. Similarly to the Hankel-DMDEnKF, the KAE EnKF produces unbiased estimates for the system's stationary parameters, and only a small lag in its estimates of the linearly increasing frequency.

\begin{figure}[htbp]
\begin{subfigure}[b]{\textwidth}
         \includegraphics[width=\textwidth]{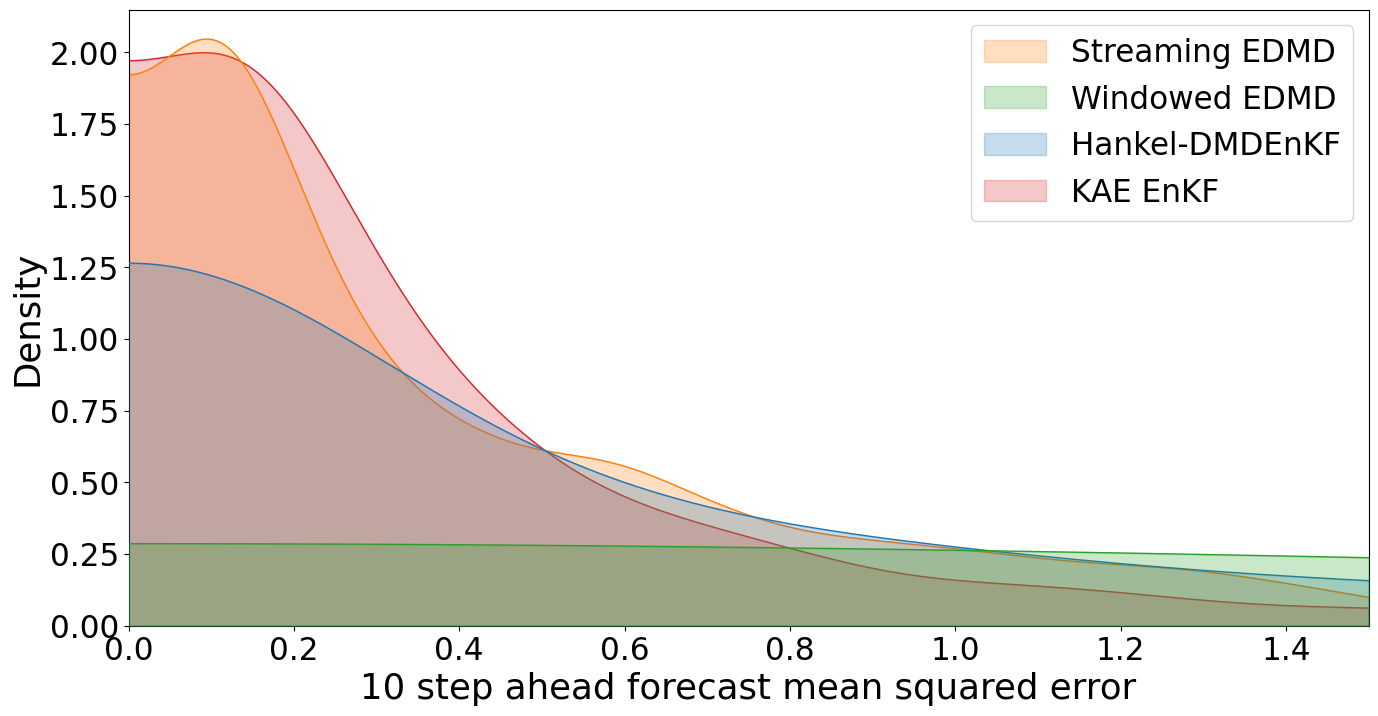}
     \end{subfigure}
     \caption{\centering\label{fig:multifreq_edmd_forecast_mse_dist}The 10-step ahead forecast mean squared error distribution for each of the nonlinear iterative DMD variants, estimated over the course of 10 runs. Streaming EDMD produces a tight error distribution, however peaks a small distance away from 0 due to the irreducible error its lag when tracking the system’s non-stationary parameter introduces. The KAE EnKF’s errors peak slightly closer to 0, and its distribution has a similar width to that of streaming EDMD. The Hankel-DMDEnKF also performs well, with its error distribution peaking at 0, although having a greater spread than for the previous two methods. Windowed EDMD is the worst performing method, with an incredibly wide error distribution.}
\end{figure}

We now evaluate the efficacy of each method when generating 10-step ahead forecast for the system's full state, with mean squared error distributions calculated over 10 runs and displayed in Figure \ref{fig:multifreq_edmd_forecast_mse_dist}. Here, the lack of stability in Windowed EDMD and to a lesser degree the Hankel-DMDEnKF as visible in Figures \ref{fig:multifreq_edmd_eigmod_dist} and \ref{fig:multifreq_edmd_eigarg_trajectory}, manifests as correspondingly wider forecast error distributions. For the same reason, streaming EDMD and the KAE EnKF have a similar spread in their forecast error distributions. A major differentiating feature however, is that the KAE EnKF and Hankel-DMDEnKF's distributions peak at a lower error value than that of streaming EDMD. This is caused by streaming EDMD's significantly lagging estimate of the non-stationary eigenvalue argument seen in \ref{fig:multifreq_sedmd_eigarg_trajectory}, which induces an irreducible error in the model's subsequent forecast, driving its error distribution away from 0. The KAE EnKF and Hankel-DMDEnKF suffer only a minor bias when tracking the time varying parameter, as shown in Figures \ref{fig:multifreq_kaeenkf_eigarg_trajectory} and \ref{fig:multifreq_dmdenkf_eigarg_trajectory}, hence the distribution of their forecast errors peak at a lower value. The KAE EnKF's tight error distribution, which peaks close to 0, indicates that when seeking stable, unbiased forecasts for nonlinear, non-stationary systems that contain multiple latent frequencies, of the models tested, the KAE EnKF is the preferred method.

\section{Video of a pendulum application}\label{s:pendulum}
\subsection{Problem setup}

\begin{figure}[htbp]
\begin{subfigure}[b]{\textwidth}
         \includegraphics[width=\textwidth]{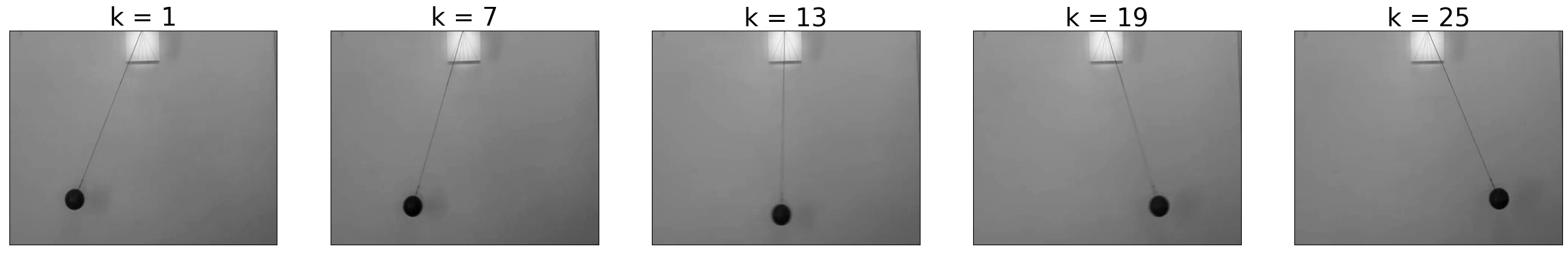}
     \end{subfigure}
     \caption{\centering\label{fig:pendulum_example_data}Image frames from the real video of a pendulum, showing the pendulum going through approximately half a full oscillation over the course of 24 time steps.}
\end{figure}

We take raw video footage of a pendulum in motion \footnote{ScienceOnline: The Pendulum and Galileo (\href{www.youtube.com/watch?v=MpzaCCbX-z4}{www.youtube.com/watch?v=MpzaCCbX-z4}).}, as shown in Figure \ref{fig:pendulum_example_data}, and use the video's frames to generate a discrete dynamical system from the pendulum's continuous dynamics. By taking the video's frames to be measurements of the full state, we thus embed a latent system with a single complex conjugate pair of eigenvalues, into a high-dimensional, nonlinear pixel space. Identifying the eigenvalues of latent systems from high-dimensional and nonlinear full state measurements, for the purposes of tracking and forecasting, is precisely the use case the KAE EnKF was designed for. Hence, we evaluate the efficacy of the KAE EnKF and compare it with that of EDMD, when applied directly to the image data visible in Figure \ref{fig:pendulum_example_data}.

The pendulum video dataset comprises 523 image frames, each consisting of $576\times720$ pixels, taken from a stationary camera, centred on a pendulum swinging uninterrupted inside a lab. As both the KAE EnKF and EDMD take a vector as their input, each image matrix is rearranged in lexicographical ordering to produce a 414720 dimensional vector. To ensure relevant information about the pendulums velocity and acceleration is encoded into the state, the flattened vectors of the frames from the preceding two time steps are then stacked underneath the data for current time step $k$. This results in a dataset comprising 521 full state vectors $\mathbf{x}_k \in{\mathbb R^{1244160}}$, of which we use the first 423 time steps to train the models, and reserve the subsequent data points for testing purposes. The only other preprocessing we apply to the pendulum video frames, is centring the data via the removal of the states in the training set's mean value. The pendulum completes one full oscillation approximately every 48 frames and experiences minimal visible damping, hence the corresponding latent states $\Tilde{\mathbf{x}}_k\in{\mathbb R^{2}}$ evolve approximately under the action of a dynamical system with constant complex conjugate eigenvalue pair $\lambda_\pm = 1e^{\pm (2\pi/48) i}$.

\subsection{Reconstructing the training data with the Koopman autoencoder and EDMD}

The training set is used to fit EDMD, with the full state measurements $\mathbf{x}_k$ input directly as the algorithm's dictionary of observables, for they are already time-delay embedded as specified in their above construction. This data is also used to train the Koopman autoencoder model during the KAE EnKF's spin-up phase. To train the Koopman autoencoder model, Bayesian hyperparameter optimization techniques \cite{practical_bayesian_hyperparameter_optimisation} are employed, using the Python Hyperopt library \cite{hyperopt_library}, to optimize the algorithm's computational parameters such as learning rate and how much each individual penalty terms contributes towards the overall loss function. The encoder and decoder networks $\mathbf{h}_\phi$ and $\Tilde{\mathbf{h}}_{\Tilde{\phi}}$ are also widened by the doubling of the hidden layer's input dimension, compared with the network architecture employed in the synthetic applications, to accommodate the stronger nonlinearity present in the real application's video data. 

\begin{figure}[ht]
\begin{subfigure}[b]{\textwidth}
         \includegraphics[width=\textwidth]{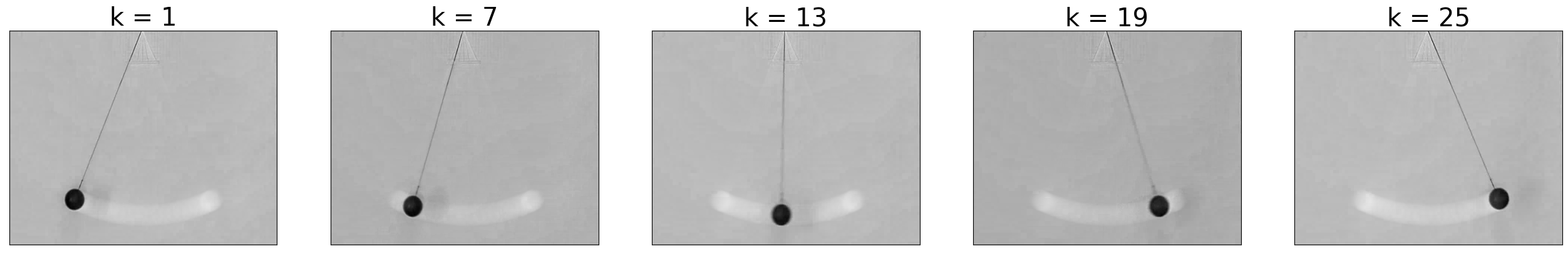}
         \caption{\centering\label{fig:pendulum_example_data_mean_removed}Frames from the video of a pendulum, with the mean, calculated over the training data, removed.}
     \end{subfigure}

\begin{subfigure}[b]{\textwidth}
         \includegraphics[width=\textwidth]{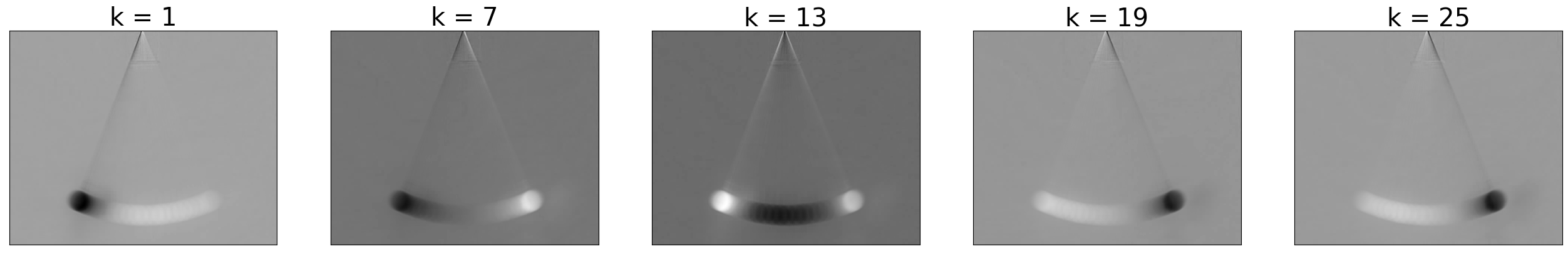}
         \caption{\centering\label{fig:pendulum_edmd_encode_decode}Reconstructions of the pendulum images, generated by mapping frames to and then back from the EDMD model’s latent space.}
     \end{subfigure}

\begin{subfigure}[b]{\textwidth}
         \includegraphics[width=\textwidth]{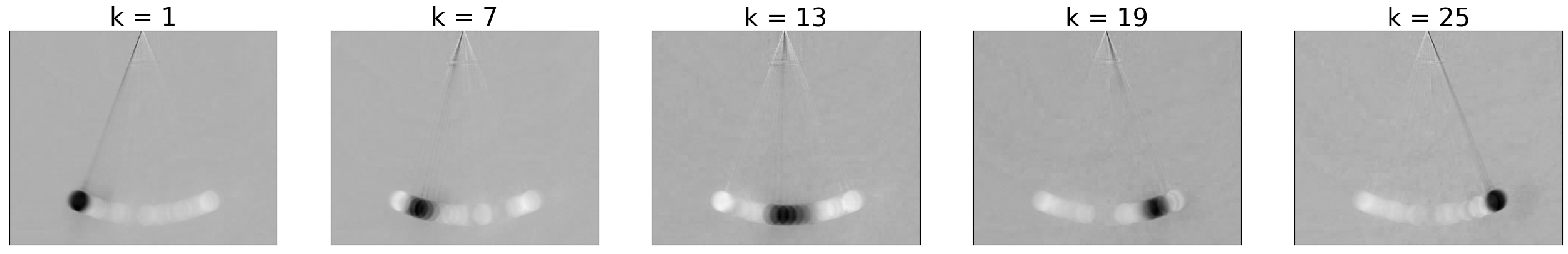}
         \caption{\centering\label{fig:pendulum_kae_enkf_encode_decode}Reconstructions of the pendulum images, generated by mapping frames to and then back from the Koopman autoencoder’s latent space.}
     \end{subfigure}

     \caption{\centering\label{fig:pendulum_encode_decode}Frames from the model’s centred input data (top), reconstructed my mapping the inputs to each model’s latent space, before mapping the latent states back to the full state space. EDMD (middle) produces accurate reconstructions for points at the ends of the pendulums arc, however generates inaccurate estimates of the pendulums position otherwise, with significant blurring at points in the trajectory where the pendulum is moving more quickly. The Koopman autoencoder (bottom) produces more accurate overall reconstructions, although also experiences some minor blurring in reconstructions as the pendulum's velocity increases.}
\end{figure}

Both models are trained on the centred full state data $\mathbf{x}_k$ for $k=(1,...,423)$, and before we begin analysing each model's performance tracking and forecasting previously unseen data, we first review their efficacy in reconstructing the training data. Each model can be broken down into two essential parts, a stationary mapping back and forth between the full state space and the model's latent space, and the estimated latent dynamical system where all the model's dynamics take place. For the Koopman autoencoder, the mappings between full and latent spaces are simply the model's encoder/decoder networks described in equation \eqref{eqn:gxk=tildexk}. For EDMD, the DMD modes $\mathbf{\Phi}$ act as a mapping from latent space back to the full state space, and their pseudoinverse $\mathbf{\Phi}^+$ send full states to their respective latent states. Each model's latent, linear dynamics are available trivially, corresponding to the Koopman approximator $\mathbf{K}_\lambda$ from equation \eqref{eqn:xtk+1=klxtk} in the Koopman autoencoder framework, and in the diagonal matrix of eigenvalues calculated when performing EDMD. We evaluate each method's ability over the training data to map full states into the model's latent space and back only, without applying the model's estimated latent dynamics, and the results are displayed in Figure \ref{fig:pendulum_encode_decode}.

Figure \ref{fig:pendulum_example_data_mean_removed} displays centred image frames from the training data, as it completes half an oscillation. EDMD's attempted reconstruction of this data is shown in Figure \ref{fig:pendulum_edmd_encode_decode}. At times $k=1$ and $k=25$, when the pendulum is at the cusp of its swing, EDMD produces fairly accurate image reconstructions. This is likely due to the pendulums reduced velocity at this stage in its oscillation, resulting in a greater confidence in the model's estimated position of the pendulum and hence a more stable reconstruction. As the pendulum's velocity increases at times $k=7$ and $k=19$, we see EDMD's estimation of the pendulum's position begin to become more blurred, with a slight bias towards the region the pendulum occupied during previous time steps. By $k = 13$, when the pendulum is travelling at its maximal velocity, the pendulum in the reconstruction is highly blurred, encompassing almost the entire pendulum's trajectory, and excluding only the oscillations end points. The varying background colour between reconstructed frames also indicates that EDMD's latent state mapping is not only capturing dynamically relevant information about the system, as the video's stationary background pixels should remain constant and uninterfered with when transitioning through the model's latent space. Figure \ref{fig:pendulum_kae_enkf_encode_decode} demonstrates the Koopman autoencoder's reconstructions, which are significantly more accurate and sharper than those produced by EDMD. The model suffers from the same blurring effect as EDMD when the pendulums speed increases, however to a much lesser degree, and the pixels corresponding to the frame's background remain stationary as intended.

We now inspect the Koopman autoencoder and EDMD's estimates of the system's dynamics in the latent space. Both models correctly identify a complex conjugate pair of eigenvalues with modulus very close to $1$, EDMD's modulus errors were of order $\mathcal{O}(10^{-3})$, whereas the Koopman autoencoder achieved errors an order of magnitude lower at $\mathcal{O}(10^{-4})$. The system's eigenvalue argument $2\pi/48 = 0.131$ was accurately estimated by the Koopman autoencoder at 0.133 with an order of error of $\mathcal{O}(10^{-3})$, whereas EDMD calculated an eigenvalue argument within order $\mathcal{O}(10^{-3})$ of 0 at 0.002. This results in an EDMD model with very slow dynamics, which when generating future forecasts, essentially predicts the system's current state as its new future state.

\begin{figure}[!htbp]
\begin{subfigure}[b]{.49\textwidth}
         \includegraphics[width=\textwidth]{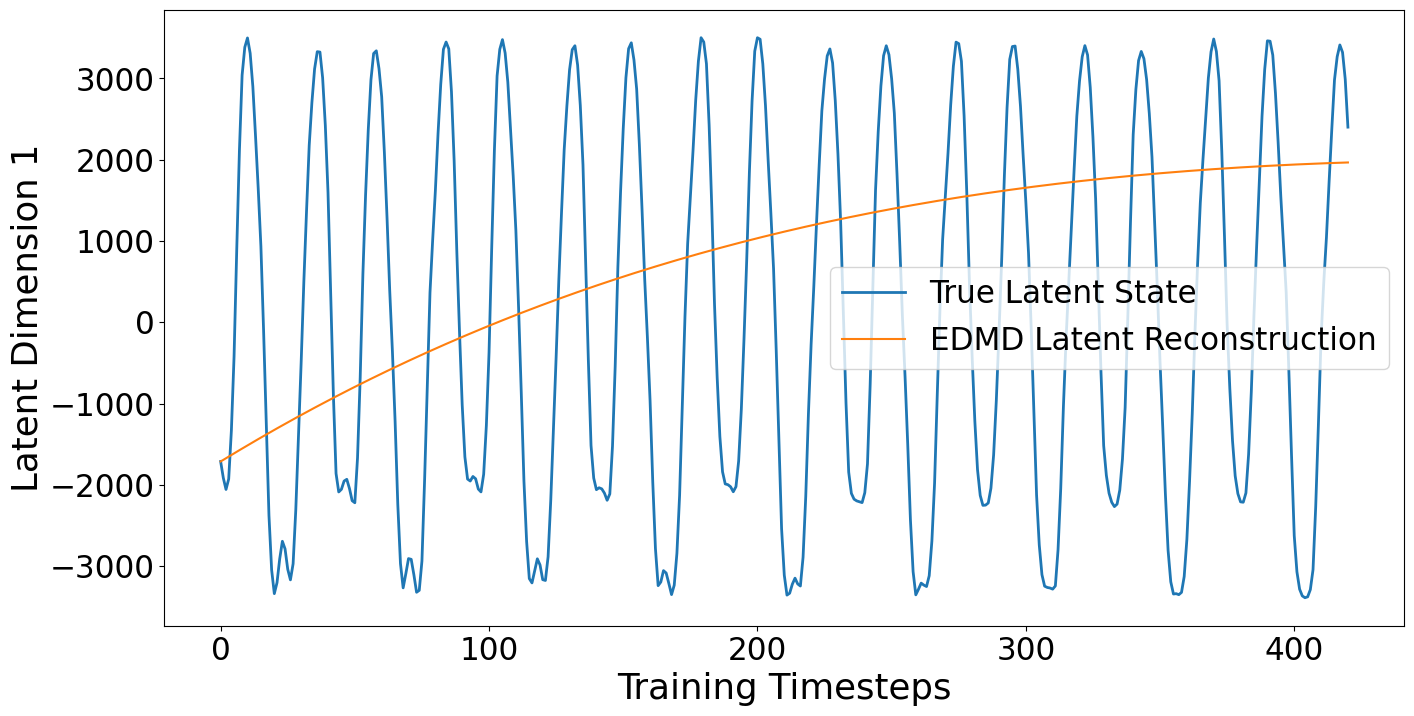}
                \caption{\centering\label{fig:pendulum_edmd_latent_reconstruction}The system's true states from the training set, mapped into EDMD's latent space, and reconstructed using its latent dynamics.}
\end{subfigure}
\hfill     
\begin{subfigure}[b]{.49\textwidth}
         \includegraphics[width=\textwidth]{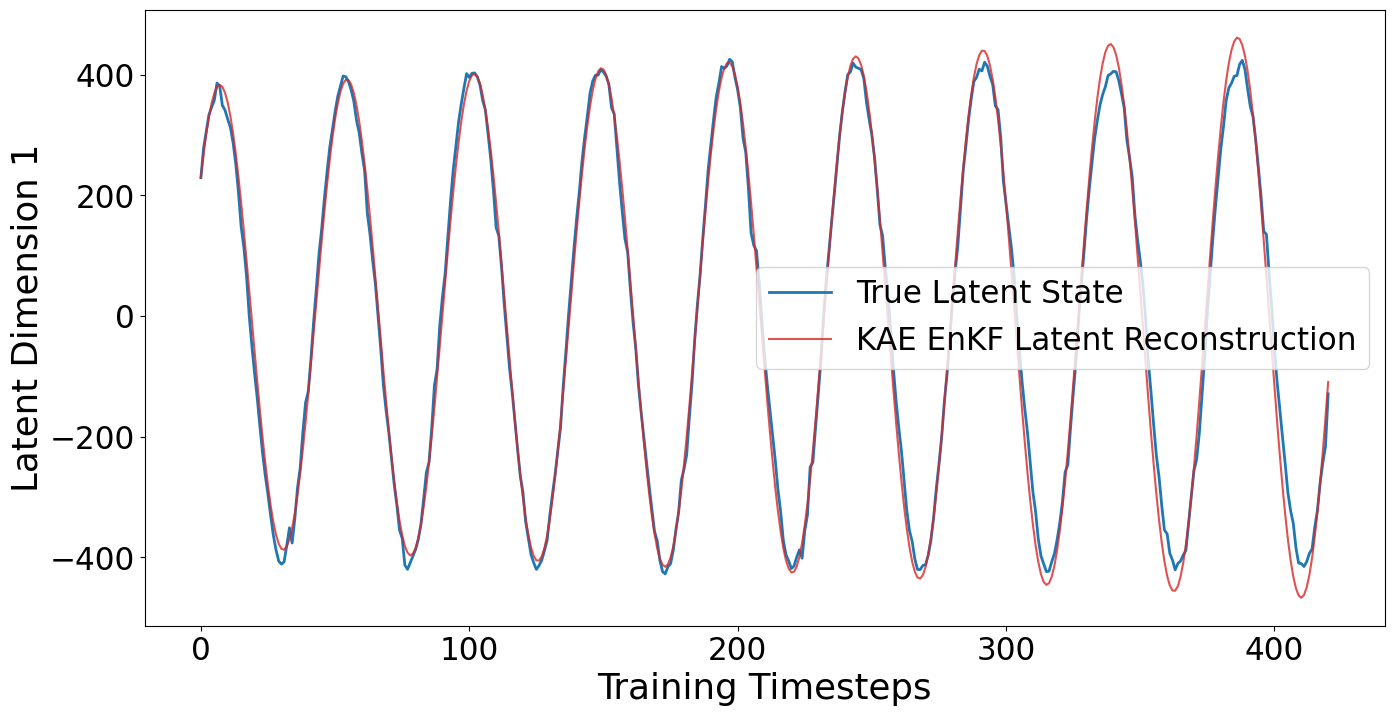}
                \caption{\centering\label{fig:pendulum_kaeenkf_latent_reconstruction}The system's true states from the training set, mapped into Koopman autoencoder's latent space, and reconstructed using its latent dynamics.}
\end{subfigure}

     \caption{\centering\label{fig:pendulum_latent_reconstruction}The system's true states from the training set, mapped into the model's latent space, and reconstructed using the model's latent dynamics, for EDMD (left), and the KAE EnKF (right). EDMD's latent representation of the data differs in multiple way from the shape of a pure sine wave, indicating its mapping to a supposedly linear latent space is flawed. The model's latent dynamics are significantly slower than the true dynamics, hence the dynamic reconstruction quickly diverges from the true system's latent states. The Koopman autoencoder produces an almost exactly sinusoidal representation of the true latent data, and the dynamic reconstruction captures this with a high level of accuracy.}
\end{figure}

We map the training data into each model's latent space, then run their latent dynamics forward from $\Tilde{\mathbf{x}}_1$ continuously, to generate reconstructions of the training data in each model's latent space. EDMD's latent states are a complex conjugate pair, due to its complex conjugate DMD modes/eigenvalues, however this eigenstructure corresponds to latent dynamics of a rotation in the basis $(\operatorname{Re}(\mathbf{\Phi}),\operatorname{Im}(\mathbf{\Phi}))$. Hence, to visualize EDMD's latent state evolution, we project the full state $\mathbf{x}_k$ onto normalized basis vectors $(\operatorname{Re}(\mathbf{\Phi})/\|\operatorname{Re}(\mathbf{\Phi})\|,\operatorname{Im}(\mathbf{\Phi})/\|\operatorname{Im}(\mathbf{\Phi})\|)$, to generate latent states $\Tilde{\mathbf{x}}_k\in{\mathbb R^{2}}$. The Koopman autoencoder and EDMD's latent embedding of the training data, alongside each model's dynamic reconstruction of the data in latent space, can be seen in Figure \ref{fig:pendulum_latent_reconstruction}.

EDMD is unable to find a fully linear latent representation of the data, as is visible in Figure \ref{fig:pendulum_edmd_latent_reconstruction}. The system's true full state is mapped to an oscillating latent state, however the time evolution of this state is not exactly sinusoidal, containing artefacts such disproportionately sharp peaks, and periodically varying amplitudes. The Koopman autoencoder on the other hand is able to find a nonlinear mapping that sends the system's full states to a latent space, within which they evolve almost exactly linearly, as demonstrated by the near perfect sine curve for the true latent states in Figure \ref{fig:pendulum_kaeenkf_latent_reconstruction}. The Koopman autoencoder's enhanced ability to identify such a mapping, is likely due to the increased flexibility of its neural network based encoder/decoders, over the DMD mode identifying approach of EDMD. When viewing the model's estimated latent dynamics propagated forward in time from $\Tilde{\mathbf{x}}_1$ for EDMD in Figure \ref{fig:pendulum_edmd_latent_reconstruction}, we see the dynamic reconstruction evolves significantly slower and quickly diverges from the system's true latent states. This is to be expected, due to the value of EDMD's complex conjugate eigenvalue's argument being much closer to 0 than the true system's frequency. The Koopman autoencoder's dynamic reconstruction in Figure \ref{fig:pendulum_kaeenkf_latent_reconstruction} matches the true latent states exceptionally closely, with slight deviations between the reconstruction and ground truth developing as the forecast horizon increases.

Having established that EDMD generates poor full to latent state mappings, and an inaccurate model of the system's latent dynamics, we now proceed only with the Koopman autoencoder, investigating how the KAE EnKF performs when tracking and forecasting previously unseen frames from the pendulum video dataset.

\subsection{Tracking and forecasting unseen data}

\begin{figure}[htb]
\begin{subfigure}[b]{\textwidth}
         \includegraphics[width=\textwidth]{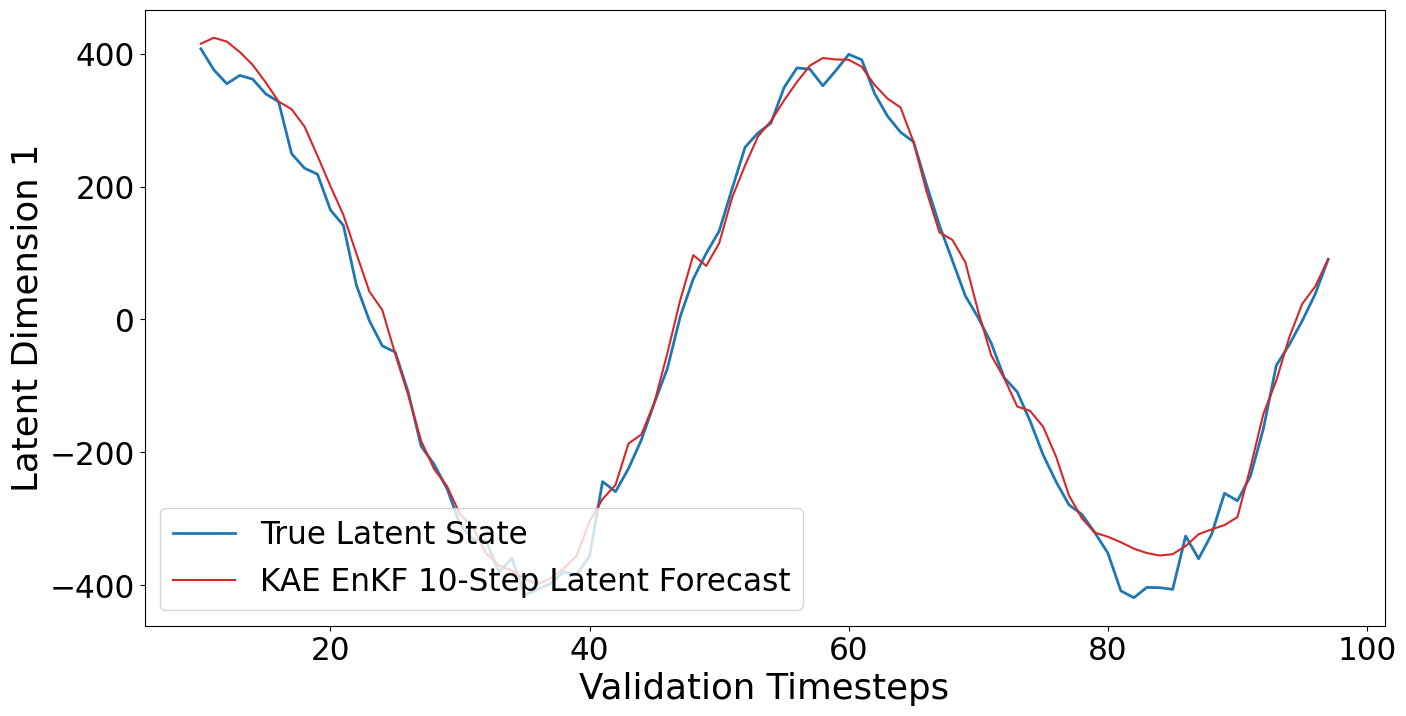}
     \end{subfigure}
     \caption{\centering\label{fig:pendulum_kaeenkf_latent_forecast}The KAE EnKF's 10-step ahead forecast, in the first dimension of the model's latent space, over the previously unseen testing data. The Koopman autoencoder's forecasts are relatively smooth, and match the true data well.}
\end{figure}

We apply the KAE EnKF’s filtering step iteratively to the remaining 98 data points in the testing set. At each time step, we propagate the latent ensemble 1 step forward using the Koopman autoencoder model, then assimilate the new data point and update our ensemble accordingly. We then generate a 10-step ahead forecast from the current time step, applying the Koopman autoencoder to the updated ensemble, and take the ensemble’s mean as a point estimate for the 10-step ahead forecasted latent state. These forecasts of the system’s latent state, performed over the testing set, are shown in Figure \ref{fig:pendulum_kaeenkf_latent_forecast}.

\begin{figure}[htb]
\begin{subfigure}[b]{\textwidth}
         \includegraphics[width=\textwidth]{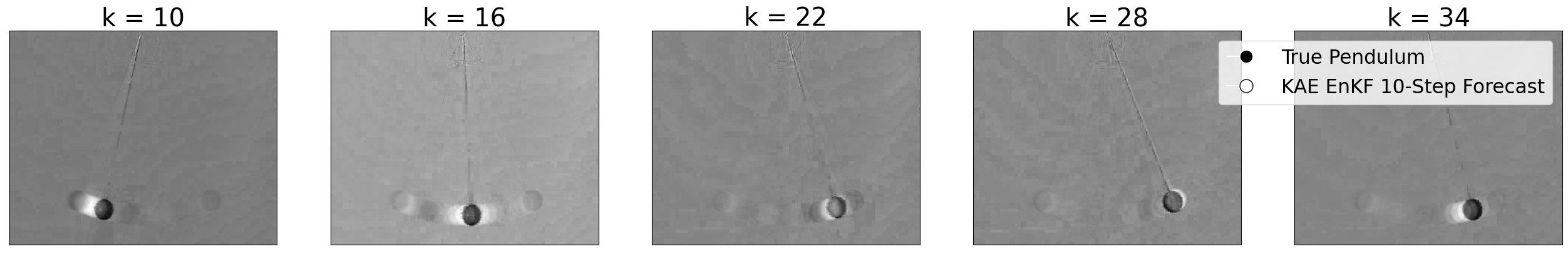}
     \end{subfigure}
     \caption{\centering\label{fig:pendulum_kae_enkf_val_forecast_subtracted}The difference between the true pendulum data in the test set, and the KAE EnKF's corresponding 10-step ahead full state forecasts. The true pendulum is represented by a black circle, and the forecast, by virtue of being subtracted from the true data, a white circle.}
\end{figure}

Figure \ref{fig:pendulum_kaeenkf_latent_forecast} shows a more zoomed in look than Figure \ref{fig:pendulum_kaeenkf_latent_reconstruction}, at how the true system’s latent state evolves. On this reduced timescale, small fluctuations in the true data from a pure sine wave become more visible, due to measurement noise and the Koopman autoencoder’s slightly imperfect mapping from the full space to the latent space. The KAE EnKF’s 10-step ahead forecast is highly accurate, and produces stable, relatively smooth forecasts, indicating the filter’s system and measurement uncertainty matrices are well calibrated. Errors in the forecast tend to grow at the extremes of the data’s peaks and troughs. We now map the 10-step ahead latent forecasts back to the full state space, and subtract them from the true data, to view how their predictions of the pendulums position corresponds to its position in the true video. The subsequently generated image frames can be seen in Figure \ref{fig:pendulum_kae_enkf_val_forecast_subtracted}.

The full state forecasts are reasonably accurate, and similarly to the training data reconstructions, experience blurring of the pendulum’s image as its speed increases, reflecting the model’s reduced confidence in the pendulums exact position. To fully test the KAE EnKF's potential, we now introduce a non-stationary parameter to the system, and examine how effectively the KAE EnKF is able to adapt to the change, and how its tracking and forecasting performance is impacted.

\subsection{Tracking and forecasting unseen data from a non-stationary system}

\begin{figure}[htbp]
\begin{subfigure}[b]{\textwidth}
         \includegraphics[width=\textwidth]{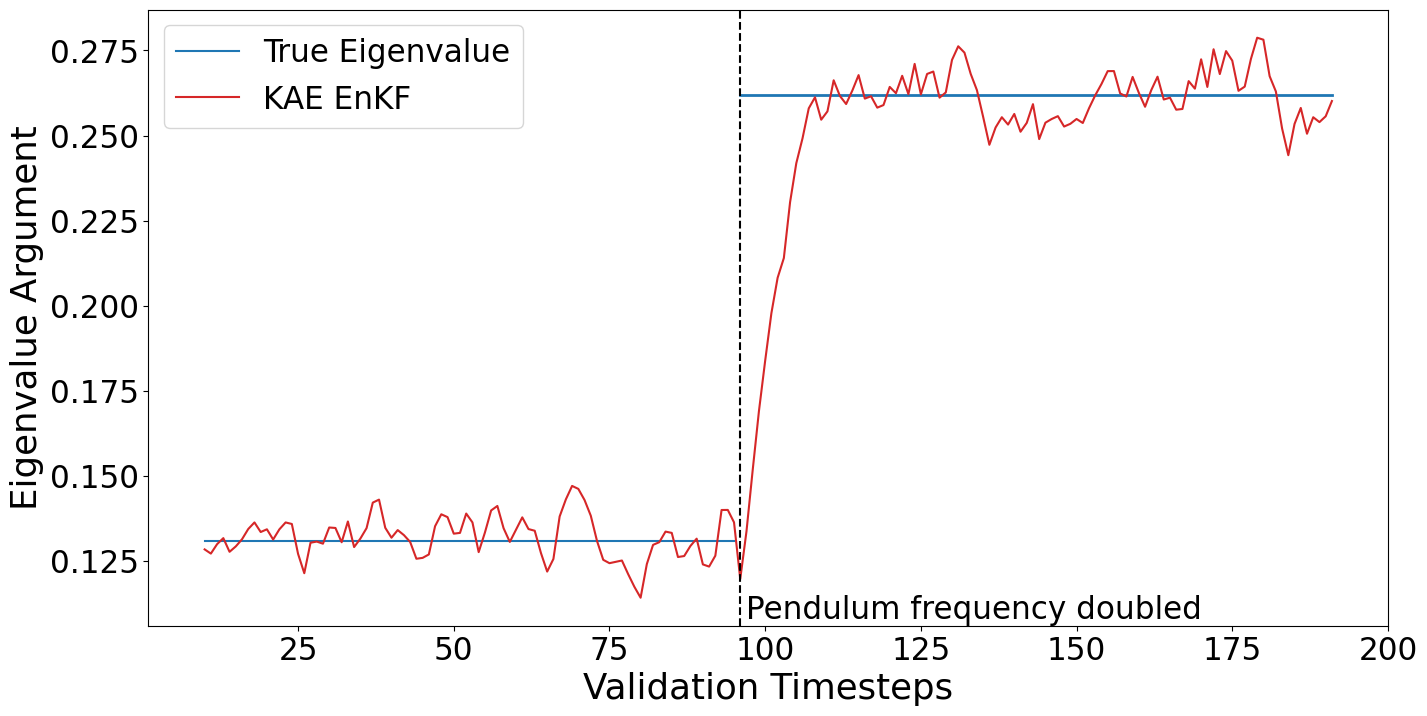}
     \end{subfigure}
                   \caption{\centering\label{fig:pendulum_kaeenkf_latent_freqdoubled_frequency}The system's true eigenvalue argument, and the KAE EnKF's eigenvalue argument estimate, over the augmented testing set, which includes a doubling of the system's frequency at $k=96$. The KAE EnKF quickly adapts to the change in the system's frequency, and maintains accurate estimates with small fluctuations at times when the system's parameter is stationary.}
\end{figure}

To augment the system to include a non-stationary parameter, we restrict the validation set to its first 96 time steps, so that the data divides exactly by the system's period. In doing so, we can stack multiple copies of the validation set together one after the other, without introducing a discontinuity in the system's state. The two copies of the validation set that we append to the original validation set also have every other frame removed, hence inducing a doubling of the system's eigenvalue argument to $2\pi/24$. The evolution of the system's eigenvalue argument over the course of the validation set can be seen in Figure \ref{fig:pendulum_kaeenkf_latent_freqdoubled_frequency}.

Figure \ref{fig:pendulum_kaeenkf_latent_freqdoubled_frequency} also shows how the KAE EnKF's eigenvalue argument estimate develops when tracking the non-stationary system. At the point where the system's frequency doubles, the KAE EnKF reacts almost immediately, rapidly increasing its eigenvalue argument estimate over the next 10 time steps until it reaches the system parameter's new value. When the system's eigenvalue argument is constant, the KAE EnKF's estimate fluctuates around the parameter’s true value, however does not converge to it, even as more time steps supporting this value of the parameter to be correct become available. This is caused by the choice to use constant covariance matrices over all data points in the filtering step. Reducing the filter’s system uncertainty makes the KAE EnKF’s parameter estimates more resistant to change, promoting more stable estimates of the eigenvalue’s argument during periods when they are constant, however inducing a further delay in adapting to their new values if they are non-stationary. As such, the filter's covariance matrices must be chosen to balance exploitation of consistent information, while maintaining the capacity for innovation in the event of new data, which contradicts the model's parameter value assumptions. One possible solution would be to condition the filter's system covariance matrix on the difference between the previous time steps forecast and measured states, also known as measurement innovation in the EnKF literature. This would lead to tighter parameter estimates as we become more confident the model of the system is accurate, while allowing for rapid alterations to the model's parameters in the event of a changing system parameter, as its forecasts and newly assimilated measurements begin to diverge.

\begin{figure}[htb]
\begin{subfigure}[b]{\textwidth}
         \includegraphics[width=\textwidth]{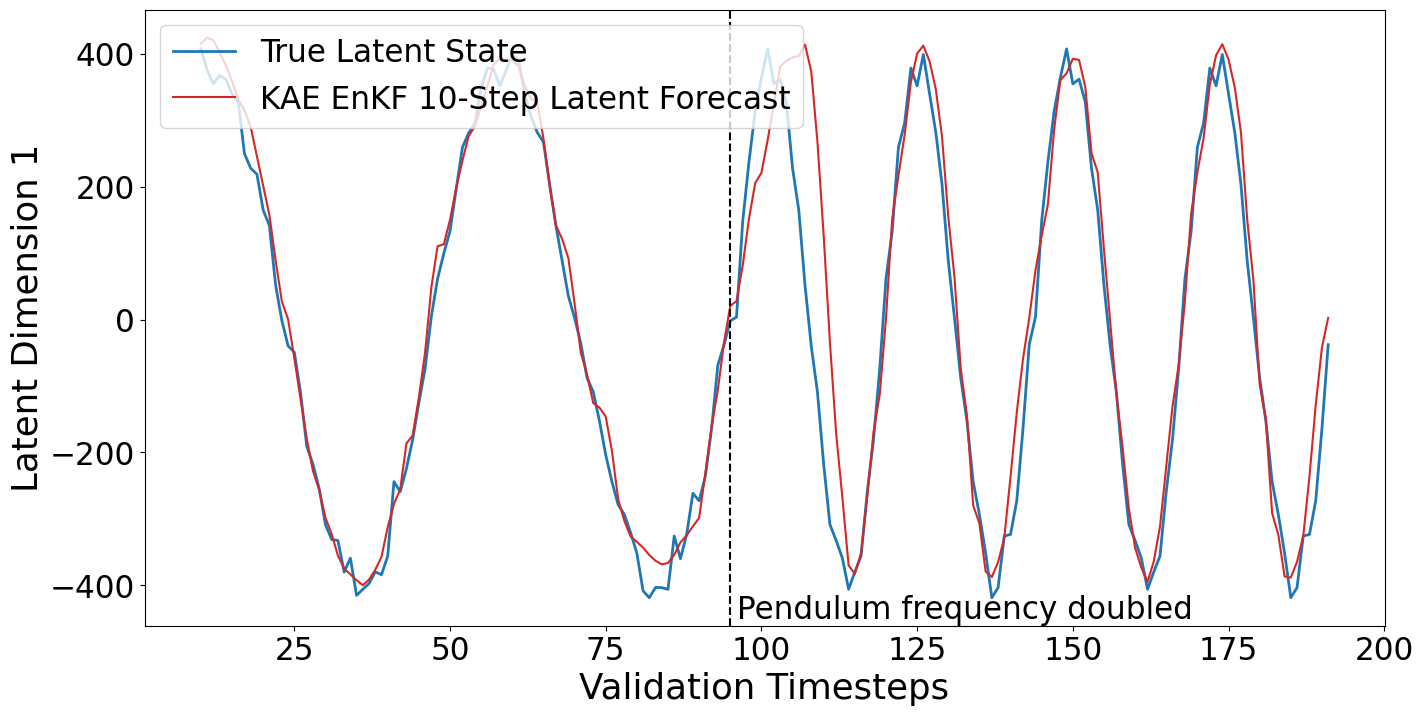}
     \end{subfigure}
     \caption{\centering\label{fig:pendulum_kaeenkf_latent_freqdoubled_forecast}The KAE EnKF’s 10-step ahead forecast, in the first dimension of the model’s
latent space, over the augmented testing set, that includes a doubling of the system's frequency at $k=96$. The KAE EnKF's forecast adapts quickly to the system's new frequency, with an additional delay the length of the forecast horizon relative to Figure \ref{fig:pendulum_kaeenkf_latent_freqdoubled_frequency}, before returning to producing accurate forecasts of the updated system.}
\end{figure}
Figure \ref{fig:pendulum_kaeenkf_latent_freqdoubled_forecast} demonstrates how the KAE EnKF's 10-step ahead latent state forecast reacts to the model's changing frequency estimate. The sine wave representing the system's true states undergoes an immediate frequency doubling at $k=96$, however the 10-step ahead KAE EnKF forecast does not re-synchronise with the system's true values until approximately 20 time steps later. This is to be expected, as Figure \ref{fig:pendulum_kaeenkf_latent_freqdoubled_frequency} shows the KAE EnKF takes around 10 time steps to adapt to the system's new eigenvalue argument, from which point forecasts generated using the corrected model would begin to return to their original accuracy levels. Hence, rectified forecasts would not present in the results until a time 10 steps plus the forecast horizon after the system's parameter change, which in this case is approximately 20 time steps later. After this point, the KAE EnKF's forecasts return to the similarly high levels of accuracy experienced before the system's frequency doubled.

\begin{figure}[htb]
\begin{subfigure}[b]{\textwidth}
         \includegraphics[width=\textwidth]{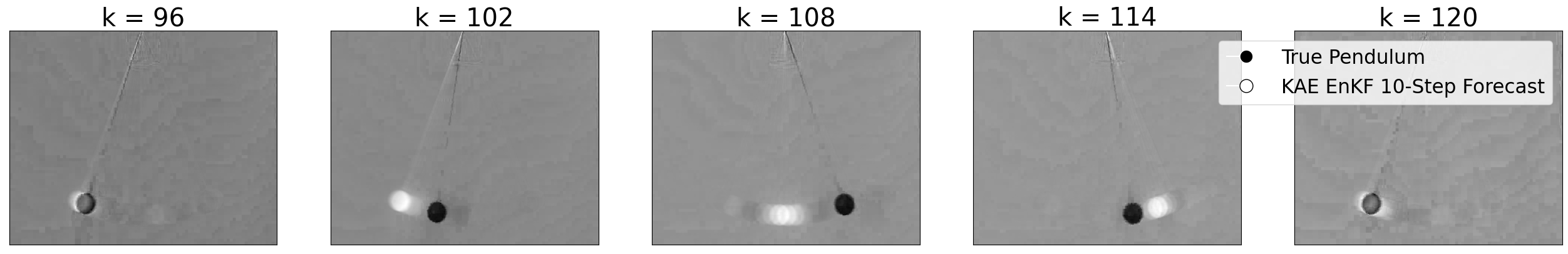}
     \end{subfigure}
     \caption{\centering\label{fig:period_doubled_pendulum_kae_enkf_val_forecast_subtracted}The difference between the true pendulum data, and the KAE EnKF's corresponding 10-step ahead full state forecasts, from the point in the augmented test set when the system's frequency doubles. The true pendulum is represented by a black circle, and the forecast, by virtue of being subtracted from the true data, a white circle. The KAE EnKF's forecast is initially accurate at $k=96$, then begins to lag behind the now twice as fast true system's pendulum for $k=102$ and $k=108$, before beginning to catch up at $k=114$, and finally re-synchronising with the true pendulums position by $k=120$, once the KAE EnKF has successfully assimilated the system's parameter change.}
\end{figure}

Figure \ref{fig:period_doubled_pendulum_kae_enkf_val_forecast_subtracted} shows the KAE EnKF's full state forecasts during the model's transition period between the system's original and doubled frequency. For $k=96$, we see the pendulum's position accurately forecast, marking the final time step of the system's original frequency. At $k=102$ and $k=108$, the pendulum can be seen accelerating away from the KAE EnKF’s forecasted position, due to its increased frequency. By $k=114$, the KAE EnKF has begun to increase its eigenvalue argument estimate, resulting in its full state forecast beginning to catch up with the true state’s value, although still trailing the true value slightly. As we reach $k=24$ the pendulum has completed a full oscillation at its new frequency, and the KAE EnKF has successfully extracted this frequency update from the assimilated data, yielding a forecast that predicts the pendulums position for the new system with a high level of accuracy.

In summary, when modelling the training pendulum video data, the KAE EnKF significantly outperformed EDMD, in both generating an effective mapping from the video frames to an appropriate latent space, as well as developing an accurate model of the system's latent dynamics. The KAE EnKF's high performance was maintained when forecasting previously unseen data in the testing set, and when the pendulum's frequency was artificially changed, the KAE EnKF was able to quickly and efficiently update its model to reflect the system's new parameter.

\section{Conclusion}\label{s:conclusion}

In conclusion, we described alterations that can be made to the Koopman autoencoder's training procedure, to encourage more efficient training that identifies globally optimal and stable eigenvalues for oscillatory data streams. We then defined a new algorithm, the KAE EnKF, that fuses Koopman autoencoders with ensemble Kalman filtering. The KAE EnKF efficiently maintains up-to-date state and Koopman eigenvalue estimates, for dynamical systems that are constantly generating new data.

When applied to a synthetic linear system with a time-varying eigenvalue, the KAE EnKF is able to track the system's eigenvalue only slightly more accurately than its linear counterpart, the DMDEnKF \cite{dmdenkf}. However, as the degree of nonlinearity in the system is increased, the KAE EnKF's estimates become vastly more stable and accurate than those produced by the DMDEnKF. The KAE EnKF's forecasting and eigenvalue tracking performance on this time-varying, nonlinear system is then compared against that of other existing iterative, nonlinear data-driven modelling techniques. The KAE EnKF's estimates are significantly more stable than estimates generated by the techniques designed to track non-stationary systems, and significantly less biased than the estimates produced by techniques designed to model stationary systems. By employing appropriate restrictions to the set of frequencies generated by the KAE EnKF when it is set to output multiple frequencies, these results also hold for a synthetic nonlinear, time-varying system with multiple underlying frequencies. By filtering the Koopman autoencoder's latent state as opposed to its full state, the KAE EnKF is quadratically more efficient in processing new data, while also producing more accurate state forecasts and parameter estimates of the synthetic system. This is caused by uncertainty being more effectively propagated in the system's dynamically relevant latent directions, when applied directly to the KAE EnKF's latent space.

The KAE EnKF was then applied to raw video footage of a pendulum in motion. The Koopman autoencoder successfully generates a low-dimensional latent state representation of the pendulum video, that captures sufficient information for accurate reconstruction of the pendulum's position, and in which the system's dynamics act linearly. Conversely, extended DMD when applied to this problem achieved significantly lower quality reconstructions of the pendulum's position, and was unable to produce an effective model of the system's latent dynamics. The KAE EnKF was able to efficiently produce short-term forecasts of the pendulum's future latent states, before decoding the forecasts back into video frames, and these forecasts in latent/pixel space were both highly accurate on previously unseen data. When external forcing was applied to double the pendulum's frequency, the KAE EnKF quickly adapted to the change, and identified the system's updated eigenvalue, within receiving only $\sim10$ frames that had been generated using the new parameter.

This work presents many exciting avenues for future research, as the KAE EnKF framework is highly generalizable by design. In recent times, Koopman operator based modelling techniques have been effectively utilized in a number of applications \cite{modern_koopman_theory}. Hence, the development of a framework in the KAE EnKF for robustly generating a data-driven model of a nonlinear dynamical system that adapts to new data as it becomes available, will likely have a significant use case in those fields where Koopman based methods have already seen success, potentially on systems too complex to be amenable to analysis by DMD alone. Additionally, the KAE EnKF provides model parameter and state uncertainty estimates, which are critical for making real world decisions based on the output of a mathematical model.

More generally, there exist systems that do not admit effective, finite Koopman representations, for example if the system contains multiple fixed points or more complicated attractors \cite{koopman_representation_problems}. In this case, the KAE EnKF framework could still be utilized, by modifying the autoencoder's constraint of linear latent dynamics, to instead enforce dynamics that are sparse \cite{autoencoder_sindy}, or physics-informed \cite{physics-informed_autoencoders}. In this way, the KAE EnKF's ability to produce uncertainty estimates, and process streaming data to improve its model on the fly, could be taken advantage of by alternative autoencoder based models of dynamical systems.

Looking more generally still, filtering the latent state of a system as opposed to its full state has huge potential in many situations where data assimilation techniques are applied \cite{latent_data_assim, latent-kalmannet}. Filtering over a reduced state size will often significantly improve the filter's efficiency, particularly as many data assimilation applications are performed over high-dimensional state spaces. As shown empirically in this paper, when applied to synthetic data, latent state filtering can also improve the filter's tracking and forecasting performance, by ensuring uncertainty is robustly represented in the most dynamically relevant directions of the system's state. The potential to intelligently assimilate data from high-dimensional systems, with improved tracking and forecasting performance, and at a fraction of the previous computational cost, presents an exciting opportunity for fields such as numerical weather prediction and other geophysical sciences.

\section*{Code availability}
Codes used to produce the results in this paper are available at: 

\href{https://github.com/falconical/KAE-EnKF}{https://github.com/falconical/KAE-EnKF}.

\section*{Acknowledgements}
This work was supported by the UKRI, whose Doctoral Training Partnership Studentship helped fund Stephen Falconers PhD. We would also like to thank Stefan Klus, for providing the pendulum video data used in this work, as well as useful insights. We also thank Nadia Smith and Spencer Thomas from the National Physical Laboratory for their valuable discussions.


\section*{Appendix}

\section{KAE Neural Network Implementation Details}

\subsection{General Implementation details}

The KAE networks were implemented in code using the Python machine learning library Pytorch \cite{pytorch}, and Adam optimizer \cite{adam_opt} was utilized during the stochastic gradient descent sections of the KAE's training. The primary hyperparameters that require determining before training the KAE are weightings for each term in the loss function from equation \eqref{eqn:kae_enkf_loss2}, and the stochastic gradient descent algorithm's learning rate.

When applied over the synthetic datasets in section 4.5, the KAE was found to be highly robust with respect to the hyperparameters chosen. Hence, a loss function weighting of $[1, 1, 1 ,1, 0.01]$ was selected, representing an equal weighting of the reconstruction, latent dynamics, full state dynamics and mode stability loss terms, with a reduced weighting of the networks' regularization loss. This configuration was by no means optimal, however was selected  as it produced a strong performing KAE, while remaining simple to represent. In section 4.6, when applying the KAE to real pendulum video data, hyperparameter selection proved significantly more important. Applying the simple loss function weightings from section 4.5 yielded poor performance, hence the Python hyperparameter optimization library Hyperopt \cite{hyperopt_library} was utilized, to intelligently select loss function weightings that yielded a KAE with minimal forecasting errors.

For stochastic gradient descent algorithms, larger learning rates typically induce faster training, however too high a learning rate can hinder training by consistently overshooting nearby minima. For this reason, a general rule of thumb is to select the maximal learning rate for which training loss still decreases, and this logic was applied when choosing a learning rate in all our applications.

\subsection{Network Architectures}
In sections 4 and 5 the KAE EnKF is applied to a range of synthetic and real datasets. The associated KAE's encoder/decoder networks are built by stacking fully connected layers, each followed by ReLU activation functions to introduce nonlinearity into the network. The only exceptions are the encoder and decoder's final layer, which consist of a fully connected layer only. This ensures elements of the model's latent/full state can take any real values, as opposed to only non-negative values as output by the ReLU functions. The KAE's decoder architecture consists of the same size and number of layers as used in the encoder, however stacked in the reverse order, as is standard when designing an autoencoder's structure.

As a starting point for each application, a simple encoder network was used, consisting of an input layer, a single hidden layer, and an output latent layer. The first dimension of the input layer was fixed by the problem's state dimension, and the latent layer's output dimension was determined by the number of frequencies being used to model the system. Hence, the only free network architecture variable was the single hidden layer's dimensions, which we initialized with an input/output size of 10. Abiding by the principle that models should be as simple as possible but no simpler, if the network performed sufficiently well in tracking the system's parameters and generating future state forecasts with this network architecture, the hidden layer's dimension were kept at 10. If the network performed poorly, the hidden layer's input and then output size was doubled to 20 neurons. If performance remained unsatisfactory after these increases to the network's width, the intention was to then increase the network's depth by adding another hidden layer, however this requirement was not met for any of the applications tested. Table \ref{table:kae_nn_arch} details the size of each layer in the KAE encoder's architecture, used in each of the KAE EnKF synthetic/real experiments.

\begin{table}[h]
 \begin{tabular}{||c c c c||} 
 \hline
 Experiment & \# Input Layer & \# Hidden Layer & \# Latent Layer \\ [0.5ex] 
 \hline\hline
     4.2 KAE EnKF vs Full State KAE EnKF& 100 x 10 & 10 x 10 & 10 x 2 \\
 \hline
   4.3.1 KAE EnKF vs EDMD - Single Frequency& 100 x 10 & 10 x 10 & 10 x 2 \\
 \hline
   4.3.2 KAE EnKF vs EDMD - Multiple Frequencies& 100 x 10 & 10 x 10 & 10 x 6 \\
 \hline
   5 Pendulum Video& 1244160 x 20 & 20 x 10 & 10 x 2 \\
 \hline 

\end{tabular}
\caption{\centering The size of each layer of the KAE's encoder, for each experiment conducted. The initial input layer's size is determined by the full state dimension, and the final latent layer's size, (after dividing by 2), dictates the number of latent frequencies in the model. The hidden layers size increase for systems with a greater degree of nonlinearity, for example when modelling the raw pendulum data in experiment 5.}
\label{table:kae_nn_arch}
\end{table}

For all synthetic experiments, the datasets contained a similar level of nonlinearity, hence requiring a hidden layer with only a width of 10 neurons. Experiment 4.3.2 modelled data generated by a system with 3 latent frequencies, meaning the output size of the encoder's latent layer needed tripling to allow for this additional frequency information to be passed to the KAE's Koopman approximator $\mathbf{K}_\lambda$. Each frame of the raw pendulum video in experiment 5 consisted of over 400,000 pixels, and 3 frames had to be passed to the network to ensure the full state also contained information about the pendulum's velocity and acceleration. This required an equally large initial input layer to read each of the state's pixels. The network then quickly contracted, with a hidden layer of only width 20 required to propagate the full state's dynamically relevant information into a single frequency latent state.

\bibliographystyle{siam}
\bibliography{mybibfile}

\end{document}